\documentclass[11pt]{amsart}

\usepackage[ngerman,british]{babel} 
\usepackage[utf8]{inputenc}
\usepackage[T1]{fontenc}
\usepackage{amsthm}
\usepackage{amssymb}
\usepackage{mathtools} 
\usepackage{amsmath}
\usepackage[outline]{contour}
\usepackage[mathscr]{euscript}
\usepackage[all,2cell,cmtip]{xy}
\usepackage{csquotes}
\usepackage[usenames,dvipsnames]{xcolor}
\usepackage[colorlinks=true,allcolors=blue,linktocpage]{hyperref}
\usepackage{microtype}

\usepackage{tikz-cd}
\usepackage{amssymb}
\usetikzlibrary{calc}
\usetikzlibrary{decorations.pathmorphing}

\usetikzlibrary{babel} 

\tikzset{curve/.style={settings={#1},to path={(\tikztostart)
    .. controls ($(\tikztostart)!\pv{pos}!(\tikztotarget)!\pv{height}!270:(\tikztotarget)$)
    and ($(\tikztostart)!1-\pv{pos}!(\tikztotarget)!\pv{height}!270:(\tikztotarget)$)
    .. (\tikztotarget)\tikztonodes}},
    settings/.code={\tikzset{quiver/.cd,#1}
        \def\pv##1{\pgfkeysvalueof{/tikz/quiver/##1}}},
    quiver/.cd,pos/.initial=0.35,height/.initial=0}

\tikzset{tail reversed/.code={\pgfsetarrowsstart{tikzcd to}}}
\tikzset{2tail/.code={\pgfsetarrowsstart{Implies[reversed]}}}
\tikzset{2tail reversed/.code={\pgfsetarrowsstart{Implies}}}
\tikzset{no body/.style={/tikz/dash pattern=on 0 off 1mm}}

\usepackage[hyperref,
maxbibnames=99,
backend=biber,
style=alphabetic-verb,
citestyle=alphabetic-verb,
giveninits=true
]{biblatex}
\DefineBibliographyStrings{english}{
  page             = {},
  pages            = {},
} 
\renewbibmacro{in:}{} 
\usepackage{graphicx}
\usepackage{cleveref}   
\usepackage{underscore} 

\usepackage{import}
\usepackage{color} 
\usepackage{xifthen}

\usepackage{rotating}

\newcommand{\mcal}[1]{\mathcal{#1}}

\newcommand{\mscc}[1]{\mathbf{#1}}
\newcommand{\mrm}[1]{\mathrm{#1}}
\newcommand{\mbf}[1]{\mathbf{#1}}
\newcommand{\mbb}[1]{\mathbb{#1}}
\newcommand{\mfrak}[1]{\mathfrak{#1}}
\newcommand{\down}[3][ ]{(#2 \downarrow #3)^{#1}}
\newcommand{\paren}[1]{\left( #1 \right)}

\newcommand{\id}{{\mathrm{id}}}
\newcommand{\colim}{\mathrm{colim}}
\newcommand{\Hom}{\mathrm{Hom}}

\newcommand{\Bred}{\mathcal{B}\mbox{-}\mathrm{red}}

\newcommand{\Vinj}{{\mscc{V}^{\mathrm{inj}}}}
\newcommand{\Vhook}{\mscc{V}^{\hookrightarrow}}
\newcommand{\Vsim}{{\mscc{V}^{\simeq}}}

\newcommand{\opp}{{\mathrm{op}}}
\newcommand{\Opp}{{\mathrm{Op}}}

\newcommand{\Cinfty}{\mscc{C}\mathrm{at}_{\infty}}

\newcommand{\cdon}{\[\begin{tikzcd}}
\newcommand{\cdoff}{\end{tikzcd}\]}
\newcommand{\ncdon}{\begin{equation}\begin{tikzcd}}
\newcommand{\ncdoff}{\end{tikzcd}\end{equation}}
\newcommand{\OO}{{\mathrm{O}}}

\newcommand{\Gr}{{\mathrm{Gr}}}

\newcommand{\Sing}{{\mathrm{Sing}}}
\newcommand{\Bbox}{\mathcal{B}^{\oplus}}

\newcommand{\Nerve}{{\mathrm{N}}}
\newcommand{\Ndelta}{\mathrm{N}^{\mbf{\Delta}}}
\newcommand{\Nhc}{\mathrm{N}^{\mathrm{hc}}}
\newcommand{\sSet}{\mathrm{sSet}}
\newcommand{\Cdelta}{{\mathcal{C}\mathrm{at}_{\mbf{\Delta}}}}
\newcommand{\pSh}{\mathrm{pSh}}
\newcommand{\EEx}{\mathbf{EX}}

\newcommand{\Phat}{\mscc{P}^{\Delta}}

\newcommand{\Fun}{\mathrm{Fun}}
\newcommand{\Path}{\mathrm{Path}}
\newcommand{\Tangent}{\mathrm{T}}
\newcommand{\Normal}{\mathrm{N}}
\newcommand{\vfr}{\mathrm{vfr}}

\newcommand{\evv}{\mathrm{ev}}

\newcommand{\btu}{\bigtriangleup}
\makeatletter
\newcommand{\btd}{\mathord{\mathpalette\raise@half\bigtriangledown}}
\newcommand\raise@half[2]{%
  \raisebox{.9\depth}{$\m@th#1#2$}%
}
\makeatother

\newcommand{\cat}[2]{\mscc{#1}\mathrm{#2}}

\numberwithin{equation}{section}

\theoremstyle{definition}
\newtheorem{definition}[equation]{Definition}
\newtheorem*{definition*}{Definition}
\newtheorem{notation}[equation]{Notation}
\newtheorem*{notation*}{Notation}
\newtheorem{construction}[equation]{Construction}
\newtheorem*{construction*}{Construction}

\theoremstyle{remark}
\newtheorem{remark}[equation]{Remark}
\newtheorem*{remark*}{Remark}
\newtheorem{example}[equation]{Example}
\newtheorem*{example*}{Example}

\newtheorem*{warning*}{Warning}

\newtheorem*{convention*}{Convention}

\theoremstyle{plain}
\newtheorem{proposition}[equation]{Proposition}
\newtheorem*{proposition*}{Proposition}
\newtheorem{theorem}[equation]{Theorem}
\newtheorem*{theorem*}{Theorem}
\newtheorem{definition-theorem}[equation]{Definition-Theorem}
\newtheorem*{definition-theorem*}{Definition-Theorem}

\newtheorem*{postulate*}{Postulate}
\newtheorem{lemma}[equation]{Lemma}
\newtheorem*{lemma*}{Lemma}

\newtheorem*{observation*}{Observation}
\newtheorem{corollary}[equation]{Corollary}
\newtheorem*{corollary*}{Corollary}

\newtheorem*{question*}{Question}

\title{{The stratified Grassmannian and its depth-one subcategories}}

\author{\"Od\"ul Tet\.{i}k}
\thanks{This research was partially supported by the NCCR SwissMAP, funded by the Swiss National Science Foundation (SNF) and by the Austrian Science Fund (FWF) through SNF Grants no.\ 200020\_192080 and 200021\_227719 and FWF Project no.\ P 37046.}
\address{Institut für Mathematik, Universität Zürich, Winterthurerstrasse 190, 8057 Zürich, Switzerland}
\curraddr{University of Vienna, Faculty of Physics, Mathematical Physics Group, Boltzmanngasse 5, 1090 Vienna, Austria}
\email{oeduel.tetik@univie.ac.at}

\date{}

\begin{document}

\maketitle

\begin{abstract}
    We introduce a tangential theory for linked smooth manifolds of depth $1$, i.e., for spans $\mathfrak{S}=(M\overset{\pi}{\twoheadleftarrow} L\overset{\iota}{\hookrightarrow}N)$ of smooth manifolds where $\pi$ is a fibre bundle and $\iota$ is a closed embedding. The tangent classifier of $\mathfrak{S}$ is given as a topological span map $\mathfrak{S}\to B\mathrm{O}(n,m)$ where $B\mathrm{O}(n,m)=(B\mathrm{O}(n)\twoheadleftarrow B\mathrm{O}(n)\times B\mathrm{O}(m)\hookrightarrow B\mathrm{O}(n+m))$. We show that this recovers and generalises the tangential theory introduced by Ayala, Francis and Rozenblyum for conically smooth stratified spaces by constructing fully faithful functors $\mathbf{EX}(B\mathrm{O}(n,m))\hookrightarrow\mathbf{V}^{\hookrightarrow}$ of quasi-categories, where $\mathbf{EX}$, introduced in a prequel, takes the exit path quasi-category of the span, and $\mathbf{V}^{\hookrightarrow}$ is a quasi-category model of the infinite stratified Grassmannian of AFR. This result has analogues for other classical structure groups and for Stiefel manifolds. We thus reduce the classification of conically smooth bundles over depth-$1$ posets to that of ordinary bundles on linked smooth manifolds.
\end{abstract}

\tableofcontents

\tikzset{Rightarrow/.style={double equal sign distance,>={Implies},->},
triple/.style={-,preaction={draw,Rightarrow}},
quadruple/.style={preaction={draw,Rightarrow,shorten >=0pt},shorten >=1pt,-,double,double
distance=0.2pt}}

\section{Introduction}

In \cite{ayala2020factorization}, Ayala, Francis and Rozenblyum (AFR) defined the tangent bundle of a conically smooth stratified space (CSS) $M$ in the form of a functor 
\[
    \Tangent M\colon\cat{E}{xit}(M)\to\Vinj,
\]
where $\cat{E}{xit}(M)$, the \emph{exit path $\infty$-category} of $M$, is the stratified homotopy type of $M$ which reduces to the singular chains complex of $M$ when it is trivially stratified. The target $\Vinj$, which we will for the moment call the \emph{stratified Grassmannian}, is the stratified analogue of the classifying space $B\OO(n)$. In fact, $n$ varies within $\Vinj$, reflecting the different dimensions the strata of $M$ may have. 
A morphism in $\cat{E}{xit}(M)$ is, roughly speaking, either a path which stays within a single stratum or a path that starts in a stratum and immediately exists into a higher adjacent stratum. As soon as $M$ has multiple adjacent strata, $\cat{E}{xit}(M)$ is not an $\infty$-groupoid but merely an $\infty$-category.
The target may be roughly thought of as the $\infty$-category whose objects are finite dimensional real vector spaces and whose morphisms are linear injections. On this view, $B\OO(n)$ may be interpreted as the sub-$\infty$-groupoid of $\Vinj$ consisting of only $n$-dimensional (real) vector spaces (inside $\mbb{R}^\infty$), through which $\Tangent M$ will factor if and only if $M$ consists of only $n$-dimensional strata. Indeed, the maximal sub-$\infty$-groupoid of $\Vinj$ is $\coprod_{n\geq0}B\OO(n)$.

With this in place, AFR then (re-)introduced\footnote{See also Ayala--Francis--Tanaka \cite{ayala2017local}.} a rich theory of tangential structures on CSSs by considering $\infty$-categories over $\Vinj$, and the associated lifting problems analogous to those in the classical theory.

This stratified tangent classifier is constructed in an indirect manner, making it difficult to unpack and work with in practice. Given the enormous current interest in stratified space theory within mathematical physics generated by the study of defects, symmetries, and functorial field theory, a simpler, hands-on approach to stratified tangent classifiers, and to stratified tangential structures generally, is desirable. The goal of the present paper is to provide such an approach, which in fact formalises a heuristic discussion by AFR themselves. We achieve this in depth $1$, leaving higher depth to future work.\footnote{By \emph{depth} we refer to the depth of the stratifying poset.}

\subsection{The linked tangent bundle} Let $n\geq0$ and $m\geq1$ be natural numbers, and consider a smooth $(n+m)$-manifold $M$ without boundary, stratified by a closed $n$-dimensional submanifold $\Sigma\subset M$. The non-invertible morphisms in $\cat{E}{xit}(M)$ are those paths $\gamma\colon[0,1]\to M$ which satisfy $\gamma(0)\in\Sigma$ and $\gamma(t)\in M\smallsetminus \Sigma$ for $t>0$. In fact, the space of such paths is equivalent, by a result of \cite{ayala2017stratified}, to the \emph{link} of $\Sigma$ in $M$, which in this case may be taken to be $L=\mbb{S}(N\Sigma)$, the sphere bundle of the normal bundle of $\Sigma$. The two strata and the link are summarised in the span
\cdon 
    & \mbb{S}(N\Sigma)\ar[dr,"{\iota}",hook]\ar[dl,two heads, "{\pi}"'] & \\
    \Sigma && M\smallsetminus\Sigma
\cdoff
where $\iota$ is the closed embedding induced by the choice of a tubular neighbourhood, and $\pi$ is the bundle projection. Now, Section 2.1.4 of op.\ cit.\ explains that the tangent classifier may be described heuristically by means of this span. The fundamental $\infty$-groupoids $\Pi_\infty(X)=\mrm{Sing}_\bullet(X)$ of $X=\Sigma, M\smallsetminus\Sigma$ embed into $\cat{E}{xit}(M)$, and the restrictions $\Tangent M|_{X}\colon \Pi_\infty(X)\hookrightarrow\cat{E}{xit}(M)\overset{\Tangent M}{\to}\Vinj$ coincide with the ordinary tangent classifiers.
Now note the following phenomenon: there are two induced bundles, $\Tangent^1=\pi^*\Tangent \Sigma$ and $\Tangent^2=\iota^*\Tangent(M\smallsetminus \Sigma)$, over $L$. Up to a contractible choice, there exists a bundle embedding $\Tangent^1\hookrightarrow \Tangent^2$ (detailed in \Cref{constr:tangent bundle of linked manifold} below). Thus, the assignment 
\[
    L\ni\ell\mapsto (\Tangent^1_\ell\hookrightarrow\Tangent^2_\ell)
\]
defines a map from $L$ to non-invertible morphisms in $\Vinj$, describing $\Tangent M$ on non-invertible morphisms.

\subsection{Main results and applications}
The main idea of the present paper is to construct a theory where the above becomes a definition. This is easily achieved: Define a \emph{linked (smooth) manifold} to be a span $\mfrak{S}=(M\xleftarrow{\pi}L\xrightarrow{\iota}N)$ where $M$, $L$, $N$ are smooth manifolds, $\pi$ is a fibre bundle and $\iota$ is a closed embedding. Every CSS $X$ stratified over $\{0<1\}$ gives rise to linked manifold $\mfrak{S}_X$. It is straightforward to construct a tangent classifier $\Tangent\mfrak{S}$ of $\mfrak{S}$ as a topological span map of type
\begin{equation*}
    \begin{tikzcd}[column sep=tiny]
        & L \ar[ddl,two heads,"{\pi}"']\ar[dr,hook,"{\iota}"] \ar[rrr,"{\pi^*\Tangent M\times \Normal}"] &&& B\mathrm{O}(n)\times B\mathrm{O}(m) \ar[dr,hook,"{\oplus}"] \\
        && N \ar[rrr,"{\Tangent N}",pos=0.25] &&& B\mathrm{O}(n+m)\\
        M \ar[rrr,"{\Tangent M}"] &&& B\mathrm{O}(n) \ar[from=uur,two heads,crossing over,"{\mathrm{pr}}"',pos=0.27] 
    \end{tikzcd}
\end{equation*}
for appropriate $n$ and $m$, where $\Normal\colon L\to B\OO(m)$  classifies the normal bundle of the analogue of the bundle embedding $\Tangent^1\hookrightarrow\Tangent^2$ discussed above. In \Cref{T6HNI7D} we observe that classifiers can be chosen that make the diagram commute. Let us write $B\OO(n,m)$ for the displayed span of Grassmannians, which we call the (infinite) \emph{$(n,m)$-Grassmannian}. The question immediately arises as to how exactly this is related to the theory of AFR. The main contribution of this paper is the following:
\begin{theorem*}[{\ref{YUXW1VO}}]
    There exists a fully faithful functor 
    \[
        \mbf{U}\colon \EEx(B\OO(n,m))\hookrightarrow\Vhook
    \]
    of $\infty$-categories.
\end{theorem*}

We call $\mbf{U}$ the \emph{unpacking map}. In other words, $\mbf{U}$ gives an equivalence 
\[
    \EEx(B\OO(n,m))\simeq\Vhook|_{n,n+m}
\]
where the right-hand side is the full sub-$\infty$-category of $\Vhook$ (or, equivalently, of $\Vinj$) generated by $n$- and $(n+m)$-dimensional vector spaces. We call $\Vhook|_{n,n+m}$ a \emph{(full) depth-$1$ sub-$\infty$-category} of $\Vhook$. 

The exit path $\infty$-category construction $\EEx(-)$ was introduced in the prequel \cite{tetik2023linked} for spans of greater generality than linked manifolds, including the $(n,m)$-Grassmannian. We recall this construction briefly in \Cref{6CW6YDL}. The target $\Vhook$ of $\mbf{U}$ is a quasi-categorical version of $\Vinj$ which differs from it in no important respect but only up to some redundancy that is the cost of a strictification, modulo which it was suggested already in AFR \cite[Remark 2.7]{ayala2020factorization} in its incarnation as a Segal space. Neither construction refers to any type of stratified space. Consequently, no stratified space theory is required or used in the present paper except as motivation.

\Cref{YUXW1VO}, combined with other results, entails that the topological span map $\Tangent \mfrak{S}$ is enough to understand the tangent classifiers of depth-$1$ CSSs of depth $1$. Moreover, the theory of tangential structures can be pulled back to the linked setting along $\mbf{U}$. We will now note the main corollaries using the following result from the prequel:

\begin{theorem*}[{\cite{tetik2023linked}}]
    There exists a fully faithful functor \(\EEx\colon\mscc{LS}\hookrightarrow\Cinfty\) from the $\infty$-category of linked spaces to the $\infty$-category of all $\infty$-categories. Moreover, for $X$ a CSS of depth $1$ and $\mfrak{S}_X$ the associated linked manifold, there exists an equivalence 
    \(\mbf{R}\colon \EEx(\mfrak{S}_X)\xrightarrow{\sim}\cat{E}{xit}(X)\) of $\infty$-categories.
\end{theorem*}

For simplicity, we will state the corollaries for CSSs $X$ over $\{0<1\}$. Let $\mfrak{S}$ denote the associated linked smooth manifold.

\begin{corollary}
    There exists a weak equivalence $\Hom_{\Cinfty}(\cat{E}{xit}(X),\Vinj|_{n,n+m})\simeq \Hom_{\mscc{LS}}(\mfrak{S},B\OO(n,m))$.
\end{corollary}

\begin{corollary}\label{7UKFW0A}
    Let the dimensions of the strata of $X$ be $n$ and $n+m$. The tangent classifier $\Tangent X$ of $X$ factors as in the diagram 
    \[
        \begin{tikzcd}
            \cat{E}{xit}(X)\ar[r,"{\Tangent X}"] & \Vhook|_{n,n+m}\\
            \EEx(\mfrak{S})\ar[r,"{\Tangent \mfrak{S}}"]\ar[u,"{\sim}","{\mbf{R}}"'] & \EEx(B\OO(n,m))\ar[u,"{\sim}","{\mbf{U}}"']
        \end{tikzcd}
    \]
    where both vertical maps are equivalences.
\end{corollary}

Recall from \cite{ayala2020factorization} that conically smooth vector bundles over a CSS $X$ are classified by functors $\cat{E}{xit}(X)\to\Vhook$.

\begin{corollary}\label{XTGZI27}
    The moduli space
    \[
        \mscc{E}_{n,m}(X)\subset\Fun^{\simeq}(\cat{E}{xit}(X),\Vinj)=\Hom_{\Cinfty}(\cat{E}{xit}(X),\Vinj)
    \] 
    of conically smooth vector bundles $E$ on $X$ satisfying $\mrm{rk}({E}|_M)=n$ and $\mrm{rk}({E}|_N)=n+m$ is weakly equivalent to the space of topological span maps $\mfrak{S}\to B\OO(n,m)$:
    \[
        \mscc{E}_{n,m}(X)\simeq \Hom_{\cat{T}{op}}(\mfrak{S},B\OO(n,m)).
    \]
    Similarly, the full moduli space admits a weak equivalence 
    \[
        \mscc{E}(X)=\Hom_{\Cinfty}(\cat{E}{xit}(X),\Vinj)\simeq\coprod_{n,m\geq0}\Hom_{\cat{T}{op}}(\mfrak{S},B\OO(n,m)).
    \]
\end{corollary}

Let us finally note an analogous corollary for conically smooth principal bundles in depth $1$. Let $\{\xi_n\colon G(n)\to \OO(n)\}_{n\geq0}$ be a good multiplicative system of topological groups over the orthogonal groups, i.e., there exist operations $\boxplus\colon G(n)\times G(m)\to G(n+m)$ which are cofibrations covering the direct sum $\oplus\colon \OO(n)\times \OO(m)\to \OO(n+m)$. In particular $BG(n,m)=(BG(n)\leftarrow BG(n)\times BG(m)\to BG(n+m))$ is a linked space over $B\OO(n,m)$. Classical structure groups provide examples, as do the Stiefel manifolds, bypassing the group level. The classifying $\infty$-category $\mscc{G}^{\hookrightarrow}$ of stratified $G$-bundles can be constructed analogously to $\Vhook$, and there exists a natural functor 
\[
\xi\colon\mscc{G}^{\hookrightarrow}\to \Vhook.
\]
The proof of \Cref{YUXW1VO} applies mutatis mutandis to give the following:
\begin{corollary}\label{QOSZIFS}
    There exists a fully faithful functor 
    \(
        \mbf{U}\colon\EEx(BG(n,m))\to\mscc{G}^{\hookrightarrow}
    \)
    of $\infty$-categories.
\end{corollary}

While conically smooth principal bundles have not been officially defined in the literature, let us define them to be classified by functors $\cat{E}{xit}\to\mscc{G}^{\hookrightarrow}$. A \emph{$G$-structure} on a conically smooth vector bundle can then be defined using $\xi$ as a homotopy lift of $\infty$-functors.

\begin{corollary}\label{IPWVX24}
    Let $G$ be as above. Then the moduli space $\mscc{G}(X)=\Hom_{\Cinfty}(\cat{E}{xit},\mscc{G}^{\hookrightarrow})$ of conically smooth $G$-bundles admits a weak equivalence 
    \[
        \mscc{G}(X)\simeq\coprod_{n,m\geq0}\Hom_{\cat{T}{op}}(\mfrak{S},BG(n,m)).
    \]
    Moreover, the space $G$-structures on a conically smooth vector bundle $E\colon\cat{E}{xit}(X)\to\Vinj$ satisfying $\mrm{rk}({E}|_M)=n$ and $\mrm{rk}({E}|_N)=n+m$ is weakly equivalent to the space of homotopy lifts
    \[
    \begin{tikzcd}
            & BG(n,m)\ar[d,"{\xi}"]\\
        \mfrak{S}\ar[ur,dotted]\ar[r,"{\mfrak{E}}"'] & B\OO(n,m)
    \end{tikzcd}
    \]
    of topological span maps, where $\mfrak{E}$ is a span map representative of $E$ along the weak equivalence $\Hom_{\Cinfty}(\cat{E}{xit}(X),\Vinj|_{n,n+m})\simeq \Hom_{\mscc{LS}}(\mfrak{S},B\OO(n,m))$.
\end{corollary}

\subsection{A synopsis of the unpacking map}

While both $\EEx(B\OO(n,m))$ and $\Vhook|_{n,m}$ are about $n$- and $(n+m)$-dimensional vector spaces and non-invertible morphisms between them, their constructions are very different. Say $V$, $K$ are vector spaces, and $\dim V=n$ and $\dim K=n+m$. A morphism $V\to K$ in both $\EEx(B\OO(n,m))$ and $\Vhook$ consists in the choice of an $m$-dimensional vector space $W$ together with a path 
\[
    W\oplus V\to K\quad  \text{within}\quad  B\OO(n+m).
\] 
However, higher-dimensional morphisms in the two $\infty$-categories are quite different, which is why the construction of $\mbf{U}$ is rather technical, and a better understanding (or version) of it remains desirable. Namely, the higher morphisms in $\EEx(B\OO(n,m))$ can be expressed succinctly in a similar manner, whereas $\Vhook$ is constructed as an under-$\infty$-category of an $\infty$-category defined as the homotopy-coherent nerve of a topological category (recalled in \Cref{app:simplicial nerves}). The most surprising aspect of the construction of $\mbf{U}$ presented here is that it seems to depend non-trivially on the classical adjunction between geometric realisation and the singular chains functor, rendering the construction not purely combinatorial (cf.\ \Cref{EEYT903}, specifically the proofs of all results after \Cref{HZBV0LB}).

To illustrate the construction of $\mbf{U}$ on an example, let us consider $2$-morphsims. A $2$-morphism of $\EEx(B\OO(n,m))$ of so-called exit index $2$ is determined by a $2$-simplex of $B\OO(n+m)$ of type 
\ncdon
\label{9ASEI61}
        & K\\
    W\oplus V \ar[r,"{\gamma_W\oplus\gamma_V}"']\ar[ur,"{\gamma_\oplus}"] & W'\oplus V'\ar[u,"{\gamma_{\oplus'}}"']
\ncdoff
where the bottom edge comes from $B\OO(n)\times B\OO(m)$. By construction this is a $2$-morphism in $\EEx(B\OO(n,m))$ of type 
\cdon[sep=small]
    &K \\
    V\ar[ur]\ar[r] & V'\ar[u]
\cdoff
with the vertical edges non-invertible. On the other hand we have $\Vhook=\ast/\Bbox\OO$, the under-$\infty$-category of the $\infty$-category $\Bbox\OO$ given as the homotopy-coherent nerve of the topological category $B^{\oplus}\OO$ whose sole object is $\ast$ and whose topological endomorphism monoid $B\OO^{\infty}_{\amalg}$ is given essentialy by $\coprod_{n\geq0}B\OO(n)$ together with the direct-sum operation $\oplus$. Thus $2$-morphisms of $\Vhook$ correspond to $3$-morphisms of $\Bbox\OO$, and $\mbf{U}$ sends the above $2$-morphism to a $3$-morphism of type
\cdon
    &2\ar[dr,"{W'}"{description}]\\
1\ar[ur,"{0}"description]   &   &3\\
    &0\ar[ul,"{V}"description]\ar[uu,"{V'}"{description,pos=0.2}]\ar[ur,"{K}"{description}]
\ar[from=ul,to=ur,crossing over,"{W}"{description,pos=0.7}]
\cdoff
As we will see, specifying such a $3$-morphism requires a path $W'\to W$ in $B\OO(m)$ in the `wrong direction,' together with $2$-simplices of type
\ncdon[sep=small] \label{LL2SEBY}
        & W'\oplus V'\ar[dr]\\
        W'\oplus V\ar[ur]\ar[rr] && K
\ncdoff 
and 
\ncdon[sep=small] \label{I2XVTKD}
    & W\oplus V\ar[dr]\\
    W'\oplus V\ar[ur]\ar[rr] && K
\ncdoff 
in $B\OO(n+m)$, whose bottom horizontal edges coincide. Thus, specifying $\mbf{U}$ at the level of $2$-morphisms already becomes a non-trivial task. The simplices of $\Vhook$ become only more complicated with increasing dimension.

The saving observation is the following. A $(k+1)$-simplex of $\Bbox\OO$ is, by the definition of the homotopy-coherent nerve, a functor $\Path[k+1]\to B^{\oplus}\OO$ of simplicial categories, where $\Path[k+1]$ is a simplicial enrichment of the standard simplex $\Delta[k+1]$. The most non-trivial part of the construction is the definition of the restriction 
\[
    \mbf{U}|\colon\Hom_{\Path[k+1]}(0,k+1)=\Nerve_{\bullet}(P)\to B\OO^{\infty}_{\amalg}
\]
with target the topological monoid of all infinite Grassmannians mentioned above, and domain the nerve of a certain poset $P$. 

Now, we specify in \Cref{HZBV0LB} a \emph{topological} $k$-simplex 
\[
    \btd\colon\Delta^{k}\hookrightarrow |\Nerve_{\bullet}(P)|
\]
such that we have a `collar-gluing'
\[
    |\Nerve_{\bullet}(P)|\cong \bigtriangleup\cup_{\partial}\btd
\]
where $\bigtriangleup$ is by definition the closure of the complement of the image of $\btd$, and $\partial$ is a certain $(k-1)$-simplex in $|\Nerve_{\bullet}(P)|$ giving the common boundary. Then we construct a new poset $P'$ together with a homeomorphism $|\Nerve_{\bullet}(P')|\cong\bigtriangleup$, which we show in \Cref{5KGP6JE} by a convexity argument. We thus reduce the construction of the restriction to $\mbf{U}|$ to two maps, one of type $\bigtriangledown\to B\OO^{\infty}_{\amalg}$ and one of type $\bigtriangleup\to B\OO^{\infty}_{\amalg}$, which agree on $\partial$. We finally show that a $k$-morphism of $\EEx(B\OO(n,m))$ determines exactly a map of type $\bigtriangledown\to B\OO^{\infty}_{\amalg}$, and that the other map can be induced canonically in a combinatorial manner. In terms of filling the simplices \eqref{LL2SEBY} and \eqref{I2XVTKD}, this construction translates to decomposing the rectangle obtained by gluing them along their common edge as 
\cdon[sep=small]
    & W'\oplus V \ar[dl]\ar[dr] \\
    W'\oplus V'\ar[dr,"{\gamma_{\oplus'}}"'] && W\oplus V \ar[ll,"{\gamma_{W}\oplus\gamma_{V}}"description] \ar[dl,"{\gamma_{\oplus}}"]\\
    & K
\cdoff 
instead, where $\bigtriangledown$ is the lower triangle, filled by the given exit $2$-morphism \eqref{9ASEI61}, the boundary $\partial$ is the horizontal edge, and $\bigtriangleup$ is the upper triangle, which can be filled canonically.  Reconstructing the original triangles from these fillers requires passing through the adjunction between geometric realisation and the singular chains functor. We show that this idea generalises functorially to all dimensions and to exit paths of all indices. In order to maintain functoriality, we construct $\mbf{U}$ inductively.

\subsubsection*{Acknowledgments} I thank K.\ \.{I}lker Berktav and Alberto S.\ Cattaneo for useful conversations, and Branko Juran for useful remarks and for making me aware of Lurie's \cite[01Z5]{kerodon}, which eventually led to \Cref{TUNKQ9R}. 

\section{Recollections}\label{app:nerves slices}

We will recall: the exit path $\infty$-category of a linked space following \cite{tetik2023linked}, the homotopy-coherent nerve of a Kan-enriched category following \cite{kerodon}, and the under-$\infty$-category construction, also following \cite{kerodon}. Our main aim is to fix notation: we will need not only the existence and basic properties of these constructions, but also their details.

\subsection{Exit paths}\label{6CW6YDL}

A \emph{linked space} is a span $\mfrak{S}=(M\xleftarrow{\pi}L\xrightarrow{\iota}N)$ of topological spaces where $\pi$ is a fibration and $\iota$ is a cofibration.
Its \emph{exit path $\infty$-category} $\EEx=\EEx(\mfrak{S})$ has the following simplices:
\begin{itemize}
    \item $\EEx_0=M_0\amalg N_0$.
    \item
    \(
        \EEx_k=M_k\amalg \Phat_{k-1} \amalg N_k
    \)
    for $k\geq1$.
\end{itemize}
Here and later, $X_{\ast}=\mrm{Sing}_*(X)$ for $X$ a topological space, and $\Phat_{k-1}$ is the set of pairs $(\gamma,e)$ where $\gamma\in N_{k}$, $e\in\{1,\dots,k\}$, such that the \emph{base} satisfies
\[
    \mbf{b}(\gamma,e)\coloneqq\gamma|_{0,\dots,e-1}\in\iota(L)_{e-1}.
\] 
We call such a pair an \emph{exit $k$-path of (exit) index $e$}. We will sometimes suppress $\iota$.

The face and degeneracy maps on the $M$- and $N$-components are given by those of $M$ and $N$. In order to define them on $\Phat_*$, suppose $(\gamma,e)\in\Phat_{k-1}$, and call, abusively, $d_i\gamma$ \emph{low} if it is in $\mbf{b}(\gamma,e)$, \emph{upper} if it is in $\gamma|_{e,\dots,k}$, and \emph{vertical} if it is neither low nor upper. Now, for $(\gamma,1)\in\Phat_0$, we set 
\begin{align*}
    d_1(\gamma,1)&=\pi (d_1\gamma)\in M_0,\\ 
    d_0(\gamma,1)&=\iota (d_0\gamma)\in N_0.
\end{align*}
For $k\geq2$, $(\gamma,e)\in\Phat_{k-1}$, and $d_i$ a face map: 
\begin{itemize}
    \item if $d_i\gamma$ is low, then we set $d_i(\gamma,e)=\pi (d_i\gamma)\in M_{k-1}$.
    \item if $d_i\gamma$ is upper, then we set $d_i(\gamma,e)=\iota (d_i\gamma)\in N_{k-1}$.
    \item if $d_i\gamma$ is vertical, then we set
    \(
        d_i(\gamma,e)=\left(d_i\gamma,\flat_{e,i}\right)\in\Phat_{k-2}. 
   \)
\end{itemize}
Finally, for $k\geq1$, $(\gamma,e)\in\Phat_{k-1}$, and $s_i$ a degeneracy, we set 
    \(
        s_i(\gamma,j)\coloneqq\left(s_i\gamma,\sharp_{e,i}\right)\in\Phat_k
    \).
Above, we have 
\[
    \flat^k_{e,j}=
        \begin{cases}
            e, & j\geq e\\
            e-1, & j<e
        \end{cases}
        \quad \text{and} \quad 
        \sharp^{k}_{e,j}=
        \begin{cases}
            e, & j\geq e\\
            e+1, & j<e.
        \end{cases}
\]
We will use this notation throughout this paper without further mention.\footnote{This is an economical version of the definition of $\EEx$ as explained in \cite[Remark 3.17]{tetik2023linked}.}
It is shown in \cite{tetik2023linked} that $\EEx$ defines a fully faithful functor from an $\infty$-category linked spaces into the $\infty$-category $\cat{C}{at}_\infty$ of all $\infty$-categories. The following facts will be relevant: we have equivalences
\[
    \Hom_{\EEx}(p,q)\simeq P_{L_p,q}
\]
for $p\in M$, $q\in N$, where $P_{L_p,q}$ is the space of paths in $N$ that start in the embedded fibre $\iota(L_p)$ and end in $q$; further, we have an equivalence
\[
    \down{M}{N}\coloneqq M\times_{\EEx^{\{0\}}}\EEx^{\Delta[1]}\times_{\EEx^{\{1\}}}N\simeq L.
\]

The following examples are of special practical interest:
\begin{enumerate} 
    \item If $M$ is a smooth manifold with boundary $\partial M$, choose a nowhere-vanishing inward-pointing normal vector field $X$ along $\partial M$ as well as a chosen time $t>0$, so that there is a closed embedding $\iota=\iota_{X,t}\colon \partial M\hookrightarrow M^\circ$. This gives rise to the linked space $\mfrak{M}_{\partial }=(\partial M\overset{\id}{\leftarrow}\partial M\overset{\iota}{\hookrightarrow}M^\circ)$. Then $\cat{E}{xit}(\partial M\subset M)\simeq\EEx(\mfrak{M}_\partial)$. The left-hand side is the exit path $\infty$-category of $M$ with its boundary stratification in the sense of \cite[\S A]{luriehigheralgebra} or, equivalently, \cite{ayala2017stratified}.
    \item If $M$ is a smooth manifold without boundary and $\Sigma\subset M$ a closed smooth submanifold, choose a tubular neighbourhood $\iota\colon\Normal\Sigma\hookrightarrow M$. Then, setting $\mfrak{M}_\Sigma=(\Sigma\leftarrow\mbb{S}(\Normal\Sigma)\overset{\iota}{\hookrightarrow}M\smallsetminus\Sigma)$, where $\mbb{S}$ takes the sphere bundle, we have $\cat{E}{xit}(\Sigma\subset M)\simeq\EEx(\mfrak{M}_\Sigma)$.
    \item Combining the previous points, let $M$ be a smooth manifold with or without boundary, stratified to depth $1$ by its boundary and a finite collection of non-intersecting closed smooth submanifolds. This gives rise to a linked space $\mfrak{M}$ such that $\cat{E}{xit}(M)\simeq\EEx(\mfrak{M})$.
\end{enumerate}

\subsection{The homotopy-coherent nerve}\label{app:simplicial nerves}

Let $\Cdelta$ denote the category of simplicial categories, that is, the category of categories enriched in $\sSet$.
We assume the reader is familiar with the nerve $\Nerve(C)=\Nerve_\bullet(C)\in\sSet$ of an ordinary category $C$.
We will recall the simplicial nerve construction (\cite[]{cordier1982notion}, though see also \cite[  00KT]{kerodon}), following \cite[{\S}1.1.5]{lurie2009higher}. We will then recall, following \cite{kerodon}, its mirror image, the homotopy-coherent nerve which will be our nerve of choice. 

Similarly to the Yoneda embedding \(\mbf{\Delta}\hookrightarrow\sSet\), \([k]\mapsto\Delta[k]\), which gives a simplicial set for each $k\in\mbb{N}$, there exists a functor 
\[
\mathfrak{C}\colon\mbf{\Delta}\rightarrow \Cdelta.
\]
\begin{definition}
    \begin{enumerate}
        \item The simplicial category $\mathfrak{C}[k]$ has the same objects as those of $[k]$, and the simplicial sets of morphisms in each $\mathfrak{C}[k]$ are given by 
        \begin{align*}
            \Hom_{\mathfrak{C}[k]}(i,j)\coloneqq\Nerve(P_{i,j}),
        \end{align*}
        where $P_{i,j}$, \(0\leq i,j\leq k\) is empty if $i>j$, and
        \[
            P_{i,j}=\{ I\subseteq\{i\leq a+1\leq\dots\leq j\}\subseteq[k] : a,b\in I\}
        \]
        if $i\leq j$. In other words, $P_{i,j}$ is the poset consisting of the subposets of $[k]$ that start at $i$ and $j$, with partial order $\preceq$ given by subset inclusions.
        
        For each triple $i\leq j\leq p$ in $[k]$, there is a map 
        \[
            P_{j,p}\times P_{i,j}\rightarrow P_{i,p}
        \]
        defined by taking unions. The ordinary nerve functor applied to these maps yields maps 
        \[
            \Hom_{\mathfrak{C}[k]}(j,p)\times\Hom_{\mathfrak{C}[k]}(i,j)\rightarrow\Hom_{\mathfrak{C}[k]}(i,p)
        \]
        of simplicial sets, which is associative since so is taking unions. 
        \item A map \(f\colon [l]\rightarrow[k]\) in $\mbf{\Delta}$ induces a map $\mathfrak{C}[l]\rightarrow\mathfrak{C}[k]$ as follows: on objects, it is given by $[l]\ni i\mapsto f(i)\in[k]$, and on the mapping posets it is given by 
        \(P_{i,j}\ni I\mapsto f(I)\in P_{f(i),f(j)}\), applying $\Nerve$ to which defines the map $f=\mathfrak{C}f\colon\mathfrak{C}[l]\rightarrow\mathfrak{C}[k]$.
    \end{enumerate}
\end{definition}

We call the $P_{i,j}$ \emph{mapping posets}, and their nerves \emph{mapping spaces}.

\begin{definition}
    The \emph{simplicial nerve} $\Ndelta(\mscc{D})=\Ndelta_\bullet(\mscc{D})$ of a simplicial category $\mscc{D}$ is the simplicial set whose set of $k$-simplices is defined by 
    \[
        \Ndelta_k(\mscc{D})\coloneqq\Hom_{\Cdelta}(\mathfrak{C}[k],\mscc{D}).
    \]
    This is contravariant in $[k]$ via the covariance of $\mathfrak{C}$. 
\end{definition}

In other words, $\Ndelta$ is the restriction of the Yoneda embedding $\Cdelta\rightarrow\pSh(\Cdelta)$ along $\mathfrak{C}\colon\mbf{\Delta}\rightarrow\Cdelta$.

\begin{definition} The \emph{homotopy-coherent nerve} $\Nhc(\mcal{D})=\Nhc_{\bullet}(\mcal{D})$ of a simplicial category $\mcal{D}$ is the simplicial set whose set of $k$-simplices is defined by
\[
    \Nhc(\mscc{D})_{k}\coloneqq\Fun(\Path[k],\mcal{A}),
\]
where
\(
    \Path[k]\coloneqq\mathfrak{C}[k]^{\mathrm{op}}.
\)
We write $\geq$ for the partial order thereon.
\end{definition}

Recall that, for any category $C$, we have  an isomorphism $\Nerve_{\bullet}(C^{\opp})\cong(\Nerve_{\bullet}C)^\opp$ (\cite[003Q]{kerodon}), where the RHS is the $\infty$-categorical opposite. Thus, 
\[
    \Hom_{\Path[k]}(i,j)\cong\Nerve(P_{i,j}^{\opp}),
\]
and we keep this superscript `$\opp$' throughout. This redundancy is meant to provide clarity. In \cite{lurie2009higher}, Lurie writes $P$ for our $P$, while in \cite{kerodon} he writes $P$ for our $P^{\opp}$. We exclusively work with $\Nhc$ (as in \cite{kerodon}), but we keep the superscript.
We call functors of type 
\(
    \Path[k]\rightarrow\mscc{A}
\) 
\emph{$k$-paths} in $\mscc{A}$.

If $\mscc{D}$ is Kan-enriched, then $\Ndelta(\mscc{D})$ is an $\infty$-category by \cite[Proposition 1.1.5.10]{lurie2009higher}. The same holds for $\Nhc(\mscc{D})$ by \cite[00LJ]{kerodon}. Both are variants of a result of Cordier and Porter \cite{cordier1986vogt}.

\subsection{Joins and (co)slices}\label{app:(co)slices}

For $f\colon K\rightarrow\mscc{C}$ a functor from a simplicial set to an $\infty$-category, there is (\cite{joyal2002quasi}, \cite[01GP]{kerodon}) a right fibration $\mscc{C}/f\rightarrow\mscc{C}$ and a left fibration $f/\mscc{C}\rightarrow\mscc{C}$, whose domains are respectively called the \emph{slice} and \emph{coslice} of $\mscc{C}$ at $f$. We will recall their definitions, but refer the reader to the op.\ cit.\ for the named lifting properties.

The following is equivalent to the more standard definition (see {\cite[0234]{kerodon}}). It is a simplicial version of Milnor's topological construction from \cite{milnor1956construction}.

\begin{definition}
    The \emph{join} \(X\star Y=(X\star Y)_\bullet\) of two simplicial sets $X=X_\bullet$, $Y=Y_\bullet$ is defined by 
    \[
        (X\star Y)_k=\{(\pi,f_-,f_+): \pi\colon\Delta[k]\rightarrow\Delta[1], f_-\colon\Delta[k]|_0\rightarrow X, f_+\colon\Delta[k]|_1\rightarrow Y\},
    \]
    where $\pi,f_-,f_+$ are maps of simplicial sets, and $\Delta[k]|_{i}=\{i\}\times_{\Delta[1]}\Delta[k]$, $i=0,1$, is defined using $\pi$. Given $\phi\colon[l]\rightarrow[k]$ in $\mbf{\Delta}$, the corresponding $\phi\colon\Delta[l]\rightarrow\Delta[k]$ defines a map $(X\star Y)_{k}\rightarrow(X\star Y)_l$ by restrictions.

\end{definition}
\begin{remark}\label{rmk:injections into joins}
    We have injections 
    \(
    \iota_0\colon X\hookrightarrow X\star Y\),  
    \(
    \iota_1\colon Y\hookrightarrow X\star Y.
    \)
    For the former, let $f\colon\Delta[k]\rightarrow X$ be a $k$-simplex of $X$. Defining 
    \(
        \pi\colon\Delta[k]\rightarrow \{0\}\hookrightarrow \Delta[1]
    \)
    and setting $f_-=f$, and necessarily $f_+\colon\emptyset\rightarrow Y$, gives a map $X_k\rightarrow (X\star Y)_k$, and similarly for $\iota_1$.
\end{remark}

\begin{remark}
    The join construction is functorial in both arguments. Given $\phi\colon X\rightarrow X'$, $\psi\colon Y\rightarrow Y' $, we write $\phi\star\psi$ for the induced map $X\star Y\rightarrow X'\star Y'$.
\end{remark}

\begin{definition}\label{dfn:slice}
    Let $K$ be a simplicial set, $\mscc{C}$ an $\infty$-category, and $f\colon K\rightarrow\mscc{C}$ a map. The \emph{slice} $\mscc{C}/f$ of $\mscc{C}$ at $f$ is the simplicial set defined by 
    \[
        (\mscc{C}/f)_n=(\Hom_{\sSet})_K(\Delta[n]\star K,\mscc{C}),
    \]
    where the subscript $K$ indicates that the set in question consists of maps 
    \(
        \phi\colon \Delta[n]\star K\rightarrow\mscc{C}
    \)
    whose precomposition 
    \(
        K\overset{\iota_1}{\hookrightarrow}\Delta[n]\star K\overset{\phi}{\rightarrow}\mscc{C}
    \)
    is $f$. 
    The face and degeneracy maps are given by precomposition and functoriality: a map $\psi\colon\Delta[m]\rightarrow\Delta[n]$ induces a map 
    \(
        \Delta[m]\star K\overset{\psi\star\id}{\rightarrow}\Delta[n]\star K\overset{\phi}{\rightarrow}\mscc{C},
    \)
    which is clearly in $(\mscc{C}/f)_m$, i.e., $(\phi\circ(\psi\star\id))|_K=f$.
    The slice is again an $\infty$-category.
    The projection $\mscc{C}/f\rightarrow\mscc{C}$ is given by precomposing $\phi\colon\Delta[n]\star K\rightarrow\mscc{C}$ with $\Delta[n]\overset{\iota_0}{\hookrightarrow}\Delta[n]\star K$.

    The \emph{coslice} $f/\mscc{C}$ is defined analogously, with $\Delta[n]\star K$ replaced by $K\star \Delta[n]$, $\iota_1$ by $\iota_0$ and vice versa, throughout. It is again an $\infty$-category.
\end{definition}

\begin{notation}
    Let $\iota_x\colon\Delta[0]\rightarrow\mscc{C}$ be given by a vertex $x\in\mscc{C}_0$. We write 
    \[
        \mscc{C}/x\coloneqq\mscc{C}/\iota_x, \quad x/\mscc{C}\coloneqq\iota_x/\mscc{C}.
    \]
    They are respectively called the \emph{over-} and \emph{under-$\infty$-category} at $x$.
\end{notation}

\begin{remark}\label{rmk:joins of simplexes}
    There are canonical isomorphisms 
    \[
        \Delta[k]\star\Delta[l]\simeq\Delta[k+1+l],
    \]
    such that the composition 
    \(
        \Delta[k]\overset{\iota_0}{\hookrightarrow}\Delta[k]\star\Delta[l]\overset{\sim}{\rightarrow}\Delta[k+1+l]
    \)
    is given by 
    \(
        [k]\hookrightarrow[k+1+l],\,  i\mapsto i,
    \) 
    and such that the composition 
    \(
        \Delta[l]\overset{\iota_1}{\hookrightarrow}\Delta[k]\star\Delta[l]\overset{\sim}{\rightarrow}\Delta[k+1+l]
    \)
    is given by 
    \(
        [l]\hookrightarrow[k+1+l],\, i\mapsto k+1+i.
    \)
\end{remark}

\begin{remark}\label{rmk:faces in an under-category}
    We should explicate the simplicial operators in an under-$\infty$-category $x/\mscc{C}$. 
    Via \Cref{rmk:joins of simplexes}, a $0$-simplex of $x/\mscc{C}$ is a $1$-simplex of $\mscc{C}$ with source $x$. Given a $1$-simplex 
    \(
        \gamma\colon\Delta[0]\star\Delta[1]\rightarrow \mscc{C}
    \)
    of $x/\mscc{C}$, its source and target $\gamma_0$, $\gamma_1$, are given, according to \Cref{dfn:slice}, by 
    \(
        \gamma_0\colon \Delta[0]\star\Delta[0]\overset{\id\star 0}{\hookrightarrow}\Delta[0]\star\Delta[1]\overset{\gamma}{\rightarrow}\mscc{C},
    \)
    and similarly with $\id\star 1$ for $\gamma_1$. The faces (and degeneracies) of simplices of all dimensions can be understood analogously: see \Cref{C4U7MUG}.
\end{remark}

\section{Linked tangent bundles}\label{sec:linked classifying maps}

In this section we will introduce tangent classifiers for linked smooth manifolds. Moreover, we will briefly consider variframed linked manifolds as a simple and illustrative example of a stratified tangential structure pulled back from AFR's context along the unpacking map, which we will assume.

\begin{definition}\label{36F1N2M}
    A linked space $\mfrak{S}=({M\xleftarrow{\pi}L\xrightarrow{\iota}N})$ is a \emph{linked (smooth) manifold} if $M$, $L$, and $N$ are smooth manifolds; $\pi$ is a fibre bundle,\footnote{For the construction of the tangent bundle of a linked manifold (\Cref{constr:tangent bundle of linked manifold}), it is enough that $\pi$ be a surjective submersion.}; and $\iota$ is a closed embedding.
\end{definition}

We will start by observing a fact that will let us give the tangent bundle of a linked manifold as an ordinary span map. This formalises an informal discussion present in \cite[{\S}2.1.4]{ayala2020factorization} in the conically-smooth setting. 
We assume all manifolds Hausdorff and paracompact so that vector sub-bundles split.

Let $\iota\colon L\hookrightarrow N$ be a closed embedding of smooth manifolds, and $E\to N$ a rank-$(n+m)$ vector bundle classified by $E\colon N\to B\OO(n+m)$, equipped with the inner product induced by that on the separable Hilbert space $H\cong\mbb{R}^{\infty}$ used to construct the Grassmannians $B\OO(k)=\Gr_{k}(H)$ as in \Cref{sec:direct sums} below. Let further $E_0$ be a rank-$n$ vector sub-bundle of $\iota^*E$, classified by $E_0\colon L\to B\OO(n)$. The pullback bundle itself is classified by $\iota^*E\colon L\hookrightarrow N\to B\OO(n+m)$. 

The normal bundle $E_{0}^{\perp} \subset \iota^*E$, classified by $E_{0}^\perp\colon L\to B\OO(m)$, satisfies $E_{0}\oplus E_{0}^\perp\cong \iota^*E$. 
It is classical that the Whitney sum is classified as follows: Consider the isomorphism \(\Phi\colon H\oplus H\cong H\) given by sending, with respect to a fixed basis of $H$ indexed over $\mbb{N}$, the first copy to odd coordinates and the second copy to even coordinates. The (abstract) direct sum precomposes with this isomorphism to give a map 
\begin{align*}
    E_{0}^\perp\widetilde{\oplus} E_{0}\colon L &\xrightarrow{E_{0}\times E_{0}^\perp}
    \Gr_{n}(H)\times \Gr_{m}(H)
    \\
    & \overset{\oplus}{\to}\Gr_{n+m}(H\oplus H)
    \\
    & \overset{\Phi}{\cong} \Gr_{n+m}(H)=B\OO(n+m).
\end{align*} 
The classifier 
\[\oplus_W\colon L\to B\OO(n+m)\] 
of the Whitney sum $E_{0}\oplus E_{0}^\perp$ is then homotopic to $E_{0}^\perp\widetilde{\oplus} E_{0}$.

\begin{lemma}\label{T6HNI7D}
    Let
    \begin{itemize} 
        \item $\iota\colon L\hookrightarrow N$ be a closed embedding of smooth manifolds,
        \item $E\to N$ be a rank-$(n+m)$ vector bundle equipped with an inner product, 
        \item and $E_0\hookrightarrow \iota^*E$ be a rank-$n$ vector sub-bundle. 
    \end{itemize}
    Then there exists a classifier $E\colon N\to \Gr_{n+m}(H\oplus H)$ of the isomorphism class of $E\to N$ such that the following diagram commutes:
    \begin{equation}\label{9IW9WA8}
        \begin{tikzcd}[column sep=large]
        L \ar[d,hook,"{\iota}"']\ar[r,"{E_{0}\times E_{0}^{\perp}}"] & \Gr_{n}(H)\times \Gr_{m}(H)\ar[d,"{\oplus}"] \\
        N \ar[r,"{E}"'] & \Gr_{n+m}(H\oplus H)
        \end{tikzcd}
    \end{equation}
\end{lemma}
\begin{proof}
    Let us concatenate the homotopy $\oplus_W\sim E^{\perp}_{0}\widetilde{\oplus}E_0$ constructed above with the standard one from $\iota^*E$ to $\oplus_W$, classifying the inverse of the bundle isomorphism $E_{0}\oplus E_{0}^\perp\cong \iota^*E$ given fibrewise by $(v,w)\mapsto v+w$, to obtain a homotopy 
\[h\colon \iota^*E\to\oplus_W\to E_{0}^\perp\widetilde{\oplus} E_{0}\] 
of maps $L\to B\OO(n+m)$ which sits in the commutative diagram 
\begin{equation*}
    \begin{tikzcd}
    L\ar[d,hook,"{\iota}"']\ar[r,"{h}"] & B\OO(n+m)^I \ar[d,"{\evv_0}"] \\
    N\ar[ur,dotted,"{H}"description] \ar[r,"{E}"']& B\OO(n+m)
    \end{tikzcd}
\end{equation*}
As $\iota$, being a closed embedding, is a cofibration, there exists a homotopy extension $H \colon N\to B\OO(n+m)^I$ as depicted. We may now consider 
\(E'\coloneqq H_{1}\colon N\to B\OO(n+m)\)
and apply the inverse isomorphism $\Phi^{-1}\colon H\cong H\oplus H$ to obtain \(\Phi^{-1}E'\colon N\to \Gr_{n}(H\oplus H).\) On the other hand, applying $\Phi^{-1}$ to $E_{0}^\perp\widetilde{\oplus} E_{0}$ recovers $E_{0}^\perp\oplus E_{0}=\oplus\circ E_{0}\times E_{0}^\perp$.  Therefore the two classifiers 
\(E^{\perp}_{0}\oplus E_{0}, \iota^*E'\colon L\to \Gr_{n+m}(H\oplus H)\)
coincide. 
\end{proof}

\begin{notation}\label{AGR9ZRZ}
    We will sometimes write simply $B\OO(k)$ for $\Gr_{k}(H^{\oplus -})\cong B\OO(k)$ for any countable number of copies of $H$ and thus abuse notation in diagrams of type \eqref{9IW9WA8}.
\end{notation}

Let now $\mathfrak{S}=(M\overset{\pi}{\leftarrow} L\overset{\iota}{\to} N)$ be a linked manifold, with each manifold riemannian. 
Given this contractible choice of metrics, we will show that there is a canonical map
\begin{equation}\label{47Y4NIX}
    \Tangent\mathfrak{S}\colon\mathfrak{S}\to B\OO(n,m)
\end{equation}
of linked spaces, which we call the \emph{tangent bundle} of $\mathfrak{S}$.

\begin{construction}\label{constr:tangent bundle of linked manifold} 
Since $\mathrm{d}\pi$ surjects, the induced linear dual map 
\(
    (\pi^*\Tangent M)^\vee\hookrightarrow (\Tangent L)^\vee,
\)
of bundles over $L$, injects. Using the metrics, this gives an injection
\(
    \pi^*\Tangent M\hookrightarrow \Tangent L.
\)
Composing with $\mathrm{d}\iota$, we have a bundle injection 
\(
    \pi^*\Tangent M\hookrightarrow\Tangent L\hookrightarrow\iota^*\Tangent N
\)
over $L$. Let us denote the normal bundle of this injection by 
\(
    \Normal\coloneqq\Normal_{N}M\coloneqq(\pi^*\Tangent M)^\perp\subset\iota^*\Tangent N.
\)
Now, in the diagram
\begin{equation}\label{0863VPP}
    \begin{tikzcd}[column sep=tiny]
        & L \ar[ddl,two heads,"{\pi}"']\ar[dr,hook,"{\iota}"] \ar[rrr,"{\pi^*\Tangent M\times \Normal}"] &&& B\mathrm{O}(n)\times B\mathrm{O}(m) \ar[dr,hook,"{\oplus}"] \\
        && N \ar[rrr,"{\Tangent N}",pos=0.25] &&& B\mathrm{O}(n+m)\\
        M \ar[rrr,"{\Tangent M}"] &&& B\mathrm{O}(n) \ar[from=uur,two heads,crossing over,"{\mathrm{pr}}"',pos=0.27] 
    \end{tikzcd}
    ,
\end{equation}
the back square
\begin{equation*}
    \begin{tikzcd}
    L \ar[d]\ar[r] & B\mathrm{O}(n)\times B\mathrm{O}(m)\ar[d] \\
    N \ar[r] & B\mathrm{O}(n+m)
    \end{tikzcd}
\end{equation*}
commutes using \Cref{T6HNI7D} and \Cref{AGR9ZRZ}, and the front square commutes trivially. This yields the span map \eqref{47Y4NIX}.
\end{construction}

\begin{remark}\label{P7AP3NA}
    Writing 
    \(
        \Normal_{L}M\coloneqq (\pi^*\Tangent M)^\perp\subset \Tangent L, 
    \)
    we have a splitting 
    \(
        \Tangent L\cong \pi^*\Tangent M\oplus \Normal_{L}M.
    \)
    Similarly, writing 
    \(
        \Normal_{N}L\coloneqq (\Tangent L)^\perp\subset \iota^*\Tangent N,
    \)
    we have a splitting 
    \(
        \iota^*\Tangent N\cong \Tangent L\oplus \Normal_{N}L\cong \pi^*\Tangent M\oplus \Normal_{L}M\oplus \Normal_{N}L.
    \)
    Thus 
    \[
        \Normal_NM\cong \Normal_LM\oplus\Normal_NL.
    \]
\end{remark}

The main examples of interest for the purpose of relating the linked tangent bundle to the conically-smooth version are those mentioned in \Cref{6CW6YDL}:

\begin{example}\label{JCLP17V}
    If $\mfrak{S}$ is induced by a closed submanifold $M\subset \overline{N}$ so that $L=\mbb{S}(\Normal_{N}M)$, the sphere bundle of the normal bundle, then $L$ has dimension $n+m-1$, $\Normal_{L}M$ has rank $m-1$, and $\Normal_LN$ has rank $1$. More specifically, in the conically-smooth context, the link (of a pair of strata) comes with an open embedding $L\times \mbb{R}\hookrightarrow N$, which is tantamount to a trivialisation $\Normal_NL\cong \varepsilon^1$, or, equivalently, to a diffeomorphism $L\times\mbb{R}\simeq \mbb{S}\left(\Normal_NM\right)\times\mbb{R}$.\footnote{For any vector bundle $E\to X$, a metric induces a diffeomorphism $E\smallsetminus X_0\cong\mbb{S}(E)\times(0,\infty)$ over $X$, where $X_0\subseteq E$ is the image of the zero-section.} This $\mbb{R}$-factor incarnates the extra $\mbb{E}_1$-structure featuring in the classification of constructible factorisation algebras (CFAs) on stratified spaces of type $M\subset \overline{N}$: for the special case of a hyperplane in euclidean space (see \cite{ayala2016factorization}).
\end{example}

\begin{example}\label{FZWKYYE}
    If $\mfrak{S}=(M\overset{=}{\leftarrow}M\hookrightarrow N)$ is induced by a manifold $\overline{N}$ with boundary $M=\partial\overline{N}$, so that $N=\overline{N}\smallsetminus\partial\overline{N}$, then $\Normal_LM=0$ and $\Normal_NL\simeq\varepsilon^1$ again. A cfa hereon in the case of a half-plane $\overline{N}=\mbb{H}^n$ is a Swiss-cheese algebra and is classified by Kontsevich's Swiss cheese conjecture.
\end{example}

\subsection{Adapting AFR-type structures}\label{T7HB1M7}

In \Cref{NVBP9DY}, we will introduce the stratified Grassmannian $\Vhook$ as discussed in the introduction, and in \Cref{YUXW1VO} construct a functor $\mbf{U}\colon\EEx(B\OO(n,m))\to\Vhook$. In this section, we will briefly discuss the implications. 

Applying $\EEx$ to $\Tangent\mfrak{S}$ and post-composing with $\mathbf{U}$, we obtain the induced map 
\[
    \Tangent\mfrak{S}\colon\EEx(\mathfrak{S})\to \Vhook,
\]
the linked version of the constructible tangent bundle of AFR \cite{ayala2020factorization}. Now, the $\infty$-category of \emph{tangential structures} according to op.\ cit.\ is the over-$\infty$-category
\(
    \Cinfty/\Vhook
\) 
.\footnote{See \cite[{\S}3]{lurie2009higher} for the $\infty$-category $\Cinfty$ of $\infty$-categories.} 
Via 
\[
\mathbf{U}^*\colon\Cinfty/\Vhook\to \Cinfty/\EEx(B\OO(n,m)),
\] 
these transfer to tangential structures on linked manifolds: given $\mathfrak{S}$, and writing $\mscc{B}_{(n,m)}\coloneqq\mathbf{U}^*\mscc{B}$, we may define the space (homotopy type) of $\mscc{B}$-structures on $\mathfrak{S}$ to be 
\[
    \Bred(\mathfrak{S})= \Hom_{/B\OO(n,m)}\left(\EEx(\mathfrak{S}),\mscc{B}_{(n,m)}\right)
\]
taken in $\Cinfty/\EEx(B\OO(n,m))$. In other words, \(\Bred(\mathfrak{S})=\Gamma\left(\left(\Tangent \mathfrak{S}\right)^*\mscc{B}_{(n,m)}\right),\) the space of homotopy-sections of $\left(\Tangent \mathfrak{S}\right)^*\mscc{B}_{(n,m)}\to \EEx(\mathfrak{S})$. 

By a \emph{smooth tangential structure} $\mathfrak{b}\colon \mscc{B}\to\Vhook$ we mean one that factors through $B\OO(k)\hookrightarrow \Vsim\hookrightarrow\Vhook$ for some $k$. 

\begin{example}\label{718DPIV}
    Let $\mathfrak{b}$ be a smooth tangential structure given by a map $BG\to B\OO(k)$ of spaces, e.g., induced by a map $G\to \OO(k)$ of topological groups, or a rank-$k$ bundle $X\to B\OO(k)$ on a space $X$. Then,
    $$
    B_{(n,m)}=
    \begin{cases}
        \EEx\left(\emptyset\leftarrow\emptyset\to B\right)=B, & n+m=k,\\
        \emptyset, & \text{else}.
    \end{cases}
    $$
    That is, a linked manifold $\mfrak{S}$ admits a $\mfrak{b}$-structure iff $\mfrak{S}=(\emptyset\leftarrow\emptyset\to N)$, $\dim N=k$, and $N$ admits a $\mfrak{b}$-structure.
\end{example}

\begin{example}[variframings]\label{NZG5TR7}
    Consider $\mbb{N}=(\mbb{N},\leq)$ with the standard order. For $k\in\mbb{N}$, let us write $\mbf{k}=\mbb{R}^k$. The tangential structure of \emph{variframing}, introduced in \cite{ayala2020factorization}, is given by 
    \begin{align*}
        \vfr\colon\mbb{N}&\to \Vhook,\\
        k&\mapsto\mbf{k},\\
        (k\leq K)&\mapsto(\mbf{k}\xhookrightarrow{-\oplus 0}\mbf{K}).
    \end{align*} 
    We read $\vfr(k\leq K)$ as the standard\footnote{Up to the choice of a pairing function $\mbb{N}\times\mbb{N}\cong\mbb{N}$: see \Cref{sec:direct sums} for its relevance.} isomorphism $\mbf{k}\oplus(\mbf{K-k})\cong\mbf{K}$. Let us restrict $\vfr$ to depth $1$ by choosing a pair $(n,m)\in\mbb{N}\times\mbb{N}$, i.e., consider $\vfr|_{n\leq n+m}\colon \{n\leq n+m\}\to \Vhook$. We have
    \begin{equation*}
        \mathbf{U}^*(\vfr|_{n\leq n+m})\simeq\EEx(\ast\leftarrow\ast\to\ast)\simeq\Delta[1],
    \end{equation*}
    so
    \(
        \mathbf{U}^*(\vfr|_{n\leq n+m})\to \EEx(B\OO(n,m))
    \) 
    is $\EEx$ of
    \begin{equation*}
        \begin{tikzcd}
         & \ast\ar[d,"{\mbf{n}\times\mbf{m}}"description]\ar[dl]\ar[dr] & \\
        \ast\ar[d,"{\mbf{n}}"description] 
        & B\OO(n)\times B\OO(m)\ar[dl,"\mathrm{pr}"description]\ar[dr,"{\oplus}"description] 
        & \ast\ar[d,"{\mbf{n+m}}"description]%
        \\
        B\OO(n) & & B\OO(n+m)
        \end{tikzcd}
        .
    \end{equation*}
    Thus, a \emph{variframing} on a linked manifold $\mathfrak{S}=(M\overset{\pi}{\twoheadleftarrow}L\overset{\iota}{\hookrightarrow} N)$ is a framing on $M$, a framing on $N$, and a framing on $\Normal_NM$, and an equivalence over $L$ between the direct sum framing on $\pi^*\Tangent M\oplus\Normal_NM$ and equivalent to framing on $\iota^*\Tangent N$ pulled back from $N$. 
\end{example}

\begin{example}[variframed point defects]\label{GZ4AU62}
    The choice of a point $p$ in a smooth manifold $N$ of dimension $m$ and a coordinate neighbourhood around it induces the linked manifold 
    \begin{equation*}
        \mathfrak{N}_p\coloneqq\left(\{p\}\leftarrow S^{m-1}\overset{\iota}{\hookrightarrow} N\smallsetminus\{p\} \right)
    \end{equation*} 
    where the sphere is the unit sphere in coordinates. The link map of 
    $\Tangent \mathfrak{N}_p$ reads 
    \begin{equation*}
        \varepsilon^0\times\left(\Tangent S^{m-1}\oplus \Normal(\iota)\right)\colon S^{m-1}\to \ast\times B\OO(m),
    \end{equation*}  
    i.e., 
    \[
        \iota^*\Tangent \left(N\smallsetminus\{p\}\right)\colon S^{m-1}\to B\OO(m).
    \]
    A $\vfr_{0\leq m}$-structure on $\mathfrak{N}_p$ is thus a framing on $N\smallsetminus{p}$, a framing on $\Tangent S^{m-1}\oplus\Normal(\iota)$, and an equivalence. If a framing of the former type exists, then it induces one of the latter type. Given one of former type, note that a variframing does \emph{not} require the standard induced framing on $\Normal(\iota)$; for instance, for $N=\mbb{R}^2$ equipped with the normal framing, not only the counter-clockwise but also the clockwise framing on $S^1$ extends to a (non-equivalent) variframing on $\mfrak{R}^2_{\{0\}}$.
\end{example}

\section{Quasi-(de)looping}\label{NVBP9DY}

\subsection{The additive Grassmannian}\label{sec:direct sums}

Let $H$ be a separable real Hilbert space of countably-infinite dimension, so, up to isometric isomorphism, the real sequence space $\ell^2$.

\begin{definition}
    For $k\in\mbb{N}$, $B\OO(k)\coloneqq\Gr_k(H)$ denotes the Grassmannian of $k$-dimensional subspaces of $H$.
\end{definition}

$B\OO(k)$ is an infinite-dimensional Hilbert manifold modelled on $H$, and, thus topologised, is homotopy equivalent to the colimit infinite Grassmannian 
\(
    \Gr_{k}(\mbb{R}^\infty)=\colim\,\Gr_k(\mbb{R}^n)
\)
along the closed embeddings 
\(
    \Gr_{k}(\mbb{R}^n)\hookrightarrow\Gr_{k}(\mbb{R}^{n+1})
\)
given by the first-coordinate inclusions $\mbb{R}^n\hookrightarrow\mbb{R}^{n-1}$.\footnote{This can be proved using results of Palais \cite[4 ff.]{palais1966homotopy}.} For our purposes, $\mbb{R}^\infty$, $H$, and $\ell^2$ are interchangeable.

\begin{notation}
    We set $B\OO_\amalg\coloneqq \coprod_{k\geq0}B\OO(k)$ and  $B\OO^+_\amalg\coloneqq\coprod_{k\geq1}B\OO(k)$.
\end{notation}

The aim of this section is to define a monoidal structure based on direct-summing of vector spaces, in the spirit of the direct sum operation
\(
    \oplus\colon\Gr_k(\mbb{R}^n)\times\Gr_{l}(\mbb{R}^m)\rightarrow\Gr_{k+l}(\mbb{R}^{n+m})
\)
Passing to infinite Grassmannians, these give maps \(B\OO(k)\times B\OO(l)\rightarrow\Gr_{k+l}(H\oplus H).\)
Choosing an isomorphism $H\oplus H\cong H$ yields a map
\(B\OO(k)\times B\OO(l)\rightarrow B\OO(k+l),\)
which defines a map
\begin{equation}\label{eq:old summing map}
    \oplus\colon B\OO_{\amalg}\times B\OO_{\amalg}\rightarrow B\OO_\amalg
\end{equation}
connected-componentwise.

The problem with this map is that there is no choice of an isomorphism $H\oplus H\cong H$ that would make the map above associative, so it does not promote $B\OO_{\amalg}$ to a topological monoid; suffice it to note that an isomorphism $H\oplus H\cong H$, or a pairing function (bijection) $\mbb{N}\times\mbb{N}\cong\mbb{N}$ cannot be associative, as this would contradict its injectivity. 
In order to attain hands-on access to the stratified Grassmannian, we have chosen to strictify $(B\OO_{\amalg},\oplus)$ instead. This gives a topological monoid at the cost of introducing some redundancy.

\begin{notation*}
    We set \(B\OO^N(k)\coloneqq \Gr_{k}(H^{\oplus N})\) and    \[B\OO^\infty_\amalg\coloneqq \{0\} \amalg \coprod_{N\geq1} \coprod_{k\geq1} B\OO^N(k).\]
\end{notation*}

\begin{remark}\label{rmk:redundant BO}
    Each $B\OO^N(k)$ is homeomorphic to $B\OO(k)=B\OO^1(k)$, but non-canonically. Given a pairing function and a choice of parenthesisation for large exponents, we have 
    \(
        B\OO^\infty_\amalg\cong \{0\}\amalg\mbb{Z}_+\times B\OO^+_\amalg.
    \)
\end{remark}

\begin{construction}\label{LJ4NG1F}
    The direct sum operation gives maps
    \(
        B\OO^N(k)\times B\OO^M(l)\rightarrow B\OO^{N+M}(k+l),
    \)
    which define an associative operation
    \[
        \oplus\colon B\OO^\infty_\amalg\times B\OO^\infty_\amalg\rightarrow B\OO^\infty_\amalg
    \]
    connected-componentwise. The zero vector space acts as the identity.
\end{construction}

    The canonical associativity of direct-summing of vector bundles on (paracompact Hausdorff) spaces translates to a monoidal structure on $B\OO_{\amalg}$ (or its stable version $B\OO$) only up to coherent homotopy. A systematic treatment in this direction is laid out in Lurie \cite{luriehigheralgebra}; see also Schwede \cite{schwede2018global}. The $E_{\infty}$-structure on $B\OO_{\amalg}$ is parametrised at arity $n$ by the (contractible) space of linear injections $H^n\hookrightarrow H$.
    \Cref{LJ4NG1F} is considered from a slightly different point of view in \cite[{\S}2]{schwede2018global}, where $(B\OO^{\infty}_{\amalg},\oplus)$ with $H$ relaxed to a vector space variable is called the \emph{additive Grassmannian}. 
    The relevant constructions and results of this paper apply immediately to the other varieties of Grassmannians such as the oriented, quaternionic, etc.: cf.\ \Cref{ADAUJCW} and \Cref{LP5KAHG}.

\begin{definition}\label{DPFKME2}
    By $B^\oplus\OO$ we denote the Kan-enriched category with a single object $\ast$, endomorphism space $B\OO^\infty_{\amalg}$, and composition $\oplus$. We write $\Bbox\OO\coloneqq\Nhc(B^\oplus\OO)$.
\end{definition}

Note that $\Bbox\OO$ is far from being an $\infty$-groupoid: only the zero vector space is invertible.

\subsection{The stratified Grassmannian}

Our definition of the \emph{stratified Grassmannian} is straightforward: it is the under-$\infty$-category of $\Bbox\OO$ under its unique object $\ast$. (See \Cref{app:(co)slices}.) It is that suggested in \cite[Remark 2.7]{ayala2020factorization} except for the strictification of $B\OO_{\amalg}$ into $B\OO^{\infty}_{\amalg}$ and for being a quasi-category rather than a Segal space.

\begin{definition}\label{def:Vhook}
    We call the $\infty$-category $\Vhook\coloneqq \ast/\Bbox\OO$ the (infinite) \emph{stratified Grassmannian}.
\end{definition}

\begin{remark}\label{TUNKQ9R}
    A theorem of Lurie \cite[01ZS]{kerodon} states that if $x\in\mscc{C}$ is an object of a Kan-enriched category $\mscc{C}$ and $x/\mscc{C}$ is the simplicial under-category as defined in \cite[01Z8]{kerodon}, then there is an equivalence of $\infty$-categories
    \[
        \Nhc(x/\mscc{C})\simeq x/\Nhc(\mscc{C})
    \]
    if for every morphism $f\colon x\to y$ and every object $z\in \mscc{C}$, pre-composition with $f$, 
    \[
        \Hom_{\mscc{C}}(y,z)\to\Hom_{\mscc{C}}(x,z),
    \]
    is a Kan fibration between Kan complexes.

    This is \emph{not} the case in our Kan-enriched category $B^{\oplus}\OO$, since 
    \begin{equation}\label{FPX1W4B}
        -\oplus V\colon B\OO^{\infty}_{\amalg}\to B\OO^{\infty}_{\amalg}
    \end{equation}
    is not a Kan fibration when $V\neq 0$. 
    Moreover, $\Vhook$ is indeed not equivalent to $\Nhc(\ast/B^{\oplus}\OO)$, as can be inferred by comparing their morphism spaces.
    Indeed, Let $V\in B\OO^m(l)$, $W\in B\OO^{n+m}(k+l)$. The objects of $\ast/B^{\oplus}\OO$ are the points of $B\OO^{\infty}_{\amalg}$, and we have, by definition, that $\Hom_{\ast/B^{\oplus}\OO}(V,W)$ is the ordinary fibre of \eqref{FPX1W4B} at $W$, so the subspace of $B\OO^{n}(k)$ consisting of $V'$ such that $V'\oplus V=W$. This is empty if, for instance, $W$ is spanned by a $(k+l)$-frame in the second summand of $H^{n}\oplus H^{m}$. If non-empty, it is a singleton. Now, by a result of Joyal (see Hebestreit--Krause \cite{hebestreit2020mapping} for a direct proof), morphism spaces in homotopy-coherent nerves coincide, up to equivalence, with those of the original Kan-enriched category, thus 
    \[
        \Hom_{\Nhc(\ast/B^{\oplus}\OO)}(V,W)\simeq\Hom_{\ast/B^{\oplus}\OO}(V,W).
    \]
    On the other hand, by Lurie \cite[Lemma 5.5.5.12]{lurie2009higher}, $\Hom_{\Vhook}(V,W)$ is the homotopy fibre of \eqref{FPX1W4B} at $W$, and so is equivalent to the space of paths in $B\OO^{n+m}(k+l)$ that start at $V'\oplus V$ for some $V'$ and end at $W$. In other words, we would have $\Vhook\simeq\Nhc(\ast/B^{\oplus}\OO)$ if and only if $B\OO^{\infty}_{\amalg}$ were equipped with the discrete topology. This justifies \Cref{def:Vhook}.
\end{remark}

Via the identification $\Delta[0]\star\Delta[n]\simeq\Delta[n+1]$, $n$-simplices of $\Vhook$ are ($n+1$)-simplices $\Bbox\OO$ with no qualification, which is to say that we have bijections
\[
    \Vhook_{n}\cong\Fun\paren{\Path[n+1],B^{\oplus}\OO}
\]
since $\Bbox\OO$ has a unique object. In particular, a $0$-simplex of $\Vhook$ is a `vector space' $V\in B\OO^{\infty}_{\amalg}$. 
The face and degeneracy maps of $\Vhook$ are simply those of $\Bbox\OO$ shifted once:

\begin{lemma}\label{C4U7MUG}
    $d_i^{\Vhook}=d_{i+1}^{\Bbox\OO}$, $s_{i}^{\Vhook}=s_{i+1}^{\Bbox\OO}$.
\end{lemma}

\begin{remark}\label{ZNETZEU}
    A $1$-morphism of $\Vhook$ with source $V_{01}$ and target $V_{02}$ is determined by a path $V_{01}\subseteq V_{12}\oplus V_{01}\overset{\gamma}{\rightarrow}V_{02}$. In this sense, morphisms of $\Vhook$ can be said to be `injections of vector spaces' if one disregards $\gamma$. Taking constant $\gamma$, and using the inner product on $H$ to choose the orthogonal complements canonically, includes proper vector space injections into the non-invertible morphisms of $\Vhook$. In view of \Cref{TUNKQ9R}, the naive identification of morphisms with injections amounts to equipping $B\OO^{\infty}_{\amalg}$ with the discrete topology.
\end{remark}

\subsection{The core of $\Vhook$}

It is desirable that $\Vhook$ contain no more invertible morphisms than the original infinite Grassmannians, so that no more information is added unnecessarily. \Cref{prop:Vsim=BO} below states just this. While this is clear from the construction, an explicit homotopy-inverse to the stated equivalence will have to wait until \Cref{EEYT903}.

\begin{notation*}
    For $\mscc{C}$ an $\infty$-category, let $\mscc{C}^\simeq$ denote its \emph{maximal sub-$\infty$-groupoid}, or \emph{core}, i.e.\ the $\infty$-groupoid whose $n$-simplices are exactly those $n$-simplices of $\mscc{C}$ whose edges are isomorphisms $\mscc{C}$. This is indeed an $\infty$-groupoid by a result of Joyal \cite{joyal2002quasi}; see also \cite[019D]{kerodon}. We will write $\Vsim\coloneqq\paren{\Vhook}^{\simeq}$.
\end{notation*}

\begin{notation}
    We will indicate (sub)posets by underlining their elements: e.g., using notation from \Cref{app:simplicial nerves}, we write $P^{\opp}_{0,2}=\{\underline{012}\geq \underline{02}\}$, and $P^{\opp}_{0,3}$ is as follows:
\cdon[sep=small]
    &\underline{0123}\ar[dl]\ar[dr]\ar[dd]\\
    \underline{013}\ar[dr] &   &\underline{023}\ar[dl]\\
        &\underline{03}
\cdoff
\end{notation}

\begin{proposition}\label{prop:Vsim=BO}
    There is an equivalence $\Vsim\simeq B\OO_{\amalg}^\infty\simeq \ast\amalg\mbb{Z}_+\times B\OO^+_{\amalg}$ of $\infty$-groupoids.
\end{proposition}

\begin{proof}
    First, given a $k$-simplex $\gamma$ of $B\OO^{\infty}_{\amalg}$, we will construct a functor 
    \[
        \Gamma\colon\Path[k+1]\to B^{\oplus}\OO
    \]
    of simplicial categories. It is necessarily trivial on objects. Let now $i\leq j\in[k+1]$ and let 
    \[
        \alpha=(\alpha^{0}\geq\cdots\geq\alpha^{n})\in\Nerve_{n}(P^{\opp}_{i,j})
    \]
    with subposets $\alpha^{x}=\underline{\alpha^{x}_{1},\dots,\alpha^{x}_{n_x}}$, $\alpha^{x}_{y}\in[k+1]$, $\alpha^{x}_{y}<\alpha^{x}_{y'}$ strictly for $y<y'$, and $\alpha^{x}_{1}=i$, $\alpha^{x}_{n_x}=j$. 
    If $j=i$, then all such sequences trivial and each $\alpha^x$ consists of $i$ alone, and therefore, by functoriality, $\Gamma(\alpha)=s^{n}_{0}(0)\in (B\OO^{\infty}_{\amalg})_n$, the $n$-fold degenerate zero vector space. Let us therefore assume $i<j$.
    If $i>0$, we also set 
    \(
        \Gamma(\alpha)\coloneqq s_0^{n}(0).
    \)
    If $i=0$ and so $j>0$ by assumption, every subposet $\alpha^{x}$ consists of at least two elements. Consider then the associated map 
    \begin{align*}
        A\colon [n]&\to [k]\\
        x&\mapsto \alpha^{x}_{2}-1.
    \end{align*}
    It is functorial since the partial order $\leq$ is defined to mean that $\alpha$ is given by subsets of $[k+1]$ satisfying $\alpha^{0}\supseteq\cdots\supseteq\alpha^{n}$, so $\alpha^{x'}_{2}\in\alpha^{x}$ and therefore $\alpha^{x}_{2}\leq \alpha^{x'}_{2}$ whenever $x\leq x'$. It is moreover well-defined since $\alpha^{x}_{2}-1\leq j-1\leq k$, and $\alpha^{x}_{2}-1>\alpha^{x}_{1}-1\geq 0$. 
    Now, the rule 
    \begin{align*}
        \Phi\colon\Nerve_{n}(P^{\opp}_{i,j})&\to\Delta[k]_{n}=\Hom_{\mbf{\Delta}}([n],[k])\\
        \alpha&\mapsto \Phi(\alpha)\coloneqq A
    \end{align*}
    is simplicial: let $\delta\colon[n']\to[n]$ be a poset map, so $(\delta^*(\alpha))^{x}=\alpha^{\delta(x)}$ for $x\in[n']$, and observe that $\Phi(\delta^*\alpha)(x)=\alpha^{\delta(x)}_{2}-1=\delta^*(\Phi(\alpha))(x)$.
    We thus obtain the maps
    \begin{align*}
        \Gamma\colon \Nerve_{n}(P^{\opp}_{i,j})&\to (B\OO^{\infty}_{\amalg})_{n}\\
        \alpha&\mapsto \Phi(\alpha)^*\gamma
    \end{align*}
    for $n\geq 0$ which assemble into maps 
    \[
        \Gamma\colon\Hom_{\Path[k+1]}(i,j)\to B\OO^{\infty}_{\amalg}
    \]
    for all pairs $i\leq j$ in $[k+1]$. The simpliciality of $\Phi$ the simpliciality of $\Gamma$ on morphism spaces.

    We will now show that $\Gamma$ is functorial. Let $\alpha\in\Nerve_{n}(P^{\opp}_{i,j})$, $\beta\in\Nerve_{n}(P^{\opp}_{j,l})$ be sequences as above, with $i\leq j\leq l$ in $[k+1]$, so 
    \[
        \beta\cup\alpha=\paren{\beta^{0}\cup\alpha^{0}\geq\cdots\geq\beta^{n}\cup\alpha^{n}}\in\Nerve_{n}(P^{\opp}_{i,l}).
    \]
    If $i=j=l$, then $\Gamma(\beta\cup\alpha)=s^{n}_{0}(0)\oplus s^{n}_{0}(0)=s^{n}_{0}(0)$; if $i=j<l$, then $\Phi(\beta\cup\alpha)(x)=(\beta^{x}\cup\alpha^{x})_{2}-1=\beta^{x}_{2}-1$ since $\alpha^{x}=(\alpha^{x}_{1})$ is degenerate and so $(\beta^{x}\cup\alpha^{x})_{1}=\alpha^{x}_{1}=\beta^{x}_{1}$, hence $\Gamma(\beta\cup\alpha)=\Gamma(\beta)=\Gamma(\beta)\oplus s^{n}_{0}0=\Gamma(\beta)\oplus\Gamma(\alpha)$; if $i<j\leq l$, then analogously $\Phi(\beta\cup\alpha)(x)=\alpha^{x}_{2}-1$, and so $\Gamma(\beta\cup\alpha)=\Gamma(\alpha)=\Gamma(\beta)\cup\Gamma(\alpha)$ because $j>0$ gives $\Gamma(\beta)=s^{n}_{0}(0)$ by construction. 

    Let us observe that the maps 
    \begin{align*}
        \Psi\colon(B\OO^{\infty}_{\amalg})_k&\to \Vhook_k\\
        \gamma&\mapsto\Psi(\gamma)\coloneqq\Gamma
    \end{align*}
    assemble into an $\infty$-functor $\Psi\colon B\OO^{\infty}_{\amalg}\to\Vhook$. By \Cref{C4U7MUG}, we must show that 
    \[
        \Psi(d_i\gamma)=d_{i+1}^{\Bbox\OO}(\Psi(\gamma)) \quad \text{and} \quad \Psi(s_{i}\gamma)=s_{i+1}^{\Bbox\OO}(\Psi(\gamma))
    \]
    for all $i\in[k]$. We will show the first and leave the second to the reader.
    Let us assume $k\geq1$, and take $j=0$, $l>0$, and 
    \(
        \alpha\in\Nerve_{n}({P^{\opp}_{0,l}})=\Hom_{\Path[k]}(0,l).
    \)
    The face map $\partial_{i}\colon[k-1]\hookrightarrow[k]$ composes with $\Phi(\alpha)\colon [n]\to [k-1]$ to give 
    \begin{align*}
        [n]&\to[k]\\
        x&\mapsto \begin{cases}
            \alpha^{x}_{2}-1, & \alpha^{x}_{2}-1\leq i-1,\\
            \alpha^{x}_{2}, & \alpha^{x}_{2}-1\geq i.
        \end{cases}\\
        &=\begin{cases}
            \alpha^{x}_{2}-1, &\alpha^{x}_{2}\leq i,\\
            \alpha^{x}_{2}, & \alpha^{x}_{2}\geq i+1
        \end{cases}
    \end{align*}
    so that $\Psi(d_{i}\gamma)(\alpha)$ is the pullback of $\gamma$ along this map. 
    On the other hand, $d_{i+1}^{\Bbox\OO}$ is given by pre-composing with 
    \begin{equation}\label{0TMHU4F}
        \partial_{i+1}\colon\Path[k]\hookrightarrow\Path[k+1],
    \end{equation}
    which on $\alpha$ reads $\partial_{i+1}(\alpha)=(\partial_{i+1}\alpha^{0}\geq\cdots\geq\partial_{i+1}\alpha^{n})\in\Hom_{\Path[k+1]}(0,\partial_{i+1}l)$ with 
    \(
        \partial_{i+1}\alpha^{x}=\underline{\partial_{i+1}(\alpha^{x}_{1}),\dots,\partial_{i+1}(\alpha^{x}_{n_{x}})}=\underline{0,\partial_{i+1}(\alpha^{x}_{2}),\dots,\partial_{i+1}l}
    \) and so 
    \[
        \partial_{i+1}(\alpha)^{x}_{2}-1=\partial_{i+1}(\alpha^{x}_{2})-1=\begin{cases}
            \alpha^{x}_{2}-1, & \alpha^{x}_{2}\leq i,\\
            \alpha^{x}_{2}, & \alpha^{x}_{2}\geq i+1
        \end{cases}.
    \]
    By the above, we obtain $\partial_i\Phi(\alpha)=\Phi(\partial_{i+1}\alpha)$ and thus 
    \begin{align*}
        \Psi(d_i\gamma)(\alpha)&=\Phi(\alpha)^*(d_i\gamma)=(\partial_i\Phi(\alpha))^*\gamma=(\Phi(\partial_{i+1}\alpha))^*\gamma\\
        &=\partial_{i+1}^*\Phi(\alpha)^*\gamma\\
        &=d_{i+1}^{\Bbox\OO}(\Psi(\gamma))(\alpha).
    \end{align*}
    Compatibility with degeneracies can be verified similarly.

    We have thus defined an $\infty$-functor 
    \[
        \Psi\colon B\OO^{\infty}_{\amalg}\to\Vsim\hookrightarrow\Vhook, 
    \]
    which necessarily factors as depicted. 
    Let now
    \[
        \sigma\colon\Path[k+1]\to B^{\oplus}\OO
    \]
    be a $k$-simplex of $\Vsim$, which is to say that the restrictions
    \[
        \sigma|_{ijl}\colon\Path\{i\leq j\leq l\}\hookrightarrow\Path[k+1]\xrightarrow{\sigma}B^{\oplus}\OO,
    \]
    as $1$-morphisms of $\Vhook$, are isomorphisms. By \Cref{C4U7MUG}, the relevant triples satisfy $i=0$, $j\geq1$, and  $\sigma|_{0jl}$ is fully determined by a path in $B\OO^{\infty}_{\amalg}$ of type 
    \[
        \sigma(\underline{jl})\oplus\sigma(\underline{0j})\to \sigma(\underline{0l})
    \]
    where $\sigma(\underline{0j})$ is the source of $\sigma|_{0jl}$ and $\sigma(\underline{0l})$ its target. However, since $\Vhook$ contains no morphism of type $W\to V$ if $\mrm{rk}(W)>\mrm{rk}(V)$, we conclude $\sigma(\underline{jl})=0$ since $\sigma|_{0jl}$ is an isomorphism. So, for any $\alpha=\underline{\alpha_1,\dots,\alpha_{n}}\in\Nerve_0({P^{\opp}_{j,l}})$ we have 
    \(
        \sigma(\alpha)=0
    \)
    by decomposing as $\alpha=\underline{\alpha_{n-1}\alpha_n}\cup\cdots\cup\underline{\alpha_1\alpha_0}$. Thus we obtain 
    \begin{equation}\label{B6Q4MHX}
        \sigma(\underline{0,\alpha})=\sigma(\underline{0,\alpha_1})
    \end{equation}
    by decomposing as $\underline{0,\alpha}=\alpha\cup\underline{0,\alpha_1}$.
    Now, as was noted in \Cref{TUNKQ9R}, $\Hom_{\Vhook}(V,W)$ is equivalent to the space of paths in $B\OO^{\infty}_{\amalg}$ that start $V'\oplus V$ for some $V'$ and end at $W$, and we have shown that within $\Vsim$ this reduces to $V'=0$, implying, for $V$ and $W$ in the same connected component of $B\OO^{\infty}_{\amalg}$, that 
    \(
        \Hom_{\Vsim}(V,W)\simeq\Hom_{B\OO^{\infty}_{\amalg}}(V,W)
    \)
    Moreover, along this equivalence, $\Psi$ maps as the identity on morphism spaces, whence it is fully-faithful. Since it is moreover a bijection on objects, we conclude  that it is an equivalence onto $\Vsim$.
\end{proof}

\begin{remark}\label{4J6WZPX}
    \Cref{B6Q4MHX} may seem to lead to the following point-set Ansatz for an inverse to $\Psi$, with the goal of promoting $\Vsim\simeq B\OO^{\infty}_{\amalg}$ to an isomorphism: Consider the sequence
    \[
        \Lambda^k=\paren{\underline{0,1,2,\dots,k+1}\geq\underline{0,2,3,\dots,k+1}\geq\cdots\geq\underline{0,k+1}}\in\Nerve_{k}\paren{P^{\opp}_{0,k+1}}.
    \]
    Its image $\sigma(\Lambda^k)\in (B\OO^{\infty}_{\amalg})_k$ is of type 
    \[
        \sigma(\underline{0,1})\to\sigma(\underline{0,2})\to\cdots\to\sigma(\underline{0,k+1}).
    \]
    The map
    \(
        \Psi^{-1}\colon\Vsim\to B\OO^{\infty}_{\amalg},\,
        \sigma\mapsto \sigma(\Lambda),
    \)
    however, is \emph{not} simplicial. To compare faces, let $i\in[k]$ so that we have $\sigma(d_i\Lambda^{k})=d_i(\sigma(\Lambda^{k}))$, while $(d_{i+1}^{\Bbox\OO}\sigma)(\Lambda^{k-1})=\sigma(\partial_{i+1}\Lambda^{k-1})$ holds by definition, so one might expect that $\sigma(d_{i}\Lambda^{k})=\sigma(\partial_{i+1}\Lambda^{k-1})$, giving compatibility with face maps.\footnote{Here, $\partial_{i+1}$ applies to $\Lambda^{k-1}$ as in \eqref{0TMHU4F}.} However, this need not hold: for instance, taking $k=3$ and $i=1$, this equation reads 
    \[
        \sigma(\underline{01234}\geq\underline{034}\geq\underline{04})=\sigma(\underline{0134}\geq\underline{034}\geq\underline{04}).
    \]
    While their vertices agree, the functoriality of $\sigma$ does not necessitate that these two simplices agree. 

Still, it is possible to give an explicit homotopy-inverse $\Vsim\to B\OO^{\infty}_{\amalg}$ using an idea we will develop later. Namely, to $\sigma$ we can associate its value on a certain topological $k$-simplex $\btd\subseteq|\Nerve_{\bullet}(P^{\opp}_{0,k+1})|$ which yields a $k$-simplex of $B\OO^{\infty}_{\amalg}$. We get back to this in \Cref{LP5KAHG}, using the construction of $\btd$ in \Cref{EEYT903}.
\end{remark}

The following suggests $\ast/\Nhc(-)$ as a means to adjoin non-invertible paths to a collection of non-stratified objects by means of a monoidal structure. 

\begin{corollary}\label{ADAUJCW}
    Let $M$ be a topological monoid whose only invertible element is its unit. Then 
    \(
        (\ast/\Nhc(BM))^{\simeq}\simeq M.
    \)
    Moreover, $(\ast/\Nhc(BM))\simeq\Nhc(\ast/BM)$ iff $M$ is discrete.
\end{corollary}
\begin{proof}
    The proof of \Cref{prop:Vsim=BO} applies mutatis mutandis. For the second statement, see \Cref{TUNKQ9R,ZNETZEU}, which also apply for the same reasons.
\end{proof}

\section{The unpacking map in low dimensions}\label{A8I4CSU}

Our next goal is to construct a map
\[
    \mathbf{U}:\EEx(B\mathrm{O}(n,m))\rightarrow\Vhook
\]
from the exit path $\infty$-category of the $(n,m)$-Grassmannian to the stratified Grassmannian. 
When restricted to $B\mathrm{O}(n)_\bullet$ and $B\mathrm{O}(n+m)_{\bullet}$ inside $\EEx\coloneqq\EEx(B\OO(n,m))$, we define it to be inclusion into the maximal sub-$\infty$-groupoid of $\Vhook$:
\[
    \EEx(B\OO(n,m))^{\simeq}\cong B\OO(n)\amalg B\OO(n+m)\subset B\OO^{\infty}_{\amalg}\overset{\Psi}{\hookrightarrow}\Vsim\subset\Vhook
\]
where $\Psi$ is defined in the proof of \Cref{prop:Vsim=BO}. 
It remains to define the restriction,
$$
\EEx_{k+1}\supset\Phat_k\rightarrow\Vhook_{k+1}\cong\Fun(\Path[k+2],B^{\oplus}\OO), 
$$
for $k\geq0$. We will explain dimensions $1$ and $2$ verbosely before giving the full construction without further motivation.

An element $(\gamma,1)\in\Phat_0$ is determined by a path $\gamma$ in $B\mathrm{O}(n+m)$ whose initial point is a direct sum $V_{12}\oplus V_{01}$ with $V_{01}\in B\mathrm{O}(n)$, $V_{12}\in B\mathrm{O}(m)$. Denoting the endpoint by $V_{02}$, $\gamma$ determines a $2$-morphism by arranging the data as
\[  
    \begin{tikzcd} 
        && \ast \\
        \ast &&&& \ast
        \arrow["{V_{01}}", from=2-1, to=1-3]
        \arrow["{V_{12}}", from=1-3, to=2-5]
        \arrow[""{name=0, anchor=center, inner sep=0}, "{V_{12}\oplus V_{01}}"', curve={height=30pt}, dashed, from=2-1, to=2-5]
        \arrow[""{name=1, anchor=center, inner sep=0}, "{V_{02}}", from=2-1, to=2-5]
        \arrow["{\, \gamma}"', shorten <=4pt, shorten >=4pt, Rightarrow, from=0, to=1]
    \end{tikzcd}
\]
which defines a $1$-morphism of $\Vhook$. 
We have thus defined 
\begin{equation}\label{NDWWNAH}
    \mathbf{U}_{\leq1}\colon\EEx_{\leq1}\rightarrow\Vhook.
\end{equation}
\begin{lemma}\label{N1P3DYW}
    The assignment $\mbf{U}_{\leq1}$ is functorial, i.e., compatible with all relevant face and degeneracy maps.
\end{lemma}
\begin{proof}
    Observe
    \(
        d^{\Vhook}_1(\mathbf{U}(\gamma,1))=V_{01}=\mathbf{U}(\mathrm{pr}(d^{B\OO(n+m)}_1\gamma))=\mathbf{U}(d^{\EEx}_1(\gamma,1)),
    \)
    and similarly
    \(
        d_0^{\Vhook}(\mathbf{U}(\gamma,1))=V_{02}=\mathbf{U}(d^{B\OO(n+m)}_{0}\gamma)=\mathbf{U}(d^{\EEx}_0(\gamma,1)).
    \)
    Compatibility with degeneracies is clear as in these dimensions there are no degenerate non-invertible paths.
\end{proof}

\begin{theorem}\label{YUXW1VO}
    The assignment $\mbf{U}|_{\leq1}$ of \eqref{NDWWNAH} extends to a fully faithful $\infty$-functor 
    \[
    \mathbf{U}\colon \EEx(B\OO(n,m))\hookrightarrow\Vhook.
    \]
\end{theorem}

We call $\mbf{U}$ the \emph{unpacking map}. We will give an inductive proof of existence. First, we will discuss how $\mathbf{U}$ must map $2$-morphisms. 

\begin{construction}\label{J08Z5ZM}
    Assignment \eqref{NDWWNAH}, $\mbf{U}_{\leq1}\colon\EEx_{\leq1}\to\Vhook$, extends by functoriality to the simplicial subset $\overline{\EEx_{\leq1}}\subset\EEx$ generated by $\EEx_{\leq1}$. That is, for any finite-order simplicial operator $\mathcal{O}=\prod\alpha_{i}$, $\alpha_{i}=s_{j_{i}}$ or $d_{j'_{i}}$ for any collection of indices $j$, $j'$ such that the application makes sense, we set 
    \(
        \mbf{U}(\mathcal{O}(\gamma,1))\mathcal{O}\mbf{U}_{\leq1}(\gamma,1)
    \)
    for any $(\gamma,1)\in\Phat_{0}$. \Cref{N1P3DYW} states exactly that this is well-defined. In effect, this gives a new definition only on degenerate simplices of  dimension higher than $1$ stemming from exit $1$-paths, so we could have considered  $\alpha=s_{j_{i}}$ only.
    We write
    \[
        \mbf{U}_{\leq1}\colon\overline{\EEx_{\leq1}}\to\Vhook
    \]
    for the resulting map of simplicial sets.
\end{construction}

We will extend $\mbf{U}$ along a filtration of type $\bigcup\overline{\EEx_{\leq k}}=\EEx$. The reader may skip to \Cref{H6QJ3HB}, referring back to \Cref{6WR2U3W} as necessary. 

\subsection{$2$-morphisms}\label{6WR2U3W} Exit paths in $\Phat_k$ come in $k+1$ classes according to their exit indices, which need to be mapped to $\Vhook$ in different ways. 

First, in order to uniquely determine a $3$-path in $B^{\oplus}\OO$, it is enough to map out of the sets $N_{\leq 2}(P_{i,j}^{\opp})$ into $B\OO_{\leq 2}\coloneqq (B\OO^{\infty}_{\amalg})_{\leq2}$, since higher dimensions are degenerate. Indeed, in general, a $\kappa$-path in $B^{\oplus}\OO$ is determined in dimensions $\leq\kappa-1$.

\begin{definition}\label{NRLTD63}
    We call those $1$-morphisms in $\Path[l]$ that are of type 
    \[
        \Nerve_{0}(P^{\mrm{op}}_{\alpha\beta})\ni\underline{\alpha\beta}=\{\alpha<\beta\}\subset[l]
    \]
    \emph{simple} and the remaining morphisms \emph{composite}; similarly for higher morphisms.
\end{definition}
\begin{example}\label{BQOMHSV}
    For instance, in $N_{1}(P^{\opp}_{1,4})$, $(\underline{1234}>\underline{124})=(\underline{234}\cup\underline{12}>\underline{24}\cup\underline{12})=(\underline{234}>\underline{24})\cup(\underline{12}>\underline{12})$ is composite, but $\underline{1234}>\underline{14}$ is simple. More generally, any arrow with target a pair is simple, and the others are composite.
\end{example}

The simple $1$-morphisms $\underline{\alpha\beta}\subset[3]$ in $\Path[3]$ are mapped by a $3$-path to $V_{\alpha\beta}\in B\OO_0$. When $\beta=\alpha+2$ (of which type there are two pairs), there are arrows $\underline{\alpha,\alpha+1,\beta}=\underline{\alpha+1,\beta}\cup\underline{\alpha,\alpha+1}>\underline{\alpha\beta}$ in $N_1(P_{\alpha\beta}^{\opp})$, which determine paths 
\(
    V_{\alpha+1,\beta}\oplus V_{\alpha,\alpha+1}\rightarrow V_{\alpha\beta},
\)
i.e.,
$
V_{12}\oplus V_{01}\rightarrow V_{02}$ {and} $V_{23}\oplus V_{12}\rightarrow V_{13}
$
in $B\OO_1$, namely two of the face $2$-paths. The remaining two faces are supplied analogously by considering $(\alpha,\beta)=(0,3)$ and the compositions $\underline{013}=\underline{13}\cup\underline{01}$ and $\underline{023}=\underline{23}\cup\underline{02}$. 
Finally, again for $(\alpha,\beta)=(0,3)$, consider $\underline{0123}=\underline{23}\cup\underline{12}\cup\underline{01}$, which is to be mapped to $V_{0123}=V_{23}\oplus V_{12}\oplus V_{01}$. Out of $N_1(P_{0,3}^{\opp})$ we receive paths $V_{0123}\rightarrow V_{03}$, $V_{0123}\rightarrow V_{013}$, $V_{0123}\rightarrow V_{023}$. The two non-degenerate elements $(\underline{0123}>\underline{023}>\underline{03})$ and $(\underline{0123}>\underline{013}>\underline{03})$ in $N_{2}(P_{0,3}^{\opp})$ are to map in $B\OO_{2}$ to
\begin{equation}\label{3-path data I}
    \begin{tikzcd}[column sep=small]
        & V_{023}\ar[dr]\\
        V_{0123}\ar[ur]\ar[rr] && V_{03}
    \end{tikzcd}
    =
    \begin{tikzcd}[column sep=small]
        & V_{23}\oplus V_{02}\ar[dr]\\
        V_{23}\oplus V_{12}\oplus V_{01}\ar[ur]\ar[rr] && V_{03}
    \end{tikzcd}
\end{equation}
and
\begin{equation}\label{3-path data II}
    \begin{tikzcd}[column sep=small]
        & V_{013}\ar[dr]\\
        V_{0123}\ar[ur]\ar[rr] && V_{03}
    \end{tikzcd}
    =
    \begin{tikzcd}[column sep=small]
        & V_{13}\oplus V_{01}\ar[dr]\\
        V_{23}\oplus V_{12}\oplus V_{01}\ar[ur]\ar[rr] && V_{03}
    \end{tikzcd}
\end{equation}
We have thus summed up the data needed to provide a functor $\Path[3]\rightarrow B^{\oplus}\OO$.

Now, let us start with paths of exit index $2\in\{1,2\}$.
Such an exit path $(\gamma,2)\in\Phat_1\subset\EEx_2$ consists of a $2$-simplex $\gamma\in B\OO(n+m)_2$ of type
\begin{equation}\label{diag:triangle index 2}
    \begin{tikzcd}
        & K\\
    W\oplus V \ar[r,"{\gamma_W\oplus\gamma_V}"']\ar[ur,"{\gamma_\oplus}"] & W'\oplus V'\ar[u,"{\gamma_{\oplus'}}"']
    \end{tikzcd}
\end{equation}
where the bottom edge comes from $B\OO(n)\times B\OO(m)$. 
If $(\gamma,2)$ is in $\overline{\EEx_{\leq1}}$, then $\mbf{U}(\gamma,2)$ is already defined by \Cref{J08Z5ZM}. Let us assume, therefore, that $(\gamma,2)$ is \emph{not} degenerate. 
The natural choice for the image, visualised as a $3$-simplex of $\Bbox\OO$, is
    \[
    \mathbf{U}(\gamma,2)=
    {\begin{tikzcd}
    &2\ar[dr,"{W'}"{description}]\\
1\ar[ur,"{0}"description]   &   &3\\
    &0\ar[ul,"{V}"description]\ar[uu,"{V'}"{description,pos=0.2}]\ar[ur,"{K}"{description}]
\ar[from=ul,to=ur,crossing over,"{W}"{description,pos=0.7}]
\end{tikzcd}}.\]
Indeed, the edges in \eqref{diag:triangle index 2} supply the face triangles. The fact that the bottom edge is of type $\gamma_W\oplus\gamma_V$ is crucial, since the summand paths supply the triangles adjacent to the edge decorated by the zero vector space. The only wrinkle is that the upper face requires a path $W'\rightarrow W$, which can be taken to be the (standard) inverse of $\gamma_W$, which we will denote by $\gamma_W^{-1}$.\footnote{Inversion in any dimension will be recalled in \Cref{X21WIW2}.} As for \eqref{3-path data I}, 
\begin{equation}\label{diag:qeoijngr}
    \begin{tikzcd}
        & W'\oplus V'\ar[dr,"{\gamma_{\oplus'}}"]\\
        W'\oplus V\ar[ur,"{\text{(i)}}"]\ar[rr,"{\text{(ii)}}"'] && K
    \end{tikzcd}
    ,
\end{equation}
note that we must yet choose the paths $W'\oplus V\rightarrow W'\oplus V'$ and $W'\oplus V\rightarrow K$ (corresponding to the arrows $\underline{0123}>\underline{023}$ and $\underline{0123}>\underline{03}$). To this end, consider the diagram
\[
     \begin{tikzcd}
        W'\oplus V' \ar[r,"{\gamma_{\oplus'}}"] & K \\
        W'\oplus V \ar[u,"{\id_{W'}\oplus\gamma_V}"] \ar[r,"{\gamma_W^{-1}\oplus\id_V}"'] & W\oplus V \ar[ul,"{\gamma_W\oplus\gamma_V}"description] \ar[u,"{\gamma_\oplus}"']
    \end{tikzcd}
\]
and choose the obvious fillers.
By $\id_A$ we mean the constant (degenerate) loop $s_{0}A$ at $A$. This suggests putting $\text{(i)}=\id_{W'}\oplus\gamma_V$, $\text{(ii)}=(\gamma_W^{-1}\oplus\id_V)\ast\gamma_\oplus$ 
(we write concatenation from left to right) whereupon the obvious filler can be chosen. Similarly, for \eqref{3-path data II}, i.e.,
\begin{equation}\label{diag:njqnj}
    \begin{tikzcd}
        & W\oplus V\ar[dr,"{\gamma_\oplus}"]\\
        W'\oplus V\ar[ur,"{\text{(i)}'}"]\ar[rr,"{\text{(ii)}'}"'] && K
    \end{tikzcd}
    ,
\end{equation}
consider 
\[
\begin{tikzcd}
    W\oplus V \ar[r,"{\gamma_{\oplus}}"] & K \\
    W'\oplus V \ar[u,"{\gamma_{W}^{-1}\oplus\id_V}"] \ar[r,"{\id_{W}\oplus\gamma_{V}}"'] & W'\oplus V' \ar[ul,"{\gamma_W^{-1}\oplus\gamma_V^{-1}}"description] \ar[u,"{\gamma_{\oplus'}}"']
\end{tikzcd}
\]
and proceed similarly. This completes the construction of $N_{\leq2}(P_{0,3}^{\opp})\rightarrow B\OO_{\leq2}$ and so of the $3$-path $\Path[3]\rightarrow B^{\oplus}\OO$ associated to the exit path $(\gamma,2)$. 

The image of an index-$1$ exit path 
\begin{equation}\label{TTXYP2R}
    {\begin{tikzcd}
        K \ar[r,"{\gamma_K}"] & K'\\
        W\oplus V \ar[u,"{\gamma_\oplus}"] \ar[ur,"{\gamma_{\oplus'}}"']
    \end{tikzcd}}
\end{equation}
is constructed analogously, with its picture as a $3$-simplex of $\Bbox\OO$ given by 
\begin{equation*}
    {\begin{tikzcd}
        &2\ar[dr,"{0}"{description}]\\
    1\ar[ur,"{W}"description]   &   &3\\
        &0\ar[ul,"{V}"description]\ar[uu,"{K}"{description,pos=0.2}]\ar[ur,"{K'}"{description}]
    \ar[from=ul,to=ur,crossing over,"{W}"{description,pos=0.7}]
    \end{tikzcd}}.
\end{equation*}
We have thus defined 
\(
\EEx_{\leq2} \rightarrow \Vhook_{\leq 2}.
\)
(A systematic account is given in \Cref{H6QJ3HB}). Writing 
\[
    \overline{\EEx_{\leq2}}\coloneqq\overline{(\EEx_{2}\smallsetminus\overline{\EEx_{\leq1}})\cup\overline{\EEx_{\leq1}}}
\] 
for the simplicial subset of $\EEx$ generated by $\overline{\EEx_{\leq1}}$ together with the elements of $\Phat_{1}$ that were not in $\overline{\EEx_{\leq1}}$, we claim that \Cref{J08Z5ZM} applies mutatis mutandis to yield a map 
\[
    \mbf{U}_{\leq2}\colon\overline{\EEx_{\leq2}}\to\Vhook,
\] 
so we must provide an analogoue of \Cref{N1P3DYW}. In particular, this will show that the contractible choices we have made along the way have had no bearing on functoriality.

\begin{lemma}\label{QX5R6Z7}
    The map $\mbf{U}_{\leq 2}$ is a well-defined extension of $\mbf{U}_{\leq 1}$.
\end{lemma}
\begin{proof}
Consider again an index-$2$ exit path $(\gamma,2)$ as in \eqref{diag:triangle index 2}. Its source edge is the path $d_2^{\EEx}(\gamma,2)=(\gamma_V\colon V\rightarrow V')\in B\OO_1\subset \EEx_1$. Since its two remaining edges are vertical, they are the elements of $\Phat_0$ induced by $\gamma_{\oplus'}$ and $\gamma_\oplus$. Now, recall that by \Cref{C4U7MUG} $\mathbf{U}(\gamma,2)$ is identified with an element of $\Bbox\OO_3\cong\Vhook_2$ via $\Delta[0]\star\Delta[2]\simeq\Delta[3]$. Accordingly, face maps apply on the factor $\Delta[2]$, i.e., $\partial$ acts as $\id_{\Delta[0]}\star\partial$. In the picture in $\Bbox\OO_3$, this means that when pulling back along a face map $\partial\colon\Delta[1]\hookrightarrow\Delta[2]$, we restrict to the triangle whose top edge is specified by $\partial$. For instance, $\partial_2\colon\Delta[1]\hookrightarrow\Delta[2]$, which skips $2$, applies to give 
\[d_2^{\Vhook}\mathbf{U}(\gamma,2)=
{\begin{tikzcd}
    &1 \ar[dr,"{0}"description]\\
    0\ar[ur,"{V}"description]\ar[rr,"{V'}"description] && 2
\end{tikzcd}}
\]
which is precisely $\mathbf{U}(\gamma_V)$. Compatibility with degeneracies holds by construction: $s_{i}\mbf{U}_{\leq2}(\gamma,2)=\mbf{U}_{\leq2}(s_{i}(\gamma,2))$ since $\mbf{U}_{\leq2}(s_{i}(\gamma,2))=s_{i}\mbf{U}_{\leq2}(\gamma,2)$ by construction, and this is well-defined \emph{if} $s_{i}(\gamma,2)\notin \overline{\EEx_{\leq1}}$, so that $\mbf{U}_{\leq2}$ does not clash with $\mbf{U}_{\leq1}$. But $s_{i}(\gamma,2)=\mathcal{O}(\gamma',1)$ would imply 
\begin{equation}\label{HE22OCQ}
    (\gamma,2)=d_{i}s_{i}(\gamma,2)=d_{i}\mathcal{O}(\gamma',1)=\mathcal{O}'(\gamma',1)\in\overline{\EEx_{\leq1}},
\end{equation}
which is precluded. This shows that $\mbf{U}_{\leq2}$ is a well-defined extension of $\mbf{U}_{\leq1}$.
We leave the analogous treatment of the remaining two faces and of the index-$1$ case to the reader. 
\end{proof}
\begin{remark}
    We should note that the top face 
    \[
        {\begin{tikzcd}
        &2 \ar[dr,"{W'}"description]\\
        1\ar[ur,"{0}"description]\ar[rr,"{W}"description] && 3
    \end{tikzcd}},
    \] 
    given by $\gamma_W^{-1}$ is \emph{not} a face in $\Vhook$, nor is $\gamma_W\in B\OO(m)_1\subset \EEx_1$ a face of $(\gamma,2)$ in $\EEx$.
\end{remark}

\begin{remark}
    We \emph{had to} assume that $(\gamma,2)$ is non-degenerate in order to apply the construction above, or else we would have broken functoriality. For instance, suppose the exit $2$-path in question is degenerate: $(\gamma,2)=s_{0}(\gamma',1)=(s_{0}\gamma',2)$. Then the construction above applied to $(\gamma,2)$ would yield $\text{(i)}=\id$ and $\text{(ii)}=\id\ast\gamma'$ in \eqref{3-path data I}, whereas the counterparts of these edges in $s_{0}(\mbf{U}(\gamma',1))$ are $\id$ and $\id$, respectively. These, as well as the filler $2$-simplices, are {not} necessarily the same. There does not seem to be a natural closed-form formula for $\mbf{U}$, hence our inductive-recursive construction.
\end{remark}

\begin{remark}\label{628I8LO}
    That a $2$-path $x\in\EEx_{2}$ is \emph{not} in $(\overline{\EEx_{\leq1}})_{2}\subseteq\EEx_{2}$ is equivalent to it not being degenerate in the usual sense. Indeed, say $x=\mcal{O}y$ for some simplicial operator $\mcal{O}$ and $y\in\EEx_{0/1}$. Then $\mcal{O}$ must contain at least one degeneracy, since otherwise $y$ would be in $\EEx_{\geq3}$ or $y=x$ already. So let $s$ be the last (left-most) degeneracy in $\mcal{O}$ and move it through the face maps in $\mcal{O}$ to the left of $s$ using simplicial identities, so that $\mcal{O}=s'\mcal{O}'$ for a resulting degeneracy $s'$. Then $x=s'\mcal{O}'y$ is degenerate. To prove this, one can use $d_{\alpha}s_{\beta}=s_{\beta-1}d_{\alpha}$ for $\alpha<\beta$, $d_{\alpha}s_{\beta}=\id$ for $\alpha\in\{\beta,\beta+1\}$ or $d_{\alpha}s_{\beta}=s_{\beta}d_{\alpha-1}$ for $\alpha>\beta+1$. After using the latter, move on to the next degeneracy in $\mcal{O}$, and iterate. If there is none, then we are again in the situation where $\mcal{O}$ only contains face maps so that $y\in\EEx_{\geq3}$, which is absurd. 
\end{remark}
    We will keep using these simplicial closures for convenience, as they simplify some arguments -- notably the latter part of the proof of \Cref{QX5R6Z7} which concerns degeneracies. 

\section{The proof of Theorem \ref{YUXW1VO}}\label{H6QJ3HB}

We will now give the general construction, starting with a systematic account of 
\[
    \Phat_1\to\left[\Path[3],B^{\oplus}\OO\right]
\]
in such a way that the ideas generalise to all dimensions.

\begin{notation*}
    We write $[A,B]\coloneqq\Fun(A,B)$.
\end{notation*}

\begin{definition}\label{NY1KJPW}
    We write $\mfrak{N}\colon B\OO(n)\times B\OO(m)\rightarrow B\OO(m)$ to denote the second coordinate projection. When we apply $\mfrak{N}$ to a low face of an exit path $(\gamma,j)$, we mean that, first, the corresponding face of $\gamma$ is to be taken, which is then (unambiguously) to be identified with a simplex of the link $B\OO(n)\times B\OO(m)$, and then $\mfrak{N}$ is to be applied. Namely, we have, by abuse of notation, a map 
    $$
    \mfrak{N}\colon\Phat_\ast\rightarrow B\OO(m)_{\ast'}
    $$
    for each $\ast\geq0$, given by the composition
    \begin{equation*}
    \begin{tikzcd}
        \Phat_\ast\ar[rr,bend right=20,dashed,"{\mfrak{N}}"']\ar[r,"{(\gamma,j)\mapsto\gamma|_{0,\dots,j-1}}"] &[4em] \mcal{L}_{\ast'}\ar[r,"{\mfrak{N}=\mrm{pr}_2}"] &[1em] B\OO(m)_{\ast'}
    \end{tikzcd}
    \end{equation*}
    where we have not omitted $\ast$ since $\Phat$ is not a simplicial set, nor is this map simplicial. The result is degenerate unless the exit index is maximal. We call 
    \[
        \mfrak{N}(\gamma,j)\in B\OO(m)_{j-1}
    \] 
    the \emph{normal component} of $(\gamma,j)$.
\end{definition}

\begin{notation}\label{X21WIW2}
    Let $X$ be a space. We denote by $\mathrm{Op}$ the canonical isomorphism $X\simeq X^{\opp}$ of Kan complexes, by which we mean $\Sing_\bullet(X)\simeq \Sing_\bullet(X)^{\opp}$. It is given by inverting $n$-simplices by pulling back along the maps
    \[
        \mrm{Op}\colon\Delta^{n}\to\Delta^{n},\ (a_{0},\dots, a_{n})\mapsto(a_{n},\dots,a_{0})
    \] 
    of the standard topological $n$-simplex, $n\geq0$ (\cite[003R]{kerodon}). For example, the `canonical inverse' of a path $\gamma$ is $\gamma^{-1}=\mrm{Op}(\gamma)=(\Delta^{1}\xrightarrow{\mrm{Op}}\Delta^{1}\xrightarrow{\gamma} X)$, so that $d_{0}\gamma^{-1}=d_{1}\gamma$, $d_{1}\gamma^{-1}=d_{0}\gamma$. More generally, for any simplicial set $S$, we have $S^{\opp}_{\ast}=S_{\ast}$ dimension-wise, and for any $\ast=n\geq0$, $d^{S^{\opp}}_{i}=d^{S}_{n-i}$ and $s_{i}^{S^{\opp}}=s_{n-i}^{S}$.
\end{notation}

\begin{notation}\label{N2HXQ2G}
    Let $\alpha_0,\dots,\alpha_\ell\in[k]$. By $\Path[\alpha_0,\dots,\alpha_\ell]\cong\Path[\ell]$ we denote the full simplicial subcategory of $\Path[k]$ generated by the objects $\alpha_0,\dots,\alpha_\ell$.
\end{notation}

\begin{construction}\label{STG38ZQ}
We define $\mbf{U}$ on $\Phat_0\rightarrow \left[\Path[2],B^{\oplus}\OO\right]$ by
 \begin{equation*}
    (\gamma,1)=
\begin{tikzcd}
    K\\
    W\oplus V\ar[u,"{\gamma}"]
\end{tikzcd}
\mapsto
\begin{cases}
    \Nerve_{0}(P_{0,a}^\opp) \to B\OO_0 , &  \underline{0,a}\mapsto \{a-1\}^*(\gamma,1)=\begin{cases}
        V, & a=1\\
        K, & a=2
    \end{cases}\\
    \Nerve_0(P_{1,2}^\opp)\to B\OO_0, & \underline{12} \mapsto \mrm{Op}\mfrak{N}(\gamma,j)=\mfrak{N}(\gamma,j)=W\\ 
    \Nerve_1(P_{0,2})\to B\OO_1, & (\underline{012}>\underline{02})\mapsto \gamma.
\end{cases}
\end{equation*}
This induces $\mbf{U}_{\leq1}\colon\overline{\EEx_{\leq1}}\to\Vhook$ as in \Cref{J08Z5ZM}.
\end{construction}

\subsection{The map $\Phat_1\rightarrow\left[\Path[3],B^{\oplus}\OO\right]$}\label{0N49067} While skipping to the induction step below is now possible, we will treat this dimension explicitly in order to illustrate the ideas. 
Let $(\gamma,j)\in\Phat_1$ and assume $(\gamma,j)\notin\overline{\EEx_{\leq1}}$.

\paragraph{{\bf Induced faces}}  We first define the faces of $\mathbf{U}(\gamma,j)$. The faces $d^{\Vhook}_{i}\mathbf{U}(\gamma,2)$, $i\in\{0,1,2\}$, are defined by $\mathbf{U}_{\leq1}$ via functoriality, i.e., by 
\begin{equation}\label{Z57UYYD}
d_i^{\Vhook}(\mathbf{U}_{\leq2}(\gamma,j))=\mathbf{U}_{\leq1}(d_i^{\EEx}(\gamma,j)).
\end{equation} 
This fixes $\mathbf{U}_{\leq2}(\gamma,j)$ by \Cref{C4U7MUG} on the subcategories $\Path[0,k,l]\cong\Path[2]$ of $\Path[3]$, for $1\leq k<l\leq3$ (\Cref{N2HXQ2G}). This is consistent due to the well-definedness of $\mbf{U}_{\leq1}$ as shown in \Cref{N1P3DYW}.

\paragraph{{\bf The top face}} The remaining face $d_0^{\Bbox\OO}(\mathbf{U}(\gamma,j))$ is the restriction to the full simplicial subcategory $\Path[1,2,3]$. Any pair $1\leq \alpha< \beta\leq 3$ specifies a restriction along $\Path[\alpha,\beta]\hookrightarrow \Path[1,2,3]$, and writing $\delta$ for the remaining element of $\{1,2,3\}$, these restrictions must by functoriality coincide with 
\begin{align*}
    d_{\delta-1}^{\Bbox\OO}d_{0}^{\Bbox\OO}(\mbf{U}(\gamma,j))&=d_{0}^{\Bbox\OO}d_{\delta}^{\Bbox\OO}(\mbf{U}(\gamma,j))
    \overset{}{=}d_{0}^{\Bbox\OO}d_{\delta-1}^{\Vhook}(\mbf{U}(\gamma,j))\\
    &=d_{0}^{\Bbox\OO}\mbf{U}_{\leq1}(d_{\delta-1}^{\EEx}(\gamma,j))
\end{align*}
where we used the simplicial identity $d_{i}d_{j}=d_{j-1}d_{i}$ for $i<j$ for the first equation and \Cref{C4U7MUG} for the second. In other words, the edges of the $d_{0}^{\Bbox\OO}$-face are determined already by $\mbf{U}_{\leq1}$. 
Therefore, only the restriction to the $1$-morphisms of $\Nerve_{\bullet}(P_{1,3}^{\opp})$, that is,
\begin{align*}
\Nerve_1(P_{1,3}^\opp) &\rightarrow B\OO_1\\
(\underline{123}=\underline{23}\cup\underline{12}>\underline{13}) &\mapsto \left(\mathbf{U}_{(\gamma,j)}(\underline{23})\oplus\mathbf{U}_{(\gamma,j)}(\underline{12})\rightarrow \mathbf{U}_{(\gamma,j)}(\underline{13})\right)\\
&=\left(\mathbf{U}(d_0d_0(\gamma,j))\oplus \mathbf{U}(d_0d_2(\gamma,j))\rightarrow\mathbf{U}(d_0d_1(\gamma,j))\right)
\end{align*}
remains to be specified. We put
\begin{equation}\label{NN7L78Y}
\mathbf{U}_{(\gamma,j)}|\Nerve_1(P_{1,3}^\opp)\coloneqq \Opp\mfrak{N}(\gamma,j).
\end{equation}

\begin{remark}
    We should note that it is immaterial that \eqref{NN7L78Y} is `not functorial' (although $\mathbf{U}$ will be). As we will see in \Cref{Q52DBDN}, the direct sums appearing in the $d_0^{\Bbox\OO}$-face are trivial in that all summands but one are zero, the non-zero one being determined by the exit index $j$.  We use $\Opp\mfrak{N}$ to supply \emph{only the path} in $B\OO(m)$. We have $\mathbf{U}_{(\gamma,j)}(\underline{23})=0$ if $j=1$, and $\mathbf{U}_{(\gamma,j)}(\underline{12})=0$ if $j=2$. If $j=1$, the edges of $d_0^{\Bbox\OO}$ are specified by simpliciality as in 
    \[
    {\begin{tikzcd}
    & 2\ar[dr,"{0}"description] \\
    1\ar[ur,"{W}"description]\ar[rr,"{W}"description] && 3 
    \end{tikzcd}},
    \]
    and $\Opp\mfrak{N}(\gamma,1)$ is $\Opp(\id_W)=\id_W\colon W=0\oplus W\to W$. Here, \eqref{NN7L78Y} happens to be functorial as $d_0^{\Bbox\OO}$ happens to lie in $\Vsim$. However, if $j=2$, the filler of 
    \[
    {\begin{tikzcd}
     & 2\ar[dr,"{W'}"description]\\
    1\ar[ur,"{0}"description]\ar[rr,"{W}"description] && 3
    \end{tikzcd}}
    \]
    is supplied by $\Opp\mfrak{N}(\gamma,2)=\Opp(\gamma_W)=\gamma_W^{-1}\colon W'=0\oplus W'\to W$. While $d_0^{\Bbox\OO}$,, being a morphism in $\Vhook$ from $0$ to $W$, is {not} invertible, the path is as desired.
\end{remark}

\paragraph{{\bf $1$-paths induced by functoriality}} Some $1$-paths in the image of $\mathbf{U}(\gamma,j)$ are determined by the data provided thus far and by imposing functoriality. We have the following decompositions: 
\begin{align*}
    \text{(i)}=(\underline{0123}>\underline{023})&=(\underline{23}\cup \underline{012}>\underline{23}\cup\underline{02})\\
    &=\id_{\underline{23}}\cup[\underline{012}>\underline{02}]\\
    &\in \mathrm{Im}\left(\cup: \Nerve_1(P_{2,3}^\opp)\times \Nerve_1(P_{0,2}^\opp) \to \Nerve_1(P_{0,3}^\opp)\right),\\
    \text{(i)}'=(\underline{0123}>\underline{013})&=(\underline{123}\cup\underline{01}>\underline{13}\cup\underline{01})\\
    &=[\underline{123}>\underline{13}]\cup\id_{\underline{01}}\\
    &\in \mathrm{Im}\left(\cup: \Nerve_1(P_{1,3}^\opp)\times \Nerve_1(P_{0,1}^\opp) \to \Nerve_1(P_{0,3}^\opp)\right).
\end{align*}
Thus, functoriality imposes 
\begin{align}
    \label{1MYHB5J}
    \mathbf{U}_{(\gamma,j)}\text{(i)}&=\id_{\mathbf{U}_{(\gamma,j)}(\underline{23})}\oplus\mathbf{U}_{(\gamma,j)}(\underline{012}>\underline{02}),\\
    \label{7FKA9DU}
    \mathbf{U}_{(\gamma,j)}\text{(i)}'&=\mathbf{U}_{(\gamma,j)}(\underline{123}>\underline{13})\oplus \id_{\mathbf{U}_{(\gamma,j)}(\underline{01})}
\end{align}
where the first non-constant summand is determined by \eqref{Z57UYYD} and the second by \eqref{NN7L78Y}. In particular, this systematises the ad hoc assignments in \eqref{diag:qeoijngr} and \eqref{diag:njqnj}.

\paragraph{{\bf The remaining $1$-path and $2$-paths}} It remains to specify $\mathbf{U}(\gamma,j)$ on the `long path' in $\Nerve_1(P_{0,3}^\opp)$ and on $\Nerve_2(P_{0,3}^\opp)$.
In contrast to paths induced by functoriality, $\underline{03}$ is simple, so 
\begin{align*}
    \text{(ii)}=(\underline{0123}>\underline{03}) 
\end{align*}
presents a genuine choice. It is considered most naturally in conjunction with the two non-degenerate elements in $N_2(P_{0,3}^\opp)$ to be mapped, as it is their common composition:
\begin{equation}\label{9681FJ4}
    \begin{tikzcd}
    & \underline{0123}\ar[dd,dotted,"{\text{(ii)}}"description]\ar[dl,"{\id_{\underline{23}}\cup[\underline{012}>\underline{02}]}"']\ar[dr,"{[\underline{123}>\underline{13}]\cup\id_{\underline{01}}}"]\\
    \underline{023}\ar[dr,"{\eqref{Z57UYYD}}"'] && \underline{013}\ar[dl,"{\eqref{Z57UYYD}}"]\\
     & \underline{03}
    \end{tikzcd}
\end{equation}
First, note that, regardless of exit index, this square decomposes into two triangles:
\begin{equation}\label{1HM4M30}
    \begin{tikzcd}
    & \underline{0123}\ar[dl,"{\id_{\underline{23}}\cup[\underline{012}>\underline{02}]}"']\ar[dr,"{[\underline{123}>\underline{13}]\cup\id_{\underline{01}}}"]\\
    \underline{023}\ar[dr,"{\eqref{Z57UYYD}}"'] && \underline{013}\ar[ll,"{d_2\gamma}"description]\ar[dl,"{\eqref{Z57UYYD}}"]\\
     & \underline{03}
    \end{tikzcd}
\end{equation}
For $j=2$ (using the labels in \eqref{diag:triangle index 2}), this reads 
\[
    \begin{tikzcd}
     & W'\oplus0\oplus V=W'\oplus V \ar[dl,"{\id_{W'}\oplus \gamma_V}"']\ar[dr,"{\gamma_W^{-1}\oplus\id_V}"] \\
    W'\oplus V'\ar[dr,"{\gamma_{\oplus'}}"'] && W\oplus V \ar[ll,"{\gamma_{W}\oplus\gamma_{V}}"description] \ar[dl,"{\gamma_{\oplus}}"]\\
     & K
    \end{tikzcd}
    ,
\]
and for $j=1$ (using the labels in \eqref{TTXYP2R}),
\[
    \begin{tikzcd}
     & 0\oplus W\oplus V=W\oplus V\ar[dl,"{\gamma_{\oplus}}"']\ar[dr,"{\id_W\oplus\id_V}"] \\
    0\oplus K=K \ar[dr,"{\gamma_K}"'] && W\oplus V \ar[ll,"{\gamma_{\oplus}}"description]\ar[dl,"{\gamma_{\oplus'}}"]\\
     & K'
    \end{tikzcd}
    .
\]
For both indices, the bottom triangle is filled by $\gamma$ itself, and the top one has a canonical filler. This suggests assigning to $\text{(ii)}$ the outer-left concatenation: 
\begin{equation}\label{VK71H7N}
    \mathbf{U}_{(\gamma,j)}\text{(ii)}=\mathbf{U}_{(\gamma,j)}(\underline{0123}>\underline{023})\ast\mathbf{U}_{(\gamma,j)}(\underline{023}>\underline{03}).
\end{equation}
Accordingly, $\mathbf{U}(\gamma,j)|_{\Nerve_2(P_{0,3}^\opp)}$ is determined by said fillers. 

Let us specify the filler of the top triangle in the case $j=2$ explicitly.\footnote{This will be systematised in \Cref{EEYT903} and can be skipped. In fact, we will prove that even this concatenation needn't be explicated. The case $j=2$ is of central importance to the proof of \Cref{VMR7DN8}.} There is an intermediate triangle
\begin{equation}\label{T4UNR1L}
    \begin{tikzcd}
    W'\oplus V\ar[dr,"{\gamma_W^{-1}\oplus\id}"]\ar[d,"{\id\oplus \gamma_{V}}"'] \\
    W'\oplus V &[1em] W\oplus V\ar[l,"{\gamma_W\oplus\gamma_V}"]
    \end{tikzcd}
\end{equation}
filled by the direct sum of the degenerate $2$-path $s_0(\gamma_V)$, i.e., 
\[
    \begin{tikzcd}
    V\ar[dr,"{\id_V}"]\ar[d,"{\gamma_V}"']\\
    V' & V\ar[l,"{\gamma_V}"]
    \end{tikzcd}
\]
and the $2$-path 
\[
    \begin{tikzcd}
    W'\ar[dr,"{\gamma_W^{-1}}"]\ar[d,"{\id}"']\\
    W' & W\ar[l,"{\gamma_W}"]
    \end{tikzcd}
\]
given by
\[
    \Gamma_W\colon\Delta^2\to B\OO(m),\ (t_0,t_1,t_2)\mapsto\gamma(t_1,1-t_1).
\] 
Indeed, we see that its edges are as desired:
\begin{align*}
    &(d_0\Gamma_W)(t_0,t_1)=\Gamma_W(0,t_0,t_1)=\gamma_W(t_0,1-t_0)=\gamma_W(t_0,t_1),\\
    &(d_1\Gamma_W)(t_0,t_1)=\Gamma_W(t_0,0,t_1)=\gamma_W(0,1)=d_1\gamma_W=W',\\
    &(d_2\Gamma_W)(t_0,t_1)=\Gamma_W(t_0,t_1,0)=\gamma_W(t_1,1-t_1)=\gamma_W(t_1,t_0)=\gamma_W^{-1}(t_0,t_1).
\end{align*}
Thus, $\Gamma_W\oplus s_0\gamma_V$ fills the top triangle. We will explain how to move systematically from these fillers, $\Gamma_W\oplus s_2\gamma_V$ and $\gamma$, to fillers of \eqref{9681FJ4}, but at this stage one may choose \emph{any} pasting without breaking functoriality, since we have identified these desired $2$-simplices as exactly those that are not required to satisfy any further conditions. This concludes the construction of $\mathbf{U}_{\leq2}\colon\overline{\EEx_{\leq2}}\to\Vhook$.

\subsection{The top face in arbitrary dimensions}

Before moving on to the induction step, we will construct $d_{0}^{\Bbox\OO}\mbf{U}(\gamma,j)$, the top face, with a closed-form formula (\Cref{2QKCQ8U}) in every dimension in a way that generalises the above constructions in dimensions $\leq 2$.

\begin{lemma}\label{Q52DBDN} Let  $(\gamma,j)\in\Phat_{k}\subset\EEx_{k+1}$ and write, as before, $\underline{0,i}\mapsto V_{0,i}=(\gamma,j)|_{i-1}$ for the edges of $\mbf{U}(\gamma,j)\in[\Path[k+2],B^{\oplus}\OO]$, and similarly $\underline{i,\ell}\mapsto V_{i,\ell}$, with $1\leq i<\ell\leq k+2$ throughout.\footnote{We need not assume that $\mbf{U}$ has been constructed, but may use instead \Cref{STG38ZQ} for these edges since the statement concerns only the restrictions of $\mbf{U}(\gamma,j)$ to $\Path[0,\alpha,\beta]\subset\Path[k+2]$.}
    \begin{enumerate}
    \item We have $V_{i,\ell}=\mfrak{N}(\gamma,j)|_{i-1}$ for $1\leq i\leq j<\ell \leq k+2$, and zero otherwise.
    \item\label{IQ9ITM0} Let $\alpha_{1}<\cdots<\alpha_{n}$ be a sequence of natural numbers within the interval $[1,k+2]$. Let $N=N_{\alpha}\in\{1,\dots,n\}$ be the smallest index such that $\alpha_{N}>j$ if it exists, and set $N=1$ otherwise. Then we have 
    \[
    V_{\alpha_{1},\dots,\alpha_{n}}=\begin{cases}
        V_{\alpha_{N-1},\alpha_{N}}=V_{\alpha_{N-1},k+2}\neq 0, & N>1,\\
        0, & N=1.
    \end{cases}
    \]
    Consequently, there are no non-trivial direct sums in the top face.
    \end{enumerate}
\end{lemma}
\begin{proof}
    First, note that every top edge $V_{\alpha\beta}$, $1\leq \alpha<\beta\leq k+2$, is within the fibre $B\OO(m)$ of the link projection, since it is the connecting edge in the restriction of $\mbf{U}(\gamma,j)$ to $\Path[0,\alpha,\beta]\subset\Path[k+2]$. It is the $2$-face 
    \[
        \begin{tikzcd}[sep=small]
        \alpha \ar[rr,"{V_{\alpha\beta}}"]
        &
        & 
        \beta
        \\
        &
        0 \ar[ul,"{V_{0,\alpha}}"]\ar[ur,"{V_{0,\beta}}"']
        \end{tikzcd}
    \]
    in $\Bbox\OO_{2}\cong\Vhook_{1}$ given by $(\gamma,j)|_{\alpha-1,\beta-1}\colon V_{0,\alpha}\to V_{0,\beta}$ underlied by $\gamma|_{\alpha-1,\beta-1}\colon V_{\alpha,\beta}\oplus V_{0,\alpha}\to V_{0,\beta}$. 
    Now, by construction, $(\gamma,j)|_{0,\dots,j-1}$ is low and $(\gamma,j)|_{j,\dots,n}$ is upper. Thus the restriction of $\mbf{U}(\gamma,j)$ to $\Path[0,i,\ell]$ for $\ell\leq j$ is wholly within $B\OO(n)$, so $V_{i,\ell}=0$. Similarly, it is wholly within $B\OO(n+m)$ for $\ell>i>j$, so then $V_{i,\ell}=0$ as well. Thus, $V_{i,\ell}\neq0$ implies $1\leq i\leq j<\ell \leq k+2$. 
    Conversely, if the inequalities are satisfied, then $V_{0,i-1}$ is low and $V_{0,\ell-1}$ is upper, so the connecting edge $V_{i,\ell}\in B\OO(m)$ is non-zero. Thus the un-equalities specify exactly the non-zero top edges. 
    
    Note that 
    \begin{equation}\label{0DF723Q}
        V_{i,\ell}=V_{i,\ell'}\ \text{if}\ j<\ell'\leq k+2\ \text{as well}, 
    \end{equation}
    since $\gamma|_{i-1,\ell-1}\colon V_{i,\ell}\oplus V_{0,i-1}\to V_{0,\ell}$ and $\gamma|_{i-1,\ell'-1}\colon V_{i,\ell'}\oplus V_{0,i-1}\to V_{0,\ell'}$ have the same source $\gamma|_{i-1}=V_{i,\ell}\oplus V_{0,i-1}=V_{i,\ell'}\oplus V_{0,i-1}\in  B\OO(n+m)_{0}$. In particular, the first statement is well-defined.
    That $V_{i,\ell}=\mfrak{N}(\gamma,j)|_{i-1}$ is immediate: $\mfrak{N}(\gamma,j)=\mrm{pr}_{2}(\gamma|_{1,\dots,j-1})$, so $\mfrak{N}(\gamma,j)|_{i-1}=\mrm{pr}_{2}(\gamma|_{i-1})$, so $\gamma|_{i,\ell}$ is of type $\gamma|_{i-1,\ell-1}\colon \mfrak{N}(\gamma,j)|_{i-1}\oplus V_{0,i}\to V_{0,\ell}$. 

    The second statement is a straightforward consequence. We have 
    \[
        V_{\alpha_{1},\dots,\alpha_{n}}=V_{\alpha_{n-1}\alpha_{n}}\oplus \cdots\oplus V_{\alpha_{N-1},\alpha_{N}}\oplus\cdots\oplus V_{\alpha_{1},\alpha_{2}}.
    \]
    By the above, the summands to either side of $V_{\alpha_{N-1},\alpha_{N}}$ are zero. The equality 
    \[
        V_{\alpha_{N-1},\alpha_{N}}=V_{\alpha_{N-1},k+2}
    \]
    follows from \eqref{0DF723Q} by setting $\ell'=k+2$. If there exists no $N$ as described, then every summand is zero since each $V_{0,\alpha_{k}}$ is in $B\OO(n)$. If $N=1$, then each $V_{0,\alpha_{k}}$ is in $B\OO(n+m)$ so that every summand is again zero.
\end{proof}

\begin{remark}
    The simplification noted in \Cref{Q52DBDN} is specific to our stratification depth being $1$ and not higher.  If it was higher, we would see nontrivial sums appearing in the top face as well.
\end{remark}

In order to proceed, we need the following fundamental fact about the simplicial category $\Path[n]$:

\begin{proposition}[{\cite[{\S}1.1.5]{lurie2009higher}; \cite[00LL]{kerodon}; \cite[00LM]{kerodon}}]\label{ECHTYNN}
    For $i,\ell\in[n]$, there is a canonical isomorphism
    \[
        \Hom_{\Path[n]}(i,\ell)\cong \left(\Delta[1]\right)^{\times(\ell-i-1)}
    \]
    of simplicial sets. Consequently, 
    there is a canonical homeomorphism
    \[
        \left|\Hom_{\Path[n]}(i,\ell)\right|\cong [0,1]^{\times(\ell-i-1)}.
    \]
\end{proposition}
\begin{proof}
    The first step is the construction of an isomorphism between $\left(\Delta[1]\right)^{\times n}$ and the nerve of the power poset $\mbf{P}(\{1,\dots,n\})$ of the set $\{1,\dots,n\}$, ordered by inclusion (and not reverse inclusion). This elementary observation was stated in \cite{kerodon} without proof, so, for completeness, we will provide one. 
    
    A vertex of $\left(\Delta[1]\right)^{\times n}$ is specified exactly by a function from $\{1,\dots,n\}$ to the $2$-element set (the vertices of $\Delta[1]$), which specifies exactly a subset of $\{1,\dots,n\}$. More generally, a $k$-simplex $\phi=(\phi_i)_{i=1}^{n}$ of $\left(\Delta[1]\right)^{\times n}$ is a collection of $n$ poset maps $\phi_{i}\colon [k]\to[1]$; 
    On the other hand, a $k$-simplex 
    \(
        \widehat{\alpha}=(\alpha^0\subseteq\cdots\subseteq\alpha^k)\in\Nerve_{k}\paren{\mbf{P}(\{1,\dots,n\})}
    \)
    is a non-decreasing sequence of subsets $\alpha^x\subseteq\{1,\dots,n\}$. Now, $\phi_i(x)\in[1]$ if and only if $i\in\{1,\dots,n\}$ belongs to $\alpha^x$. Thus, the function 
    \begin{align*}
        (\Delta[1])^{\times n}_k&\to \mbf{P}(\{1,\dots,n\})_k,\\
        \phi=(\phi_i)_{i=1}^{n}&\mapsto (\alpha^x_{\phi})_{x=0}^{k}
    \end{align*}
    where 
    \[
        \alpha^x_{\phi}=\left\{i\in\{1,\dots,n\} : \phi_i(x)=0\right\}
    \]
    is a bijection. It is well-defined since each $\phi_i$ is a poset map. Varying $k$, these functions are easily seen to assemble into an isomorphism of simplicial sets.
    The subsequent isomorphism 
    \[
        \Nerve_{\bullet}\paren{\mbf{P}(\{1,\dots,n\})}\cong\Nerve_{\bullet}(P_{0,n+1}^{\opp})
    \]
    is given by taking complements of subsets and thereby un-un-reversing the order. It is spelled out in \cite{kerodon}, as is the homeomorphism that is the second statement. Suffice it to say that a vertex $\beta\in\Nerve_0(P_{0,n+1}^\opp)$ is a subposet of $[n+1]$ of type $\beta=\underline{0,\beta',n+1}$ with $\beta'\subseteq\{1,\dots,n\}$, and so the rule 
    \begin{align*}
        \Nerve_0\paren{P_{0,n+1}^{\opp}}&\to\Nerve_{0}\paren{\mbf{P}}\\
        \beta&\mapsto \{1,\dots,n\}\smallsetminus\beta'
    \end{align*}
    defines a bijection. It is easily seen to lift to an isomorphism of simplicial sets.
\end{proof}

\begin{remark}\label{5D3UZA4}
    Because, in the proof of \Cref{ECHTYNN}, the value $0$ is (necessarily) taken to give the affirmative, the resulting cubes, when depicted in the standard way (mapping the vertices of $\paren{\Delta[1]}^{\times n}$ according to $\phi=(\phi_{i})_{i=1}^{n}\mapsto \sum_i \phi_{i}(0)\mbf{i}$ with $\mbf{i}\in\mbb{R}^n$ the $i$'th standard basis vector), will differ from those in \Cref{KSCTBU1,9ZWYRMA} below, but only up to a change of basis. For our purposes, this is not a problem: a precise choice of basis will be used only in the proof of \Cref{5KGP6JE} which is a convexity argument after which the choice may be reverted. Convexity is preserved under linear transformations.
\end{remark}

Occasionally, we will call the underlying posets $P_{i,\ell}^{(\opp)}$ `cubes' as well. When we do, \Cref{ECHTYNN} will be understood.

\begin{example}\label{KSCTBU1}
    \Cref{Q52DBDN} fully specifies the vertices of $P_{1,k+2}^\opp$ under $\mbf{U}$. The $k$-cube $P_{1,k+2}^{\opp}$ can be depicted for $k=3$ as follows:
    \[
        \begin{tikzcd}[sep=small]
        & \underline{12345} \arrow[dl] \arrow[rr] \arrow[dd] & & \underline{1235} \arrow[dl] \arrow[dd] \\
        \underline{1345} \arrow[rr, crossing over] \arrow[dd] & & \underline{135} \\
        & \underline{1245} \arrow[dl] \arrow[rr] & & \underline{125} \arrow[dl] \\
        \underline{145} \arrow[rr] & & \underline{15} \arrow[from=uu, crossing over]\\
    \end{tikzcd}
    \]
    If $j=1$, the image of this cube under $\mbf{U}(\gamma,1)$ must have all vertices equal to $V_{1,5}=\mfrak{N}_{0}\coloneqq\mfrak{N}(\gamma,j)|_{0}$, so non-degenerate sequences $S\in\Nerve_{3}(P_{1,5}^{\opp})$ of arrows from (the image of) $\underline{12345}$ to $\underline{15}$ must all be of type 
        \(
            \mfrak{N}_0\to\mfrak{N}_0\to\mfrak{N}_0\to\mfrak{N}_0.
        \)
        \item If $j=2$, the image under $\mbf{U}(\gamma,2)$ must be of type 
        \[
        \begin{tikzcd}[sep=small]
        & \mfrak{N}_{1} \arrow[dl] \arrow[rr] \arrow[dd] & & \mfrak{N}_{1} \arrow[dl] \arrow[dd] \\
        \mfrak{N}_{0} \arrow[rr, crossing over] \arrow[dd] & & \mfrak{N}_{0} \\
        & \mfrak{N}_{1} \arrow[dl] \arrow[rr] & & \mfrak{N}_{1} \arrow[dl] \\
        \mfrak{N}_{0} \arrow[rr] & & \mfrak{N}_{0} \arrow[from=uu, crossing over]\\
        \end{tikzcd}
        \]
        hence non-degenerate sequences as above must be of the following types: 
        \begin{gather*}
            \mfrak{N}_1\to\mfrak{N}_0\to\mfrak{N}_0\to\mfrak{N}_0\\
            \mfrak{N}_1\to\mfrak{N}_1\to\mfrak{N}_0\to\mfrak{N}_0\\
            \mfrak{N}_1\to\mfrak{N}_1\to\mfrak{N}_1\to\mfrak{N}_0
        \end{gather*}
        On the other hand $\Opp\mfrak{N}(\gamma,2)$ is a $1$-path of type $(\mfrak{N}_1\to\mfrak{N}_0)$.
        If $j=3$, we have 
        \[
        \begin{tikzcd}[sep=small]
        & \mfrak{N}_{2} \arrow[dl] \arrow[rr] \arrow[dd] & & \mfrak{N}_{2} \arrow[dl] \arrow[dd] \\
        \mfrak{N}_{2} \arrow[rr, crossing over] \arrow[dd] & & \mfrak{N}_{2} \\
        & \mfrak{N}_{1} \arrow[dl] \arrow[rr] & & \mfrak{N}_{1} \arrow[dl] \\
        \mfrak{N}_{0} \arrow[rr] & & \mfrak{N}_{0} \arrow[from=uu, crossing over]\\
        \end{tikzcd}
        \]
        so non-degenerate sequences as above must be of the following types: 
        \begin{gather*}
            \mfrak{N}_2\to\mfrak{N}_2\to\mfrak{N}_0\to\mfrak{N}_0\\
            \mfrak{N}_2\to\mfrak{N}_2\to\mfrak{N}_2\to\mfrak{N}_0\\
            \mfrak{N}_2\to\mfrak{N}_1\to\mfrak{N}_0\to\mfrak{N}_0\\
            \mfrak{N}_2\to\mfrak{N}_1\to\mfrak{N}_1\to\mfrak{N}_0\\
            \mfrak{N}_2\to\mfrak{N}_2\to\mfrak{N}_1\to\mfrak{N}_0
        \end{gather*}
        On the other hand, $\Opp\mfrak{N}(\gamma,3)$ is a $2$-path of type 
        \[
            \begin{tikzcd}[sep=small]
            &\mfrak{N}_1\ar[dl]\\ 
            \mfrak{N}_{0}&&\mfrak{N}_{2}\ar[ul]\ar[ll]\end{tikzcd}.
        \]
\end{example}

An examination of \Cref{KSCTBU1} is sufficient to reach the following Ansatz:

\begin{construction}\label{2QKCQ8U}
    Let $(\gamma,j)\in\Phat_{k}$ and let $n\geq 0$ be a natural number. To every sequence $\alpha=(\alpha^0>\cdots>\alpha^n)\in\Nerve_{n}(P_{i,j}^\opp)$ where $1\leq i\leq j<\ell\leq k+2$, we associate the map
    \begin{align*}
        A\colon[n]&\to[j-1],\\
        t&\mapsto j-1-({\alpha^t_{N_{t}-1}-1})
    \end{align*}
    in $\mbf{\Delta}$, where 
    \[
        N_{t}\coloneqq N_{\alpha^{t}}
    \] 
    is defined with respect to $\alpha^t=(\alpha_{1}^{t}<\cdots<\alpha^t_{n_{t}})$ as in \Cref{Q52DBDN}.
    Since $\Opp\mfrak{N}(\gamma,j)\in B\OO(m)_{j-1}$, we may pull it back along $A$ for any $\alpha\in\Nerve_{n}=\Nerve_{n}(P_{i,j}^\opp)$ to obtain, with a slight abuse of notation, a map 
    \begin{align*}
        \Opp\mfrak{N}(\gamma,j)\colon\Nerve_{n}&\to B\OO(m)_n,\\
        \alpha&\mapsto A^*\Opp\mfrak{N}(\gamma,j).
    \end{align*}
\end{construction}

\begin{lemma}\label{KXIJ7LO}
    For $1\leq i\leq j<\ell\leq k+2$, the map $\Opp\mfrak{N}(\gamma,j)\colon\Nerve_{\bullet}(P_{i,\ell}^{\opp})\to B\OO(m)_{\bullet}$ of \Cref{2QKCQ8U} is an $\infty$-functor.
\end{lemma}
\begin{proof}
    We must first prove that the map $A\colon [n]\to[j-1]$ of \Cref{2QKCQ8U} is well-defined and monotone. Since $N_t$ is the smallest index such that $\alpha^t_{N_{t}}>j$ (if it exists), we have $\alpha^{t}_{N_{t}-1}\leq j$, and so $A(t)=j-\alpha^{t}_{N_{t}-1}\geq0$. Moreover, $\alpha^{t}_{N_{t}-1}\geq1$ since $\alpha^{t}$ is a sequence of numbers that starts at $i\geq1$, hence $A(t)\leq j-1$. We can never have $N_t=1$ since $i\leq j<\ell$.
    As for monotonicity, we must show $\alpha^t_{N_{t}-1}\geq\alpha^{t'}_{N_{t'}-1}$ for $t\leq t'$ in $[n]$. That $\alpha^{t}>\alpha^{t'}$ means $\alpha^{t'}\subseteq\alpha^{t}$ as posets. We have $\alpha^{t}_{N_{t}-1}=\max\{x\in\alpha^{t} : x\leq j\}$ and since $\{x\in\alpha^{t'}: x\leq j\}\subseteq\{x\in\alpha^{t} : x\leq j\}$, we obtain 
    \[
        \alpha^{t'}_{N_{t'}-1}=\max\{x\in\alpha^{t'}: x\leq j\}\leq\max\{x\in\alpha^{t} : x\leq j\}=\alpha^{t}_{N_{t}-1}.
    \]
    This also shows $N_{t}\leq N_{t'}$: there can only be more indices $\ast$ in $\{1,\dots,n_{t'}\}$ where $\alpha^{t'}_\ast$ is at most $j$, not fewer.

    Finally, observe that the map $\phi\colon\Nerve_n\to\Hom_{\mbf{\Delta}}([n],[j-1])$, $\alpha\mapsto \phi(\alpha)\coloneqq A$ is manifestly simplicial, and therefore so is $\Opp\mfrak{N}(\gamma,j)$. Indeed, let $\delta\colon[n']\to[n]$ be a poset map. Then $(\delta^*\alpha)^t=\alpha^{\delta(t)}$ for $t\in[n']$, so 
    \[
        \phi(\delta^*\alpha)(t)=(\delta^*\alpha)^{t}_{N_{(\delta^*\alpha)^{t}}-1}-1=\alpha^{\delta(t)}_{N_{\alpha^{\delta(t)}}-1}-1=A(\delta(t))=\delta^*(\phi(\alpha))(t).\qedhere
    \]
\end{proof}

\begin{remark}\label{J8KZ1VE}
    \Cref{2QKCQ8U} immediately gives 
    \[
        \Opp\mfrak{N}(\gamma,j)|_{t}
        =\mfrak{N}(\gamma,j)|_{\alpha^{t}_{N_{t}-1}-1},
    \]
    since (cf.\ \cite[003M]{kerodon}) $\Opp$ reverses the operation $j-1-$. Consequently, 
    \[
        V_{\alpha}=\Opp\mfrak{N}(\gamma,j)|_{\alpha}
    \]
    for $\alpha=(\alpha_1<\cdots<\alpha_n)$ in the situation of \Cref{Q52DBDN}.
\end{remark}

\begin{lemma}\label{U0KWQ13}
    The maps of \Cref{KXIJ7LO} lift to a function 
    \[
        \Opp\mfrak{N}\colon \Phat_{k}\to\left[\Path[1,\dots,k+2],B^{\oplus}\OO\right].
    \]
\end{lemma}
\begin{proof}
    For pairs $i,\ell\in\{1,\dots,k+2\}$ that do \emph{not} satisfy $1\leq i\leq j<\ell\leq k+2$, we let $\Opp\mfrak{N}(\gamma,j)\colon\Nerve_\bullet(P_{i,\ell}^{\opp})\to B\OO_{\bullet}$ be the constant map to the zero vector space. This defines $\Opp\mfrak{N}(\gamma,j)$ on all morphism spaces. It remains to verify functoriality, which holds for trivial reasons: if $\alpha\colon i\to\ell$ and $\beta\colon\ell\to\ell'$, then either 
    \begin{itemize}
        \item $\Opp\mfrak{N}(\gamma,j)(\alpha)\neq 0$, in which case $\Opp\mfrak{N}(\gamma,j)(\beta)=0$ since $\ell,\ell'>j$, hence 
        \[
            \Opp\mfrak{N}(\gamma,j)(\beta\cup\alpha)=0\oplus\Opp\mfrak{N}(\gamma,j)(\alpha)=\Opp\mfrak{N}(\gamma,j)(\alpha),
        \]
        which holds since $N_{\alpha^*}=N_{(\beta\cup\alpha)^{*}}$ follows immediately from the definition. Appending $\beta$ to the head of $\alpha$ does not change the first index for which the sequence becomes larger than $j$.
        \item or $\Opp\mfrak{N}(\gamma,j)(\beta)\neq0$ in which case $\Opp\mfrak{N}(\gamma,j)(\alpha)=0$ since $\ell\leq j$, hence 
        \[
            \Opp\mfrak{N}(\gamma,j)(\beta\cup\alpha)=\Opp\mfrak{N}(\gamma,j)(\beta)\oplus 0=\Opp\mfrak{N}(\gamma,j)(\beta),
        \]
        which holds since $\beta_{N_{\beta}-1}=(\beta\cup\alpha)_{N_{\beta\cup\alpha}-1}$. Appending $\alpha$ to the foot of $\beta$ does not change the last element in the sequence before it grows larger than $j$, since that element is already within $\beta$.
        \item or both $\Opp\mfrak{N}(\gamma,j)(\alpha)$ and $\Opp\mfrak{N}(\gamma,j)(\beta)$ are zero, in which case it will suffice to show that $\Opp\mfrak{N}(\gamma,j)(\beta\cup\alpha)=0\oplus 0=0$ as well. This is clear if $i,\ell,\ell'\leq j$ or if $i,\ell,\ell'>j$. But $i,\ell\leq j$ and $\ell,\ell'>j$ cannot coincide, so these are all the cases.
    \end{itemize}
    This argument applies to unions of chains of posets to show functoriality on higher morphisms.
\end{proof}

\begin{proposition}\label{2FX6ABB}
    The functions 
    \[
        \Opp\mfrak{N}\colon\Phat_*\to\left[\Path[1,\dots,*+2],B^{\oplus}\OO\right]
    \] 
    of \Cref{U0KWQ13} extend 
    \[
        \mbf{U}_{\leq 1}\colon\overline{\EEx_{\leq1}}\to\Vhook.
    \]
\end{proposition}
\begin{proof}
    The only overlap in dimensions $\leq1$ is within $\Phat_0\subset\EEx_1$. Here, the equality of the two maps is trivial, but we include the proof for completeness. For $(\gamma,1)\in\Phat_0$ and so for $i=1$ and $\ell=2$, we have that $\mbf{U}(\gamma,1)|_{\Nerve_0(P_{1,2}^{\opp})}\colon\underline{12}\mapsto W=\mfrak{N}(\gamma,j)=\mrm{pr}_{2}(\gamma_0)$ is the normal component of $(\gamma,1)$. On the other hand, taking $\alpha=\underline{12}\in\Nerve_0(P_{1,2}^\opp)$ yields the constant $A=\id\colon[0]\to[0]$ and so $\Opp\mfrak{N}(\gamma,1)(\underline{12})=\id^*\Opp\mfrak{N}(\gamma,j)=\Opp(\mrm{pr}_{2}(\gamma_0))=\mrm{pr}_{2}(\gamma_0)$. The compatibility of the vertex values with any extension of $\mbf{U}_{\leq1}$ in higher dimensions was noted in \Cref{J8KZ1VE}.

    Now, for $k\geq 1$, let $S$ be a simplicial operator of type $S=\Sigma^*$ for 
    \[
        \Sigma\colon[k+1]\to[1]
    \] 
    in $\mbf{\Delta}$. For $(\gamma,1)\in\Phat_0$ again, assume that $S(\gamma,1)\in\Phat_{k}\subset\EEx_{k+1}$ is also vertical. This amounts to assuming that $\Sigma$ is surjective. In this situation, $\Opp\mfrak{N}$ and $\mbf{U}_{\leq1}$ may be compared, and we must show that
    \begin{equation}\label{3U6CABJ}
        \mbf{U}_{\leq1}(S(\gamma,1))|_{\Path[1,\dots,k+2]}\overset{\text{def}}{=}S\mbf{U}(\gamma,1)|_{\Path[1,\dots,k+2]}=\Opp\mfrak{N}(S(\gamma,1))
    \end{equation}
    where the LHS is the restriction to $\Path[1,\dots,k+2]\subset\Path[k+2]$ of 
    \[
        S\mbf{U}(\gamma,1)\colon\Path[k+2]\to B^{\oplus}\OO.
    \] 
    
    We have $\mfrak{N}(S(\gamma,1))=\mfrak{N}(S\gamma,\sharp^S1)=\mrm{pr}_{2}(S\gamma|_{0,\dots,\sharp^{S}1})$, where $\sharp^{S}$ applies $\sharp$'s and $\flat$'s to the exit index $1$ according to $S$. The underlying simplex map $\Sigma$ is determined by a unique index 
    \[
        e_{\Sigma}\in\{1,\dots,k+1\}
    \]
    such that $\Sigma(x)=0$ for $x\leq e_{\Sigma}-1$ and $\Sigma(x)=1$ for $x\geq e_{\Sigma}$. It is straightforward to see that 
    \[
        \sharp^S1=e_{\Sigma},
    \]
    so 
    \[
        \Opp\mfrak{N}(S(\gamma,1))=\Opp\,\mrm{pr}_{2}(S\gamma|_{0,\dots,e_{\Sigma}-1})=\mrm{pr}_{2}(\Opp(S\gamma|_{0,\dots,e_{\Sigma}-1}))\in B\OO(m)_{e_{\Sigma}-1}.
    \]
    But
    \begin{align*}
        (S\gamma)|_{0,\dots,e_{\Sigma}-1}&=([e_{\Sigma}-1]\hookrightarrow[k+1]\xrightarrow{\Sigma}[1])^{*}(\gamma)\\
        &=([e_{\Sigma}-1]\to[0]\hookrightarrow[1])^{*}\gamma\\
        &=(s_{0})^{e_{\Sigma}-1}d_{1}\gamma\\
        &=(s_0)^{e_{\Sigma}-1}\gamma_{0},
    \end{align*}
    and $\Opp s_{0}=s_0 \Opp$ in dimension $0$, and $\mrm{pr}_{2}$ commutes with $\Opp$ as well as with any simplicial operator, so we have 
    \begin{align*}
        \Opp\mfrak{N}(S(\gamma,1))\colon \alpha\mapsto &A^*\Opp\,\mrm{pr}_{2}(S\gamma|_{0,\dots,e_{\Sigma}-1})
        =\mrm{pr}_{2}A^*\Opp(s_0)^{e_{\Sigma}-1}\gamma_0\\
        &=\mrm{pr}_{2}A^*(s_{0})^{e_{\Sigma}-1}\Opp\gamma_0
        =\mrm{pr}_{2}A^*(s_{0})^{e_{\Sigma}-1}\gamma_0\\
        &=\mrm{pr}_{2}(s_{0})^{n}\gamma_0
        =(s_{0})^{n}\mrm{pr}_{2}\gamma_0
    \end{align*}
    for $\alpha\in\Nerve_{n}(P_{i,\ell}^\opp)$ and $1\leq i\leq e_{\Sigma}<\ell\leq k+2$, and (the $n$-fold degenerate) zero otherwise.

    On the other hand, $\Sigma\colon[k+1]\to[1]$ induces a map 
    \[
        \Sigma_{+1}\colon[k+2]\to[2]
    \]
    defined by 
    \[
        \Sigma_{+1}(0)=0,\quad \text{and}\quad  \Sigma_{+1}(i)=\Sigma(i-1)+1
    \]
    for $1\leq i\leq k+2$. Along the under-$\infty$-category identifications $\Vhook_{*}\cong (\Bbox\OO)_{*+1}$, we have that $S(\mbf{U}(\gamma,1))\in\Vhook_{k+1}$ corresponds to 
    \(
        S_{+1}(\mbf{U}(\gamma,1))\coloneqq\Sigma_{+1}^*(\mbf{U}(\gamma,1))
    \)
    in $(\Bbox\OO)_{k+2}$, where we identified $\mbf{U}(\gamma,1)$ with the corresponding element \(
        \mbf{U}(\gamma,1)\colon\Path[2]\to B^{\oplus}\OO\) in  \((\Bbox\OO)_{2}.
    \)
    The restriction of the induced map $\Sigma_{+1}\colon\Path[k+2]\to\Path[2]$ to $\Path[1,\dots,k+1]$ factors through $\Path[1,2]\subset\Path[2]$ by construction and defines the LHS of \eqref{3U6CABJ} as the composition
    \[
        \Path[1,\dots,k+2]\xrightarrow{\Sigma_{+1}|}\Path[1,2]\xrightarrow{\mbf{U}(\gamma,1)|}B^{\oplus}\OO.
    \]
    Let now $e_{\Sigma}$ be as above, so that $e_{\Sigma}+1\in\{2,\dots,k+2\}$ fulfills the analogous function for the restriction of $\Sigma_{+1}$: We have $\Sigma_{+1}(x)=1$ for $1\leq x\leq e_{\Sigma}=e_{\Sigma}+1-1$ and $\Sigma(x)=2$ for $x\geq e_{\Sigma}+1$. Suppose, then, that $1\leq i<\ell\leq k+2$. We have $S\mbf{U}(\gamma,1)\colon\underline{i,\ell}\mapsto \mbf{U}(\gamma,1)(\id_{1})=\id_{B^{\oplus}\OO}=0\in B\OO$ if $\Sigma_{+1}(i)=\Sigma_{+1}(\ell)=1$, i.e., if $i,\ell\leq e_{\Sigma}$. Similarly, $S\mbf{U}(\gamma,1)\colon\underline{i,\ell}\mapsto\mbf{U}(\gamma,1)(\id_{2})=0$ if $\Sigma_{+1}(i)=\Sigma_{+1}(\ell)=2$, i.e., if $i,\ell>e_{\Sigma}$. If, however, $i\leq e_{\Sigma}$ and $\ell>e_{\Sigma}$, then $S\mbf{U}(\gamma,1)\colon\underline{i,\ell}\mapsto\mbf{U}(\gamma,1)(\underline{12})=W=\mrm{pr}_{2}(\gamma_{0})$. In both cases, the result coincides with the value of the RHS of \eqref{3U6CABJ}. 

    Generalising this observation to higher dimensions is straightforward: let $\alpha\in \Nerve_{n}(P_{i,\ell}^{\opp})$. If $i,\ell\leq e_{\Sigma}$ or $i,\ell>e_{\Sigma}$, then $S\mbf{U}(\gamma,1)\colon\alpha\mapsto (s_{0})^n0$, and if $i\leq e_{\Sigma}$ and $\ell>e_{\Sigma}$, then $S\mbf{U}(\gamma,1)\colon\alpha\mapsto (s_{0})^n\mrm{pr}_{2}\gamma_0$.
\end{proof}

\subsection{The induction step}\label{EEYT903}

Let us write 
\[
    \overline{\EEx_{\leq k+1}}\coloneqq\overline{(\EEx_{\leq k+1}\smallsetminus\overline{\EEx_{\leq k}})\cup\overline{\EEx_{\leq k}}}
\]
for the simplicial subset of $\EEx$ generated by $\EEx_{\leq k}$ together with the non-degenerate  exit $(k+1)$-paths in $\Phat_{k}\subset\EEx_{k+1}$ and let us assume a map 
\[
    \mbf{U}_{\leq k}\colon \overline{\EEx_{\leq k}}\to\Vhook
\]
is given which satisfies 
\begin{equation}\label{81ITXUL}
    \mbf{U}_{\leq k}|_{\leq 1}=\mbf{U}_{\leq 1}
\end{equation}
and
\begin{equation}\label{PTBY5XA}
    d_{0}^{\Bbox\OO}\mbf{U}_{\leq k}(\gamma,j)\overset{\text{def}}{=}\mbf{U}_{\leq k}(\gamma,j)|_{\Path[1,\dots,*+2]}=\Opp\mfrak{N}(\gamma,j)
\end{equation}
for all $(\gamma,j)\in\Phat_{*}\subset(\overline{\EEx_{\leq k+1}})_{*+1}$, with $\Opp\mfrak{N}$ defined as in \Cref{2QKCQ8U}. \Cref{2FX6ABB} states that this equality holds in the base case $k=1$. Moreover, this is consistent by \Cref{Q52DBDN}: as noted there, the condition that $\mbf{U}_{\leq k}$ extend $\mbf{U}_{\leq 1}$ fixes the spaces $V_{i,\ell}$ for $1\leq i<\ell\leq k+2$ through the induced restrictions to the subcategories $\Path[0,\alpha,\beta]\subset\Path[k+2]$, and these spaces coincide with the values of $\Opp\mfrak{N}$ wherever they overlap.

There is a final inductive assumption which we will formulate and justify in the proof of \Cref{VMR7DN8}: see immediately after \eqref{8RDCLNK}. It involves a construction that becomes necessary for the first time within said proof, which is why we have chosen to formulate it there.

Now, having fixed $\mbf{U}_{\leq k}(\gamma,j)$ on every $\Bbox\OO$-face, i.e., on every proper simplicial subcategory of $\Path[k+2]$, we have fixed it on every hom-space $\Hom_{\Path[k+2]}(i,\ell)$ for $(i,\ell)\neq (0,k+2)$. 
It remains then to specify it on all of $\Hom_{\Path[k+2]}(0,k+2)=\Nerve_\bullet P_{0,k+2}^{\opp}$. 
The initial object of $P_{0,k+2}^{\opp}$ is $[k+2]=\underline{0,\dots,k+2}$, and its final object is $\underline{0,k+2}$. All arrows with domain $[k+2]$ are composite except for the arrow $[k+2]>\underline{0,k+2}$, so this is the only $1$-morphism of $\Path[k+2]$ whose image is not determined by the inductive hypothesis. For instance, $\underline{01234}>\underline{014}=(\underline{01}>\underline{01})\cup(\underline{1234}>\underline{14})$, and the value of both factors is determined already by $\mbf{U}_{\leq k}$. Generally, any $1$-morphism $[k+2]>[k+2]\smallsetminus \beta$ with $\beta\subset\{1,\dots,k+1\}$ a proper subset is determined by functoriality in the same way, as can be seen by decomposing the target nontrivially and using the initiality of both factors in the corresponding decomposition of $[k+2]$. Namely, take $\delta\in\{1,\dots,k+1\}\smallsetminus\beta$, and write 
\[
    [k+2]\smallsetminus\beta=(\underline{0,\dots,\delta}\cap([k+2]\smallsetminus\beta))\cup(\underline{\delta,\dots,k+2}\cap([k+2]\smallsetminus\beta))=\beta_L^{\mrm{c}}\cup\beta_R^{\mrm{c}}.
\]
Then 
\[
    [k+2]\smallsetminus\beta>\underline{0,k+2}=(\underline{0,\dots,\delta}>\beta_L^{\mrm{c}})\cup(\underline{\delta,\dots,k+2}>\beta_R^{\mrm{c}}).
\]

Consequently, the images of all higher morphisms of whom $[k+2]>\underline{0,k+2}$ is an edge are likewise undetermined. Since we are mapping from the nerve of $P_{0,k+2}^\opp$, this means (the images of the) sequences with long edge $[k+2]>\underline{0,k+2}$. By the \emph{long edge} of an $n$-chain $S\in\Nerve_n$, $n\geq1$, we mean its pullback along $[1]\hookrightarrow[n]$, $0\mapsto 0$; $1\mapsto n$. Clearly, if $[k+2]>\underline{0,k+2}$ is an edge of $S$, then it is also its long edge.

It will therefore suffice to provide construct images for the simplices of $\Nerve_\bullet(P_{0,k+2}^\opp)$ with long edge $[k+2]>\underline{0,k+2}$, as well as an image for this $1$-morphism itself. We have seen that the image of the latter need obey no other condition. Our strategy will be to generalise the idea of the decomposition in \eqref{1HM4M30}. To this end, we will first identify where exactly the path $\gamma\in B\OO(n+m)_{k+1}$ lies in the image of the $(k+1)$-cube $P_{0,k+2}$. 
Let us write $\mfrak{N}_*\coloneqq \mfrak{N}(\gamma,j)|_{*}$.

\begin{lemma}\label{3077OTN}
    In the situation of \Cref{Q52DBDN}, let $\alpha=(\alpha_{1}<\cdots<\alpha_{n})$ be a sequence within $[1,k+2]$. Then we have
    \[
        \mbf{U}_{(\gamma,j)}
        (\underline{0,\alpha})=
        \begin{cases}
        \mfrak{N}_{\alpha_{N-1}-1}\oplus (\gamma,j)|_{\alpha_1-1}, & N>1,\\
        (\gamma,j)|_{\alpha_{1}-1}, & N=1.
        \end{cases}
    \]
\end{lemma}
\begin{proof}
    This follows immediately from \Cref{Q52DBDN} and \Cref{2FX6ABB} after decomposing the space as $V_{0,\alpha}=V_{\alpha}\oplus V_{0,\alpha_1}$.
\end{proof}

\begin{corollary}\label{OXY4ORA}
    In the situation of \Cref{Q52DBDN}, for $0<i<k+2$, we have 
    \(
        \mbf{U}_{(\gamma,j)}(\underline{0,i,k+2})=\gamma|_{i-1},
    \)
    and 
    \(
        \mbf{U}_{(\gamma,j)}(\underline{0,k+2})=\gamma|_{k+1}.
    \)
\end{corollary}
\begin{proof}
    The second equality is clear. Now, \Cref{3077OTN} yields 
    \[
        \mbf{U}_{(\gamma,j)}(\underline{0,i,k+2})=
        \begin{cases}
        \mfrak{N}_{i-1}\oplus(\gamma,j)|_{i-1}, & i\leq j\\
        (\gamma,j)|_{i-1}, & i>j.
        \end{cases}
    \]
    If $i\leq j$, then $(\gamma,j)|_{i-1}$ is low, so $\gamma|_{i-1}=\mfrak{N}_{i-1}\oplus\pi\paren{\gamma|_{i-1}}$, and if $i>j$, then $(\gamma,j)|_{i-1}$ is upper, so $\gamma|_{i-1}=(\gamma,j)|_{i-1}$ (recall that we suppress $\iota$). 
\end{proof}

All of $\Nerve_\bullet(P_{0,k+2}^\opp)$ is to map to $B\OO(n+m)=\mrm{Sing}_{\bullet}B\OO(n+m)$ under $\mbf{U}(\gamma,j)$. Therefore, by the adjunction between $|-|$ and $\mrm{Sing}_{\bullet}$, this is equivalent to mapping out of $|\Nerve_{\bullet}|$ to the space $B\OO(n+m)$ instead. We will specify the location of $\gamma$ inside (the image of) $\Nerve_{\bullet}$ by means of an embedding $\Delta^{k+1}\hookrightarrow|\Nerve_{\bullet}|$.

\begin{construction}\label{HZBV0LB}
    The vertices featuring in \Cref{OXY4ORA} specify a (topological) $(k+1)$-simplex within $|\Nerve_{\bullet}|$. Namely, since, by \Cref{ECHTYNN}, $|\Nerve_{\bullet}|$ is canonically homeomorphic to the $(k+1)$-cube, we may define 
    \[
        \btd\colon\Delta^{k+1}\hookrightarrow |\Nerve_{\bullet}|
    \]
    to be (i.e., map homeomorphically onto) the subset of $|\Nerve_{\bullet}|$ given by the convex hull of the vertices 
    \[
        \{\underline{0,i,k+2} : 0<i<k+2\}\cup\{\underline{0,k+2}\}\subset|\Nerve_{\bullet}|\cong[0,1]^{\times(k+1)}.
    \]
    This is indeed a topological $(k+1)$-simplex: shifting and rotating the cube so that the (image under this homeomorphism of) $\underline{0,k+2}$ lies at the origin, it is easily checked (vis-\`a-vis the proof of \Cref{ECHTYNN}, keeping in mind \Cref{5D3UZA4}) that the points (given by the images of) $\underline{0,i,k+2}$ give exactly the standard unit vectors.

    Writing 
    \[
        \bigtriangleup\coloneqq\overline{|\Nerve_{\bullet}|\smallsetminus\btd}
    \]
    for the closure within $[0,1]^{\times(k+1)}$ of the complement, we obtain the decomposition 
    \begin{equation}\label{L45NFI3}
        |\Nerve_{\bullet}|\cong \bigtriangleup\cup_{\partial}\btd
    \end{equation}
    where the common boundary $\partial$ is the convex hull of the vertices $\{\underline{0,i,k+2} : 0<i<k+2\}$. 
\end{construction}

\begin{example}\label{9ZWYRMA}
    Consider an exit $3$-path $(\gamma,3)\in\Phat_2$ of index $3$:
    \begin{equation*}
        \begin{tikzcd}[sep=small]
            &   & K\\
            \\
            & W_1\oplus V_1\ar[uur,dotted]\ar[dr,dotted]   & \\
        W_0\oplus V_0\ar[uuurr,bend left]\ar[ur,dotted]\ar[rr]    && W_2\oplus V_2\ar[uuu]
        \end{tikzcd}
    \end{equation*}
    The diagram depicts the underlying path $\gamma\in B\OO(n+m)_3$. 
    The edges of $\mathbf{U}(\gamma,3)\in \Bbox\OO_4$ are given, due to \eqref{81ITXUL}, as follows: 
    \begin{equation*}
        \begin{cases}
            \underline{01}\mapsto V_0,\ \underline{02}\mapsto V_1,\ \underline{03}\mapsto V_2,\ \underline{04}\mapsto K,\\
            \underline{14}\mapsto W_0,\ \underline{24}\mapsto W_1,\ \underline{34}\mapsto W_2
        \end{cases}
    \end{equation*}
    and the remaining edges are zero. Now, the image of the $3$-cube $P_{0,4}^\opp$
    \begin{equation*}
        \begin{tikzcd}[sep=small]
            & \underline{01234} \arrow[dl] \arrow[rr] \arrow[dd] & & \underline{0124} \arrow[dl] \arrow[dd] \\
            \underline{0234} \arrow[rr, crossing over] \arrow[dd] & & \underline{024} \\
            & \underline{0134} \arrow[dl] \arrow[rr] & & \underline{014} \arrow[dl] \\
            \underline{034} \arrow[rr] & & \underline{04} \arrow[from=uu, crossing over]\\
        \end{tikzcd}
    \end{equation*}
    under $\mathbf{U}(\gamma,j)$ looks as follows:
    \begin{equation*}
        \begin{tikzcd}
            & W_2\oplus V_0 \arrow[dl,"{\id\oplus\gamma_{V01}}"description] \arrow[rr,"{\gamma_{W12}^{-1}\oplus\id}"description] \arrow[dd,equal] & & W_1\oplus V_0 \arrow[dl,"{\id\oplus\gamma_{V01}}"description] \arrow[dd,"{\gamma_{W01}^{-1}\oplus\id}"description] \\
            W_2\oplus V_1 \arrow[rr,crossing over,"{\id_{W12}^{-1}\oplus\id}"description] \arrow[dd,"{\id\oplus\gamma_{V12}}"description] & & W_1\oplus V_1  \\
            & W_2\oplus V_0 \arrow[dl,"{\id\oplus\gamma_{V02}}"description] \arrow[rr,"{\gamma_{W02}^{-1}\oplus\id}"description] & & W_0\oplus V_0 \arrow[dl,"{\gamma_{\oplus0}}"description] \\
            W_2\oplus V_2 \arrow[rr,"{\gamma_{\oplus2}}"description] & & K \arrow[from=uu, crossing over,"{\id_{\oplus1}}"description,near end]\\
        \end{tikzcd}
    \end{equation*}
    Painting in the outer-left concatenation chosen to be the image of $\underline{01234}>\underline{04}$ according to the discussion in \Cref{0N49067} in green and the edges of the `lower' tetrahedron given by $\gamma\in B\OO_3$ itself in blue, we see that 
    \begin{equation*}
        \begin{tikzcd}[]
        & W_2\oplus V_0 \arrow[dl,color=ForestGreen] \arrow[rr] \arrow[dd,equal] & & W_1\oplus V_0 \arrow[dl] \arrow[dd] \\
        W_2\oplus V_1 \arrow[rr, crossing over] \arrow[dd,color=ForestGreen] & & W_1\oplus V_1 \\
        & W_2\oplus V_0 \arrow[dl] \arrow[rr] & & W_0\oplus V_0\ar[ul,color=blue,"{\gamma_{W01}\oplus\gamma_{V01}}"description],\ar[dlll,color=blue,"{\gamma_{W02}\oplus\gamma_{V02}}"description] \arrow[dl,color=blue] \\
        W_2\oplus V_2 \arrow[rr,color=cyan]\ar[from=uurr,crossing over,color=blue,bend left=10,"{\gamma_{W12}\oplus\gamma_{V12}}"description,near start] & & K \arrow[from=uu, crossing over,color=blue]\\
        \end{tikzcd}
    \end{equation*}
    homotopy-commutes by inspection. In terms of \Cref{HZBV0LB}, the $3$-simplex with the blue/cyan edges geometrically-realises to $\btd$, and the cube with $\btd$ cut off gives $\btu$. The $2$-dimensional instantiation of this idea was depicted in \eqref{1HM4M30}.
\end{example}

As \Cref{9ZWYRMA} suggests, the next step in our strategy is to define a new poset $\overline{P^{\opp}_{0,\widehat{k+2}}}$ such that the full subposet $P^{\opp}_{0,k+2}\smallsetminus\{\underline{0,k+2}\}$ of the original poset given by removing its final object is embedded into it, 
\begin{equation}\label{SL2GYB3}
    P^{\opp}_{0,k+2}\smallsetminus\{\underline{0,k+2}\}\subset \overline{P^{\opp}_{0,\widehat{k+2}}},
\end{equation}
and such that we have a homeomorphism
\begin{equation}\label{IZSUMZ0}
    \left|\Nerve_{\bullet}\left(\overline{P_{0,\widehat{k+2}}^\opp}\right)\right|\cong\btu.
\end{equation}

\begin{definition}\label{NTR664O}
    By $\overline{P^{\opp}_{0,\widehat{k+2}}}$ we denote the poset whose objects are the same as those of $P_{0,k+2}^{\opp}\smallsetminus\{\underline{0,k+2}\}$, and whose arrows are those of the latter together with the new primitive arrows 
    \[
        \underline{0,i,k+2}\to\underline{0,\ell,k+2}
    \]
    whenever $0<i<\ell<k+2$. 
\end{definition}

\begin{lemma}\label{5KGP6JE}
    \Cref{NTR664O} satisfies \eqref{SL2GYB3} and \eqref{IZSUMZ0}.
\end{lemma}
\begin{proof}
That \eqref{SL2GYB3} is satisfied is clear by the construction. Let us observe now that \eqref{SL2GYB3} lifts as in the diagram 
\[
    \begin{tikzcd}
    \left|\Nerve_{\bullet}\left(P_{0,k+2}^{\opp}\smallsetminus\{\underline{0,k+2}\}\right)\right|\ar[d,hook]\ar[r,hook] & \left|\Nerve_{\bullet}\left(\overline{P_{0,\widehat{k+2}}^\opp}\right)\right|\ar[dl,hook,dotted]\\
    {[0,1]}^{\times(k+1)}
    \end{tikzcd}
\]
to an embedding into the $(k+1)$-cube. Indeed, any topological simplex in $\left|\Nerve_{\bullet}\paren{\overline{P^\opp_{0,\widehat{k+2}}}}\right|$ can be sent to the convex hull of the images of its vertices within $[0,1]^{\times(k+1)}$. This makes the diagram above commute, since the map 
\(
    \left|\Nerve_{\bullet}\paren{P_{0,k+2}^\opp\smallsetminus\{\underline{0,k+2}\}}\right|\hookrightarrow[0,1]^{\times(k+1)}
\)
is given by restricting $\left|\Nerve_{\bullet}\paren{P_{0,k+2}}\right|\hookrightarrow[0,1]^{\times(k+1)}$, and the latter is easily seen to be itself defined by sending simplices to the convex hulls of the images of their vertices. The images of the vertices are fixed explicitly in the proof of \Cref{ECHTYNN}. 

Note now that the intersection of the image $\widetilde{\btu}$ of $\left|\Nerve_{\bullet}\left(\overline{P_{0,\widehat{k+2}}^\opp}\right)\right|$ and $\btd$ is given exactly by the convex hull of $\{\underline{0,i,k+2} : 0<i<k+2\}$, which is by construction the image of the geometric realisation of the $k$-simplex 
\[
    (\underline{0,1,k+2}\to\underline{0,2,k+2}\to\cdots\to\underline{0,k+1,k+2})
\]
from $\Nerve_{\bullet}\paren{\overline{P_{0,\widehat{k+2}}^{\opp}}}$. Therefore, we can glue with $\btd$ and stay within the cube: 
\[
    \widetilde{\btu}\cup_{\partial}\btd\subseteq [0,1]^{\times(k+1)}.
\]
It remains to show the reverse inclusion, which will imply that $\widetilde{\btu}$, and therefore $\left|\Nerve_{\bullet}\left(\overline{P_{0,\widehat{k+2}}^\opp}\right)\right|$, are homeomorphic to $\btu$.

The $(k+1)$-cube itself being the convex hull of its corners, it is equal to the convex hull of the union of $\widetilde{\btu}$ and $\btd$. It is the union of all convex combinations of these two sets:\footnote{See e.g.\ Rockafellar \cite[Theorem 3.3]{rockafellar1970convex}. It is {not} true in general that the convex hull of a collection of points in euclidean space is the union of the convex hulls of two subsets that partition the collection with non-empty intersection.}
\[
    [0,1]^{\times(k+1)}=\bigcup_{\substack{\lambda_i\geq 0,\\ \lambda_1+\lambda_2=1}}\left(\lambda_1\widetilde{\btu}+\lambda_2\btd\right)
\]
This means that it suffices to show that given points $x\in\widetilde{\btu}$ and $y\in\btd$, the line segment $L$ connecting $x$ and $y$ lies within $\widetilde{\btu}\cup\btd$. 

This is a consequence of the intermediate value theorem. Shifting and rotating the cube such that the image $\underline{0,k+2}$ is at the origin and the image $\mbf{i}=(0,\dots,0,1,0,\dots,0)\in\mbb{R}^{k+1}$ of $\underline{0,i,k+2}$ is the $i$'th unit vector among the $k+1$ coordinate axes (cf.\ the discussion in \Cref{HZBV0LB}), we place $L$ as well as all of the $(k+1)$-cube inside the non-negative orthant of $\mbb{R}^{n+1}$. We may assume that neither point lies on the common boundary 
\(
    \partial=\left\{\sum_{i=1}^{k+1}\lambda_i\mbf{i} : \sum\lambda_i=1\right\},
\)
since otherwise $L$ lies within (at least) one of the two sets by their convexity, and we are done. Now, $\btd$ is the convex hull of $\{0,\mbf{i}\}_{i=1}^{k+1}$, so 
\(
    y=\sum_{i=1}^{k+1}\lambda_i\mbf{i}
\) 
with
\(
    \sum\lambda_{i}<1.
\)
On the other hand, $\widetilde{\btu}$ is the convex hull of 
\(
    \left\{\sum_{i=1}^{k+1}\mu_{i}\mbf{i} : \mu_i\in\{0,1\},\ \vec{\mu}\neq0\right\},
\)
with $2^{k+1}-1$ possibilities for $\mu=(\mu_i)_{i=1}^{k+1}\in \{0,1\}^{\{1,\dots,k+1\}}\smallsetminus \{\mrm{const}_0\}$. We obtain 
\(
    x=\sum_{\substack{j\in\{1,\dots,2^{k+1}-1\}\\ \sum\lambda_j=1}}\lambda_{j}\sum_{\substack{i\in\{1,\dots,k+1\}\\ \mu^{j}_{i}\in\{0,1\}\\ \sum_i\mu^{j}_{i}\geq1}}\mu_{i}^{j}\mbf{i}
\)
and, writing 
\(
    \lambda'_{i}\coloneqq\sum_{\substack{j\in\{1,\dots,2^{k+1}-1\} \\ \sum\lambda_{j}=1}}\lambda_{j}\mu^{j}_{i},
\)
we have 
\(
    x=\sum_{i=1}^{k+1}\lambda'_{i}\mbf{i}
\)
with 
\(
    \sum_{i}\lambda'_{i}\neq1
\)
due to $x\notin\partial$ by assumption. But since $\sum_i\mu_i^{j}\geq1$ and $\sum\lambda_j=1$, we have $\sum\lambda_i'\geq1$ in any case, so we must have
\(
    \sum\lambda'_i>1.
\)
Thus, the continuous function 
$
    L\to\mbb{R}$, 
$    \sum\nu_{i}\mbf{i}\mapsto1-\sum\nu_i$
is negative at $x$ and positive at $y$, and so is zero at some $Z\in L$. But then $Z\in\partial$, so, writing
\[
    L=L_{yZ}\cup_Z L_{Zx}
\]
where $L_{\alpha\beta}$ is the line segment connecting the points $\alpha$ and $\beta$, we have $L_{yZ}\subseteq\btd$ and $L_{Zx}\subseteq\widetilde{\btu}$ by convexity, since $Z\in\widetilde{\btu}\cap\btd$. Therefore $L\subseteq\widetilde{\btu}\cup\btd$, which implies $[0,1]^{\times(k+1)}\subseteq\widetilde{\btu}\cup\btd$.
\end{proof}

Let us summarise the resulting strategy of proof for \Cref{YUXW1VO}:

\begin{lemma}\label{YGELJPS}
    Providing, for each $(\gamma,j)\in\Phat_{k}\subset\EEx_{k+1}\smallsetminus\overline{\EEx_{\leq k}}$, an $\infty$-functor 
    \[
        \Nerve_{\bullet}\paren{\overline{P^{\opp}_{0,\widehat{k+2}}}}\to B\OO(n+m),
    \]
    whose 
    \begin{itemize}
        \item restriction to $\Nerve_{\bullet}\paren{P_{0,k+2}}\smallsetminus\{{\underline{0,k+2}}\}$ agrees with the restriction of $\mbf{U}_{\leq k}$, and whose
        \item value on 
        \[
            \paren{\underline{0,1,k+2}\to\underline{0,2,k+2}\to\cdots\to\underline{0,k+1,k+2}}\in\Nerve_{k}\paren{\overline{P_{0,\widehat{k+2}}}}
        \]
        is 
        \(
            d_{k+1}(\gamma)\in B\OO(n+m)_{k}
        \)
    \end{itemize}
    yields an extension 
    \(
        \mbf{U}_{\leq k+1}\colon\overline{\EEx_{k+1}}\to\Vhook
    \)
    of $\mbf{U}_{\leq k}$.
\end{lemma}
\begin{proof}
    First, let us summarise what we have already proved. 
    Given such an exit path $(\gamma,j)$, it suffices, by the discussion at the beginning of this section, to provide an extension of $\mbf{U}_{\leq k}$ on $\Nerve_{\bullet}\paren{P_{0,k+2}^{\opp}}$. Using the adjunction between geometric realisation and the singular chains functor and combining \Cref{5KGP6JE} with \eqref{L45NFI3}, we see that it suffices to provide two maps of type
    \begin{align*}
        &\left|\Nerve_{\bullet}\paren{\overline{P^{\opp}_{0,\widehat{k+2}}}}\right|\to B\OO(n+m)\\
        &\btd\cong\Delta^{k+1}\to B\OO(n+m)
    \end{align*}
    such that 
    \begin{itemize}
        \item they agree on $\partial\cong\Delta^k$, the convex hull of $\{\underline{0,i,k+2} : 0<i<k+2\}$; and
        \item they extend $\mbf{U}_{\leq k}$.
    \end{itemize}

    Now, we can define $\btd$ to map to $\gamma\in B\OO(n+m)_{k+1}$. Since $\btd\cong\Delta^{k+1}$ is given by identifying $\underline{0,i,k+2}$ with $i-1$ and $\underline{0,k+2}$ with $k+1$, we see that the sequence $\paren{\underline{0,1,k+2}\to\cdots\to\underline{0,k+1,k+2}}\in\Nerve_{k}\paren{\overline{P_{0,\widehat{k+2}}}}$ is exactly the $(k+1)$'st face, so its value is $d_{k+1}\gamma$ by construction. 
    
    The statement will follow once we show that $\gamma\colon\btd\to B\OO(n+m)$ itself is compatible with $\mbf{U}_{\leq k}$.
    The only definitional overlap is at the vertices and those edges that are of type $\underline{0,i,k+2}>\underline{0,k+2}$, since the other edges $\underline{0,i,k+2}\to\underline{0,\ell,k+2}$, $i<\ell$, are not in $P^{\opp}_{0,k+2}$. This verifies that there are no overlapping higher simplices. Now, \Cref{OXY4ORA} states exactly that the values of the vertices agree. The value of $\mbf{U}_{\leq k}$ on the edges is fixed by the inductive hypothesis \eqref{81ITXUL}. Here, we see that there is agreement by a direct check: writing $\gamma'\coloneqq \paren{\Delta\{i-1,k+1\}\hookrightarrow\Delta[k+1]}^*\gamma$, we apply the definition from \Cref{STG38ZQ} and see $\mbf{U}_{\leq 1}(\gamma',1)(\underline{012}>\underline{02})=\gamma'$.
\end{proof}

\begin{notation*} We write
    \(
        \btu=\btu_{\bullet}\coloneqq\Nerve_{\bullet}\paren{\overline{P_{0,\widehat{k+2}}}},
    \)
    and 
    \[
        \partial\coloneqq\paren{\underline{0,1,k+2}\to\cdots\to\underline{0,k+1,k+2}}\in\btu_{k}.
    \]
\end{notation*}

By \Cref{YGELJPS}, the following concludes the proof of \Cref{YUXW1VO}.

\begin{proposition}\label{VMR7DN8}
    For each     $(\gamma,j)\in\Phat_{k}\subset\EEx_{k+1}\smallsetminus\overline{\EEx_{\leq k}}$, there exists an $\infty$-functor
    \[
        \mbf{U}^{\btu}\colon\btu\to B\OO(n+m)
    \]
    that extends $\mbf{U}_{\leq k}$ and is such that $\mbf{U}^{\btu}(\partial)=d_{k+1}(\gamma)$.
\end{proposition}

\begin{proof}
    All simplices in $\btu$ that have at most one vertex from $\partial$ are already determined by $\mbf{U}_{\leq k}$ and functoriality, as noted at the beginning of this section. Consequently, the non-degenerate $(k+1)$-simplices of $\btu$ are exactly those that possess non-degenerate edges in $\partial$. Setting $\mbf{U}^{\btu}(\partial)=d_{k+1}(\gamma)$, we will exhibit natural fillers for these, generalising the observations in \Cref{0N49067}. 
    Let 
    \begin{equation}\label{KSQPKID}
        \underline{0,\alpha^0,k+2}\geq\cdots\geq\underline{0,\alpha^r,k+2}\geq\underline{0,\beta_{r+1},k+2}\to\cdots\to\underline{0,\beta_{r+s},k+2}
    \end{equation}
    be an element of $\btu_{r+s}$ which we will denote by $X$, where, with a slight abuse of notation,
    \[
        [1,k+1]\supseteq\alpha^{0}\supseteq\cdots\supseteq\alpha^{r}\supseteq\beta_{r+1}\leq\cdots\leq\beta_{r+s}
    \]
    with non-empty sequences $\alpha^{i}=\paren{\alpha^i_j}_{j=1}^{r_i}$ and elements $\beta_i$, $s\geq1$. Observe that if $N_{\alpha^{i}}=1$, then also $N_{\alpha^{i'}}=1$ for all $i'\geq i$. 
    We have $\alpha^{0}_{1}\leq\cdots\leq\alpha^{r}_{1}$ due to the subset inclusions, and therefore 
    \begin{equation}\label{F9JN924}
       \alpha^{0}_{1}\leq\cdots\leq\alpha^{r}_{1}\leq\beta_{n+1}\leq\cdots\leq\beta_{r+s}. 
    \end{equation}
    Let, then, 
    \[
        I\in\{-1,0,\dots,r+s\}
    \]
    be the the smallest index such that $\alpha^{I'}_{1}>j$ or else $\beta_{I'}>j$ for all $I'\geq I+1$. If the former, this implies that $N_{\alpha^{I'}}=1$, and so $\alpha^{I'}_{1}>j$, for all $I'\geq I.$ Now \Cref{3077OTN} and \Cref{OXY4ORA} imply that (the $1$-skeleton of) $\mbf{U}^{\btu}_{(\gamma,j)}(X)$ must be of type 
    \begin{align*}
        &\mfrak{N}_{\alpha^{0}_{N-1}-1}\oplus(\gamma,j)|_{\alpha^{0}_{1}-1}\to \cdots\to \mfrak{N}_{\alpha^{I}_{N-1}-1}\oplus (\gamma,j)|_{\alpha^{I}_{1}-1}\to\\
        &(\gamma,j)|_{\alpha^{I+1}_{1}-1}\to\cdots\to(\gamma,j)|_{\alpha^{r}_{1}-1}\to\\
        &\gamma|_{\beta_{r+1}-1}\to\cdots\to\gamma|_{\beta_{r+s}-1}.
    \end{align*}

    If $I=-1$, then $\alpha^{i}_1-1\geq j$ for all $i\geq 0$ and so every $(\gamma,j)|_{\alpha^{i}_{1}-1}$ is upper. We set
    \begin{equation}\label{VAAZR52}
        \mbf{U}^{\btu}_{(\gamma,j)}(X)\coloneqq \Xi^*(\gamma)
    \end{equation}
    using the map
    \begin{align*}
        \Xi\colon[r+s]&\to[k+1],\\
        i&\mapsto\begin{cases}
            \alpha^{i}_{1}-1, & i\leq r\\
            \beta_{i}-1, & r+1\leq i\leq r+s
        \end{cases}
    \end{align*}
    which is monotone by \eqref{F9JN924}, and observe that it is of the desired type. One such case occurs when $X=\partial$ (with $r=-1$ and $s=k+1$) and reproduces $\mbf{U}^{\btu}_{(\gamma,j)}(\partial)=d_{k+1}(\gamma).$

    If $I\geq0$, then the construction is naturally partitioned into cases depending on $(\gamma,j)$.
    Since $k+1$ is the final vertex, $d_{k+1}\gamma$ is either low or vertical. 
    
    It is low iff $j=k+1$, in which case $I=r+s$. This implies that $\mbf{U}^{\btu}_{(\gamma,k+1)}(X)$ has no terms of type $(\gamma,j)|_*$, and that $\gamma|_{\beta_{t}-1}=\mfrak{N}_{\beta_{t}-1}\oplus\pi(\gamma|_{\beta_{t}-1})$ (as noted in the proof of \Cref{OXY4ORA}). Similarly, we have $(\gamma,k+1)|_{\alpha^{t}_{1}-1}=\pi(\gamma|_{\alpha^{t}_{1}-1})$. In sum, $\mbf{U}^{\btu}_{(\gamma,k+1)}(X)$ is to be of type 
    \begin{align*}
        &\mfrak{N}_{\alpha^{0}_{N-1}-1}\oplus\pi(\gamma|_{\alpha^{0}_{1}-1})\to \cdots\to \mfrak{N}_{\alpha^{r}_{N-1}-1}\oplus \pi(\gamma|_{\alpha^{r}_{1}-1})\to\\
        &\mfrak{N}_{\beta_{r+1}-1}\oplus\pi(\gamma|_{\beta_{r+1}-1})\to\cdots\to\mfrak{N}_{\beta_{r+s}-1}\oplus\pi(\gamma|_{\beta_{r+s}-1}).
    \end{align*}
    Let us write 
    \[
        \delta^{T}\coloneqq
        \begin{cases}
            \alpha^{T}_{N-1}-1, &0\leq T\leq r\\
            \beta_{T}-1, & r+1\leq T\leq r+s.
        \end{cases}
    \]
    We claim that 
    \begin{equation}\label{M535MPW}
        \mbf{U}^{\btu}_{(\gamma,k+1)}(X)\coloneqq\Gamma_W\oplus\Gamma_V
    \end{equation}
    does the job, where, first,
    \begin{align*}
        \Gamma_W\colon\Delta^{r+s}&\to B\OO(m)\\
        \paren{t_0,\dots,t_{r+s}}&\mapsto \mfrak{N}\paren{\sum_{\delta^{T}=0}t_T,\dots,\sum_{\substack{\delta^{T}=k}}t_T}
    \end{align*}
    where an empty sum is understood to give $0$.
    In other words, $t_T$ is a summand of (and only of) entry $\delta^{T}$. Observe that $\Gamma_W$ is well-defined: since the exit index in this case is $j=k+1$, $\mfrak{N}$ is a $k$-path. Immediately from the definition of $I$, we have $\delta^{T}=\alpha^{T}_{N-1}-1\leq j-1=k$ for $0\leq T\leq r$, and similarly $\delta^{T}=\beta_T-1\leq j-1=k$ for $r+1\leq T\leq r+s$, so that the coordinate expression makes sense. Finally, the sum of all entries is equal to $\sum_{T=0}^{r+s} t_T=1$ since, by construction, every $t_T$ appears exactly once.

    Secondly for \eqref{M535MPW}, we set 
    \[
        \Gamma_V\coloneqq \Xi^*\pi(d_{k+1}\gamma)
    \]
    using the map 
    \(
        \Xi\colon[r+s]\to[k]\), \(
        i\mapsto
        \begin{cases}
            \alpha^{i}_{1}-1, & i\leq r\\
            \beta_{i}-1, & r+1\leq i\leq r+s
        \end{cases}
    \)
    as with \eqref{VAAZR52} except that the target is different.
    
    The compatibility of \eqref{M535MPW} with the inductive assumption \eqref{81ITXUL} holds trivially. As for \eqref{PTBY5XA}, we may assume, without loss of generality, that $\alpha^{T}_{1}=i$ for some arbitrary but fixed $i\in\{1,\dots,j\}=\{1,\dots,k+1\}$ for all $T$, and it suffices to consider the restriction of $X$ along $\Delta\{0,\dots,r\}\hookrightarrow \Delta[r+s]$, since for the range from $r+1$ to $r+s$ we would have to assume $\beta_{T}=i$ for all $T$ as well, whence that range is, for our purposes in this case, degenerate. We may now decompose as 
    \[
        X|_{0,\dots,r}=\paren{\underline{\alpha^{0},k+2}\geq\cdots\geq\underline{\alpha^{r},k+2}}\cup \underline{0,i}
    \]
    with the left factor in $\Nerve_{\bullet}\paren{P^{\opp}_{i,k+2}}$. In order to prove compatibility, let  $T\in[r]$ and observe that 
    \[
        \Gamma_W|_{T}(1)=\Gamma_W(0,\dots,1,\dots,0)=\mfrak{N}|_{\alpha^{T}_{N-1}-1}(1)=\Opp\mfrak{N}|_{T}(1),
    \]
    by \Cref{J8KZ1VE}, where in the second term $1$ appears in entry $T$.\footnote{This $1$ is the one of $\Delta^{0}=\{1\}\subset\mbb{R}$.} More generally, given $(t_0,\dots,t_r)\in\Delta^r$, note that $\alpha^{T}_{N-1}\geq\alpha^{T'}_{N-1}$ for $T\leq T'$ as we know from the proof of \Cref{KXIJ7LO}, and so 
    \begin{align*}
        \Gamma_W|_{0,\dots,r}(t_0,\dots,t_r)&=\Gamma_W(t_0,\dots,t_r,0,\dots,0)\\
        &=\mfrak{N}|_{\alpha^{r}_{N-1}-1,\dots,\alpha^{0}_{N-1}-1}\paren{\sum_{\alpha^{T}_{N-1}=\alpha^{r}_{N-1}}t_{T},\dots,\sum_{\alpha^{T}_{N-1}=\alpha^{0}_{N-1}}t_{T}}\\
        &=\Opp\mfrak{N}|_{0,\dots,r}(t_0,\dots,t_r)
    \end{align*}
    directly by the definition of the latter in \Cref{2QKCQ8U} and by \Cref{J8KZ1VE}, since a poset map $\phi\colon[K]\to[L]$ is mapped under geometric realisation to the map 
    \(
        (t_0,\dots,t_K)\mapsto \paren{\sum_{\phi(T)=0}t_{T},\dots,\sum_{\phi(T)=L}t_{T}}.
    \)

    We must also check compatibility with the condition $\mbf{U}^{\btu}(\partial)=d_{k+1}(\gamma)$, but this is straightforward. The relevant range is from $r+1$ to $r+s$, and here we observe that 
    \begin{align*}
        \Gamma_W|_{r+1,\dots,r+s}(t_{r+1},\dots,t_{r+s})&=\Gamma_W(0,\dots,0,t_{r+1},\dots,t_{r+s})\\
        &=\mfrak{N}|_{\beta_{r+1}-1,\dots,\beta_{r+s}-1}(t_{r+1},\dots,t_{r+s})
    \end{align*}
    since $\beta_*$ is non-decreasing. This implies 
    \begin{align*}
        \Gamma_W\oplus\Gamma_V&=\mfrak{N}_{\beta_{r+1}-1,\dots,\beta_{r+s}-1}\oplus\pi(d_{k+1}\gamma)|_{\beta_{r+1}-1,\dots,\beta_{r+s}-1}\\
        &=d_{k+1}\gamma|_{\beta_{r+1}-1,\dots,\beta_{r+s}-1}
    \end{align*}
    since this face is low in this case by assumption.

    If, lastly, $d_{k+1}(\gamma)$ is vertical, or equivalently $j\leq k$ which covers the case $0\leq I<n+m$, then 
    \begin{equation}\label{SVZ4YUC}
        d_{k+1}(\gamma,j)=(d_{k+1}\gamma,j)
    \end{equation}
    since $k+1> j$ implies $\flat_{j,k+1}=j$.
    We first claim that
    \begin{equation}\label{4HYZ2ZB}
        \mbf{U}_{(\gamma,j)}(\underline{0,\alpha^{T},k+2})=
        \begin{cases}
            \mbf{U}_{d_{k+1}(\gamma,j)}(\underline{0,\alpha^{T},k+1}), &\alpha^{T}_{r_{T}}\leq k,\\
            \mbf{U}_{d_{k+1}(\gamma,j)}(\underline{0,\alpha^{T}}), &\alpha^{T}_{r_{T}}=k+1
        \end{cases}
    \end{equation}
    Indeed, recall that by \eqref{0DF723Q} from the proof of \Cref{Q52DBDN} we have 
    \[
        \mbf{U}_{(\gamma,j)}(\underline{i,\ell})=\mbf{U}_{(\gamma,j)}(\underline{i,\ell'})
    \]
    whenever $1\leq i\leq j<\ell,\ell'\leq k+2$. This immediately implies the claim since by \eqref{SVZ4YUC} the exit indices on both sides coincide. Namely, the latter implies
    \[
        N_{\underline{\alpha^{T},k+2}}=
        \begin{cases}
            N_{\underline{\alpha^{T},k+1}}, & \alpha^{T}_{r_{T}}\leq k\\
            N_{\alpha^{T}}, & \alpha^{T}_{r_{T}}=k+1
        \end{cases}
    \]
    since the first index on either side that exceeds $j\leq k$ is not affected by whether the full sequence ends with $k+1$ or $k+2$. Therefore, if $N$ (which we can thus employ unambiguously) is $1$, then, by \Cref{3077OTN}, \eqref{4HYZ2ZB} becomes 
    \begin{equation}\label{FZKSY0I}
        (\gamma,j)|_{\alpha^{T}_{1}-1}=(d_{k+1}(\gamma,j))|_{\alpha^{T}_{1}-1}=(d_{k+1}\gamma,j)|_{\alpha^{T}_{1}-1}
    \end{equation}
    which holds by simpliciality since $\alpha^{T}_{1}-1\leq k$ as $\alpha^{T}\subseteq[1,k+1]$, whence $\alpha^{T}_{1}-1\in\mrm{Im}(\partial_{k+1})$. If $N>1$, then again using \Cref{3077OTN} together with \eqref{FZKSY0I} we see that \eqref{4HYZ2ZB} is tantamount to 
    \[
        \mfrak{N}(\gamma,j)(\underline{\alpha^{T},k+2}_{N-1}-1)=
        \begin{cases}
            \mfrak{N}(d_{k+1}(\gamma,j))(\underline{\alpha^{T},k+1}_{N-1}-1), & \alpha^{T}_{r_{T}}\leq k,\\
            \mfrak{N}(d_{k+1}(\gamma,j))(\alpha^{T}_{N-1}-1), & \alpha^{T}_{r_{T}}=k+1.
        \end{cases}
    \]
    which holds similarly. In the same way we obtain
    \begin{equation}\label{U22OE9L}
        \mbf{U}_{(\gamma,j)}(\underline{0,\beta_{T},k+2})=
        \begin{cases}
            \mbf{U}_{d_{k+1}(\gamma,j)}(\underline{0,\beta_{T},k+1}), &\beta_{T}\leq k\\
            \mbf{U}_{d_{k+1}(\gamma,j)}(\underline{0,k+1}), &\beta_{T}=k+1
        \end{cases}
    \end{equation}
    using \Cref{OXY4ORA}, which is to say, as a special case of the above. 

    Now, the inductive hypothesis provides 
    \[
        \mbf{U}_{\leq k}(d_{k+1}(\gamma,j))\colon\Path[k+1]\to B^{\oplus}\OO
    \]
    and so in particular a map
    \begin{equation}\label{VDAPNI1}
        \mbf{U}_{\leq k}(d_{k+1}(\gamma,j))|_{\Hom(0,k+1)}\colon\Nerve_{\bullet}\paren{P^{\opp}_{0,{k+1}}}\to B\OO(n+m).
    \end{equation}
    Consider the projection 
    \begin{align*}
        \Pi\colon\overline{P^{\opp}_{0,\widehat{k+2}}}&\twoheadrightarrow \overline{P^{\opp}_{0,\widehat{k+1}}},\\
        \underline{0,\alpha,k+2}&\mapsto \begin{cases}
            \underline{0,\alpha,k+1}, &\alpha_{n}\leq k\\
            \underline{0,\alpha}, &\alpha_{n}=k+1
        \end{cases}
    \end{align*}
    with $\alpha=(\alpha_{1}\leq\cdots\leq\alpha_{n})$, $n\geq0$, within $[1,k+1]$, where $n=0$ is understood to give the empty sequence. It is clearly functorial. Using the adjunction between geometric realisation and the singular chains functor one more time, we can write \eqref{VDAPNI1} as a continuous map of type $|\Nerve_{\bullet}(P^{\opp}_{0,k+1})|\to B\OO(n+m)$, and, applying \Cref{5KGP6JE}, obtain the restriction 
    \[
        \left|\Nerve_{\bullet}\paren{\overline{P^{\opp}_{0,\widehat{k+1}}}}\right|\hookrightarrow |\Nerve_{\bullet}(P^{\opp}_{0,k+1})| \to B\OO(n+m).
    \]
    Finally, un-applying the adjunction yields the further-restricted $\infty$-functor
    \[
        \mbf{U}_{\leq k}(d_{k+1}(\gamma,j))\colon \Nerve_{\bullet}\paren{\overline{P^{\opp}_{0,\widehat{k+1}}}}\to B\OO(n+m).
    \]
    We can thus compose and obtain 
    \[
        \Pi^*\mbf{U}_{\leq k}(d_{k+1}(\gamma,j))\colon \btu\twoheadrightarrow \Nerve_{\bullet}\paren{\overline{P^{\opp}_{0,\widehat{k+1}}}}\to B\OO(n+m)
    \]
    and consequently set
    \begin{equation}\label{8RDCLNK}
        \mbf{U}^{\btu}_{(\gamma,j)}(X)\coloneqq\Pi^*\mbf{U}_{\leq k}(d_{k+1}(\gamma,j))(X).
    \end{equation}
    The equalities \eqref{4HYZ2ZB} and \eqref{U22OE9L} state precisely that \eqref{8RDCLNK} is compatible with $\mbf{U}_{\leq k}$. 
    
    On the other hand, we may append the compatibility of \eqref{M535MPW} and $\mbf{U}_{\leq k}$ to the inductive hypothesis. Namely, we assume that \emph{the map induced (by repeated use of the adjunction between geometric realisation and the singular chains functor) by $\mbf{U}_{\leq k}$ itself on $\Nerve_{\bullet}\paren{\overline{P^{\opp}_{0,\widehat{k'+2}}}}$, for all $k'\leq k-1$, is given by \eqref{M535MPW} on any exit path whose $(k'+1)$-face is low.}
    
    This assumption is justified since it holds in the base case $k=1$, as we will now observe. The case of interest is where the exit index is $j=2$, and there we have given the filler of \eqref{T4UNR1L}, which depicts exactly $\overline{P^{\opp}_{0,\widehat{3}}}$, by 
    \(
        \Gamma\oplus s_{0}(\gamma_V)
    \)
    and we see that 
    \[
        \Gamma(t_{0},t_{1},t_{2})=\gamma_W(t_{1},1-t_{1})=\gamma_{W}(t_{1},t_{0}+t_{2})=\Gamma_{W}(t_{0},t_{1},t_{2}).
    \]
    Similarly, on $X=(\underline{0123}\geq\underline{013}\to\underline{023})$ we have that $\Xi\colon [r+s]=[0+2]\to[1]$ maps $0\mapsto \alpha^{0}_{1}-1=0$, $1\mapsto\beta_{1}-1=0$, $2\mapsto\beta_{2}-1=1$, and so $\Xi=\sigma_{0}$. Thus 
    \[
        s_{0}(\gamma_{V})=\Xi^{*}\pi(d_{2}\gamma)=\Gamma_{V}.
    \]

    Consequently, \eqref{8RDCLNK} itself is also compatible with \eqref{M535MPW} by virtue of being compatible with $\mbf{U}_{\leq k}$ by the inductive assumption.
\end{proof}

The proof of \Cref{VMR7DN8} cannot be read off from the examples in \Cref{0N49067} and \Cref{9ZWYRMA} alone. Let us therefore give two final examples that illustrate the novel cases treated in the proof.

\begin{example}
    Let us consider a case where $d_{k+1}\gamma$ is vertical. 
    Suppose $k=2$ and $j=2$, so that $\gamma$ (as visualised within the $3$-cylinder) is of type
    \[
        \begin{tikzcd}[sep=small]
        & K\ar[dr] & \\
        && K'\\
        &W_{1}\oplus V_{1}\ar[uu]\ar[ur] &\\
        W_{0}\oplus V_{0}\ar[uuur]\ar[ur]
        \end{tikzcd}
    \]
    where we omitted the edge $W_{0}\oplus V_{0}\to K'$. The face $d_{k+1}\gamma$ is
    \[
        \begin{tikzcd}[sep=small]
        & W_{1}\oplus V_{1}\ar[dr] & \\
        W_{0}\oplus V_{0}\ar[ur]\ar[rr] && K
        \end{tikzcd}
        .
    \]
    The image of $P^{\opp}_{0,4}$ under $\mbf{U}$, with $\gamma$ painted in blue into $\btd$ within the geometric realisation, is then of type
    \begin{equation*}
        \begin{tikzcd}[sep=small]
        & W_1\oplus V_0 \arrow[dl] \arrow[rr] \arrow[dd,equal] & & W_1\oplus V_0 \arrow[dl] \arrow[dd] \\
        W_1\oplus V_1 \arrow[rr, crossing over] \arrow[dd] & & W_1\oplus V_1 \\
        & W_0\oplus V_0 \arrow[dl] \arrow[rr] & & W_0\oplus V_0\ar[ul,color=blue],\ar[dlll,color=blue] \arrow[dl,color=blue] \\
        K \arrow[rr,color=blue]\ar[from=uurr,crossing over,color=blue,bend left=10] & & K' \arrow[from=uu, crossing over,color=blue]\\
        \end{tikzcd}
        .
    \end{equation*}
    A path $X$ in $\Nerve_{\bullet}\paren{\overline{P^{\opp}_{0,\widehat{4}}}}$ as in \eqref{KSQPKID} that is novel is given, for instance, by $\underline{0124}\geq\underline{024}\to\underline{034}$, whose image under $\mbf{U}^{\btu}$, as can be read off the cube, is to be of type
    \[
        W_{1}\oplus V_{0}\to W_{1}\oplus V_{1}\to K.
    \]
    The filler thereof provided by the proof is exhibited by first factoring $X$ through the projection $\Pi$, which yields $\Pi(X)=\paren{\underline{0123}\to\underline{023}\to\underline{03}}$, which by $\mbf{U}(d_{3}(\gamma,2))$ is mapped to a $2$-simplex of type 
    \[
        \begin{tikzcd}[sep=small]
        & W_{1}\oplus V_{1}\ar[dr] & \\
        W_{1}\oplus V_{0} \ar[ur]\ar[rr] && K
        \end{tikzcd}
    \]
    provided exactly by the left triangle in \eqref{9681FJ4} in the case $j=2=\flat_{2,3}$.
\end{example}

\begin{example}
    In the situation of \Cref{9ZWYRMA}, i.e., with $k=2$ and $j=3$ so that $d_{k+1}\gamma$ is low, a novel path $X$ in  $\Nerve_{\bullet}\paren{\overline{P^{\opp}_{0,\widehat{4}}}}$ as in \eqref{KSQPKID} is given, for instance, by $\underline{01234}\geq\underline{0124}\geq\underline{024}\to\underline{034}$. Its image under $\mbf{U}^{\btu}$ is to be of type 
    \begin{equation}\label{G5EBZJV}
        W_{2}\oplus V_{0}\to W_{1}\oplus V_{0}\to W_{1}\oplus V_{1}\to W_{2}\oplus V_{2}.
    \end{equation}
    We read off
    \[
        \Gamma_W(t_{0},t_{1},t_{2},t_{3})=\gamma_{W}(0,t_{1}+t_{2},t_{0}+t_{3})
    \]
    where the normal path $\gamma_W$ is of type 
    \[
        \begin{tikzcd}[sep=small]
        & W_{1}\ar[dr] & \\
        W_{0}\ar[ur]\ar[rr] && W_{2}
        \end{tikzcd}
        .
    \]
    It is a direct check to see that $\Gamma_W$ fills the normal component of \eqref{G5EBZJV}.
\end{example}

\begin{remark}
    None of the results and constructions above depend on the properties of infinite Grassmannians. For any topological monoid $M$ as in \Cref{ADAUJCW}, assume that its operation $\odot$ is a cofibration, and that $M=\coprod M_i$ over some index set. Then we may consider two `strata' $M_1$, $M_2$, and supposing that the operation thereon factors through another stratum, $\odot\colon M_1\times M_2\to M_{i_{12}}$, obtain the linked space $\left(M_1\leftarrow M_{1}\times M_2\to M_{i_{12}}\right)$. Then the unpacking map gives a fully-faithful $\infty$-functor 
    \[
        \EEx\left(M_1\leftarrow M_{1}\times M_2\to M_{i_{12}}\right)\hookrightarrow \ast/\Nhc(BM)=\mscc{M}^{\hookrightarrow},
    \]
    giving a particularly simple equivalent description of the  full depth-$1$ sub-$\infty$-categories of $\mscc{M}^{\hookrightarrow}$. \Cref{7UKFW0A,XTGZI27,QOSZIFS,IPWVX24} follow.
\end{remark}

\begin{example}\label{QN6S8C6} The Stiefel manifolds $E\OO(n)\to B\OO(n)$ give rise to an $\infty$-category $\mscc{F}^{\hookrightarrow}\coloneqq\ast/\Nhc(E^{\oplus}\OO)$, the \emph{stratified Stiefel manifold} or the \emph{stratified universal frame bundle}, such that for any pair $n,m$ of natural numbers the unpacking map gives a fully faithful embedding 
\[
    \EEx(E\OO(n,m))=\EEx(E\OO(n)\overset{\mrm{pr}}{\twoheadleftarrow}E\OO(n)\times E\OO(m)\overset{\oplus}{\hookrightarrow}E\OO(n+m))\hookrightarrow\mscc{F}^{\hookrightarrow}.
\]
Moreover, the square 
\cdon 
    \EEx(E\OO(n,m))\ar[d,two heads]\ar[r,hook] &\mscc{F}^{\hookrightarrow}\ar[d,two heads]\\
    \EEx(B\OO(n,m))\ar[r,hook] & \Vhook
\cdoff
commutes.
\end{example}

We will conclude with a construction promised in \Cref{4J6WZPX}.

\begin{remark}\label{LP5KAHG}
    Let $\Gamma\colon\Path[k+1]\to B^{\oplus}\OO$ be a $k$-simplex of $\Vhook$, and $\btd$ as in \Cref{HZBV0LB}. Then the rule\footnote{where we use the restriction $\Gamma\colon\Hom_{\Path[k+1]}(0,k+1)=\Nerve_{\bullet}({P^{\opp}_{0,k+1}})\to B\OO^{\infty}_{\amalg}$, take the corresponding continous map $|\Nerve_{\bullet}(P^{\opp}_{0,k+1})|\to B\OO^{\infty}_{\amalg}$, and finally pull it back along $\btd\colon\Delta^{k}\hookrightarrow\left|\Nerve_{\bullet}\right|$}
    \begin{align*}
        \Vhook_{k}&\to (B\OO^{\infty}_{\amalg})_{k}\\
        \Gamma&\mapsto \Gamma(\btd),
    \end{align*}
    restricts on the core $\Vsim\subset\Vhook$ to an inverse to the map $\Psi\colon B\OO^{\infty}_{\amalg}\to\Vhook$ from the proof of \Cref{prop:Vsim=BO}.
    In contrast to the putative rule $\Psi^{-1}$ in \Cref{4J6WZPX}, this is functorial: $\btd=\btd^{k}\colon\Delta^{k}\hookrightarrow|\Nerve_{\bullet}(P^{\opp}_{0,k+1})|$ is the convex hull of the corners $\underline{0,i,k+1}$ and $\underline{0,k+1}$ within the $k$-cube (along \Cref{ECHTYNN}), so evidently 
    \[
        (d_{i+1}^{\Bbox\OO}\Gamma)(\btd^{k-1})=d_{i}\Gamma(\btd^{k})
    \]
    and 
    \[
        (s_{i+1}^{\Bbox\OO}\Gamma)(\btd^{k+1})=s_{i+1}\Gamma(\btd^{k})
    \] hold for $i\in\{0,\dots,k\}$, proving functoriality. 
\end{remark}

\printbibliography

\end{document}